\documentclass{article}
\usepackage{graphicx}
\graphicspath{images/}

\usepackage{blindtext}

\usepackage[utf8]{inputenc}
\usepackage{latexsym, amssymb, latexsym, amsfonts, amsmath, esint, graphicx, caption, subcaption, dsfont, units, fullpage, float, setspace, cases, polynom, multirow , stmaryrd , textcomp, tikz, mathtools, amsthm, epstopdf, epsfig, hyperref}
\usepackage[first=1, last=10]{lcg}

\usepackage{listings,enumerate,geometry,csquotes,url}
\usepackage{etoolbox}
\usepackage{pdfpages,xspace,xparse}
\usepackage{import,pst-plot}
\usepackage[numbers,sort]{natbib}
\usepackage[toc,page]{appendix}

\newcommand*\pFq[2]{{}_{#1}F_{#2}}
\newcommand*\ptFq[2]{{}_{#1}\tilde{F}_{#2}}

\newcommand{\beq}{\begin{equation}}
\newcommand{\eeq}{\end{equation}}
\newcommand{\ba}{\begin{array}}
\newcommand{\ea}{\end{array}}
\newcommand{\bea}{\begin{eqnarray*}}
\newcommand{\eea}{\end{eqnarray*}}
\newcommand{\bc}{\begin{center}}
\newcommand{\ec}{\end{center}}
\newcommand{\bt}{\begin{table}}
\newcommand{\et}{\end{table}}

\newcommand{\la}[1]{\label{#1}}

\newcommand{\no}{\noindent}

\newcommand{\rf}[1]{(\ref{#1})}
\newcommand{\beqno}{\begin{displaymath}}
\newcommand{\eeqno}{\end{displaymath}}

\newcommand{\been}{\begin{enumerate}}
\newcommand{\een}{\end{enumerate}}

\newcommand{\ra}{\rightarrow}

\newcommand{\C}{\mathbb{C}}
\newcommand{\R}{\mathbb{R}}

\newcommand{\Ai}{\mathrm{Ai}}

\newcommand{\Arg}{\mathrm{Arg}}
\renewcommand{\Re}{\mathrm{Re}}
\renewcommand{\Im}{\mathrm{Im}}

\newcommand{\erf}{\mathrm{erf}}

\newcommand{\figcl}[4]{\begin{figure}[tb]\begin{center} \includegraphics[scale=#1]{#2} \end{center} \caption{#3} \label{fig:#4} \end{figure}}

\newcommand{\case}[2]{\left\{ \hspace*{-0.15in} \arraycolsep=0.15in\def\arraystretch{#1}\begin{array}{ll} #2 \end{array}\right.}

\newlength{\myheight}
\newlength{\mylength}

\newcounter{saveeqn}

\newtheorem{theorem}{Theorem}

\def\XXint#1#2#3{{\setbox0=\hbox{$#1{#2#3}{\int}$}
     \vcenter{\hbox{$#2#3$}}\kern-.5\wd0}}

\def \imsize {0.6}

\hypersetup{
    colorlinks=true,
    linkcolor=blue,
    citecolor=red
}
\usepackage{subfiles} 

\title{The analytic extension of solutions to initial-boundary value problems outside their domain of definition}
\author{
Matthew Farkas, Jorge Cisneros, Bernard Deconinck\\
Department of Applied Mathematics\\
University of Washington\\
Seattle, WA 98195-2420
}
\date{\today}

\begin{document}

\maketitle

\begin{abstract}
\sloppypar We examine the analytic extension of solutions of linear, constant-coefficient initial-boundary value problems outside their spatial domain of definition. We use the Unified Transform Method or Method of Fokas, which gives a representation for solutions to half-line and finite-interval initial-boundary value problems as integrals of kernels with explicit spatial and temporal dependence. These solution representations are defined within the spatial domain of the problem. We obtain the extension of these representation formulae via Taylor series outside these spatial domains and find the extension of the initial condition that gives rise to a whole-line initial-value problem solved by the extended solution. In general, the extended initial condition is not differentiable or continuous unless the boundary and initial conditions satisfy compatibility conditions. We analyze dissipative and dispersive problems, and problems with continuous and discrete spatial variables.
\end{abstract}



\section{Introduction}

We examine the analytic continuation, for $x\in \mathbb{R}$, of solutions of linear, constant-coefficient initial-boundary value problems (IBVPs)  outside their spatial domain of definition. This is a classical question: for the heat equation on $x>0$ with Dirichlet boundary data, 
\begin{subequations} \la{heatdirichletibvp}
\begin{align}
u_t&=u_{xx}, &&\hspace*{-1.0in} x>0, \, t>0,\hspace*{1.0in}\\
u(x,0)&=u_0(x), &&\hspace*{-1.0in} x>0,\\ \label{fig:heatdirichletbc}
u(0,t)&=f_0(t), &&\hspace*{-1.0in} t>0,
\end{align}
\end{subequations}
\no the solution can be written using the Fourier sine transform, see \rf{classicalheat}. This leads immediately to an odd extension of the solution for $x<0$, which cannot provide an analytic extension of the solution unless, for starters, the Dirichlet data $f_0(t)\equiv 0$. 
To illustrate our goal, consider the Fourier cosine series of the function $f(x) = (x-1/2)^2$, 
\beq \label{eqn:FS}
f_\text{FS}(x) = \frac{1}{12} + \sum_{n=1}^\infty \frac{\cos(2n\pi x)}{n^2\pi^2},
\eeq
\no shown in Figure~\ref{fig:FS}. While $f_\text{FS}(x)$ converges uniformly to $f(x)$ in $[0,1]$, outside the interval, it converges to the even, periodic extension. There is no clear way to analytically extend the Fourier cosine series in $(0,1)$ to obtain $(x-1/2)^2$ outside of $(0,1)$, since
using \rf{eqn:FS} the higher-order derivatives are not defined. 

\figcl{0.7}{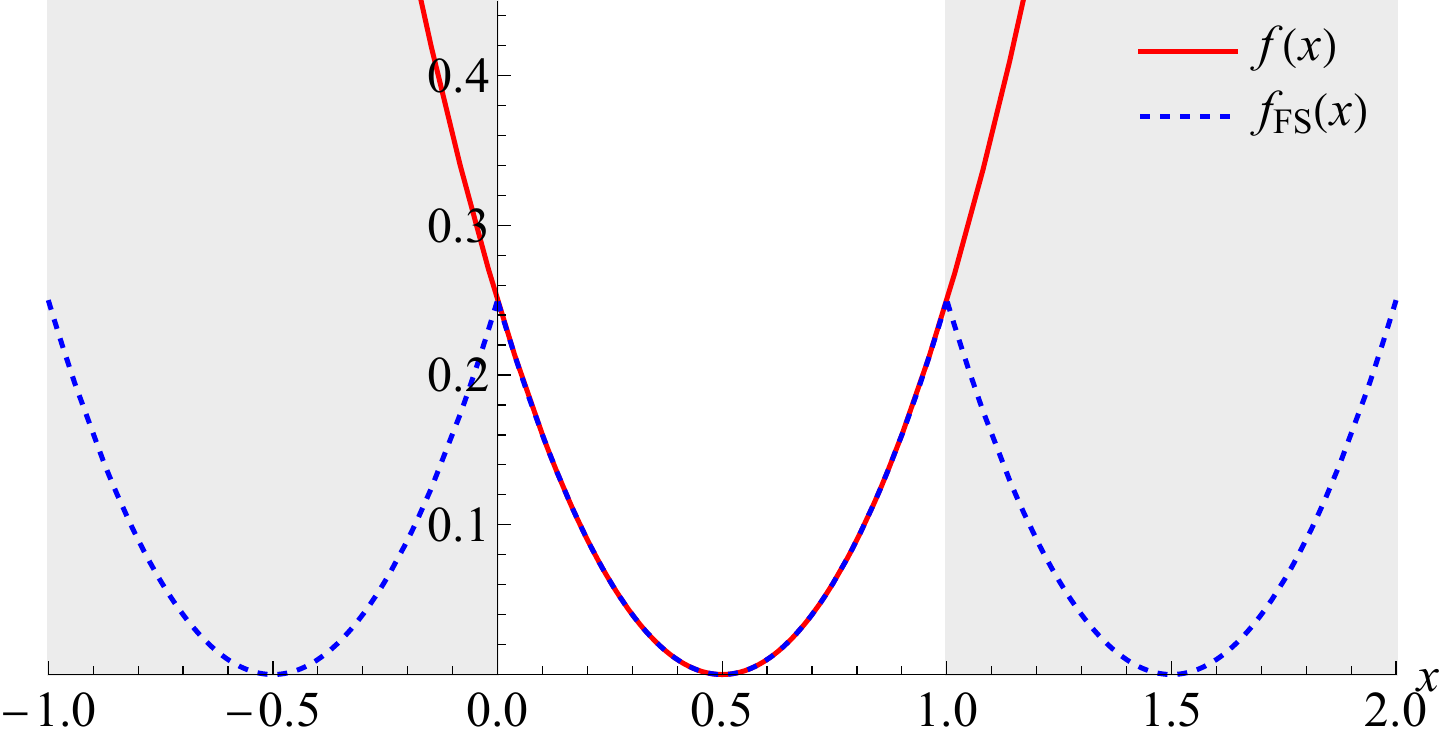}{The function $f(x) = (x-1/2)^2$ shown in red and its Fourier cosine series $f_\text{FS}(x)$ \rf{eqn:FS}, shown in dashed blue.}{FS}

In the rest of this introduction we use the IBVP \rf{heatdirichletibvp} to fix the notation, but we examine problems far more general in this manuscript, including problems of higher order and problems with discrete spatial variables. For those situations where the question of an analytic extension is reasonable ({\em e.g.}, for \rf{heatdirichletibvp}, $u_0(x)$ and $f_0(t)$ are analytic for $x>0$ and have sufficient decay, see Theorem~\ref{thm:IC}), we derive explicit representations for the analytic extension $u_\text{ac}(x,t)$, $x\in \mathbb{R}$, of the solutions. Here $u_\text{ac}(x,t)\equiv u(x,t)$ for $x>0$ in \rf{heatdirichletibvp}. Some of these extensions are obtained through expressions relating $u(x,t)$, $x>0$ and $u_\text{ac}(x,t)$, $x<0$. Others are fully explicit in that they give the analytic extension $u_\text{ac}(x,t)$ directly in terms of the given initial and boundary data. 

Our approach uses the Unified Transform Method (UTM) or Method of Fokas \cite{JC_bernard_fokas, JC_fokas_book}, as this method allows for the solution of problems with continuous and discrete spatial dependence, of arbitrary order. Further, it results in solution expressions with well-defined derivatives 
in their spatial domain. Even using the UTM, there is confusion about the solution outside its domain of definition, see for instance \cite[Section~115]{JC_zwillinger}. During the first steps of the method, a Fourier transform is used where the solution is assumed to be zero outside its domain of definition. However, when the expression for the representation of the solution is obtained, the solution no longer satisfies this assumption. This is a consequence of the elimination (using Jordan's lemma) of integral contributions in the solution expression that are identically zero in the IBVP's domain of definition, but not so outside of it. A secondary aim of this manuscript is to clarify this confusion. 

We have several reasons for wanting to extend the solution to outside its original domain of definition. The first one arises from numerical analysis: spatial finite-difference methods remain among the most popular methods for numerically solving partial differential equations. However, the application of finite-difference stencils near boundaries often requires the value of the solution at so-called ghost points, grid points outside of the physical  domain of definition ({\em e.g., }\cite{JC_randy}). The extended solutions calculated here provide an answer to this, as one can simply evaluate $u_\text{ac}(x,t)$ at the desired negative values of $x$. Second, when testing numerical or other methods for IBVPs, it is common to start from a whole-line problem (where exact solutions may be more readily available) and restrict it to a smaller domain, using as boundary conditions the values of the whole-line solution on the boundaries. We can reverse this: using our approach, we can ask what the whole-line problem is whose restriction is the solution of the original IBVP. An example of a physical extension question one could ask is the following: suppose we have a bi-infinite rod, whose temperature at $t=0$ is known for $x>0$, and whose temperature at $x=0$ is measured for all $t>0$. What is the temperature for $x<0$? Indeed, for the heat equation, our formulas have been derived before, see \cite{burggraf}.

It should be noted that the whole-line problem whose restriction is the original problem may not always be of physical interest. For instance, extending from $x>0$ to $x\in \mathbb{R}$, it may be the case that the extended solution grows as $x\rightarrow -\infty$ or even blows up at a finite $x_0<0$, see Figure~\ref{fig:FS}. Further, if the original initial and boundary conditions are incompatible, $\lim_{t\rightarrow 0^+}u_\text{ac}(x,t)$ will differ from the analytic continuation of the original initial condition. Indeed, an analytic continuation of the original initial condition can only lead to a whole-line problem whose restriction would have compatible initial and boundary data. 

Our approach may be thought of as a generalization of the Method of Images (see \cite{JC_haberman, JC_jackson}, for instance). The Method of Images uses a whole-line (or whole-plane/space, for multi-dimensional settings) problem that reduces to the given half-line problem with the given boundary conditions. For instance, the method of images for a homogeneous, Dirichlet (or Neumann) half-line problem uses a whole-line ``image'' problem with odd (or even) extension of the initial condition.  For this paper's boundary-to-initial maps, we wish to construct a whole-line problem that restricts to the given half-line problem with the corresponding {\em nonhomogeneous} boundary conditions.
%
It is perhaps especially surprising to see such formulas for third-order problems or other problems that do not have the symmetry $x\rightarrow -x$. 

As is often the case in function theory, it is important to distinguish between the abstract concept of a solution and its (different) explicit representations. However, in context it is often clear what is meant, and we may blur the line between solutions, extensions, and their representations frequently, to avoid introducing extra notation. 

In Section~\ref{sec:heat}, we consider extensions of the solutions of the half-line Dirichlet and Neumann IBVPs for the heat equation, as well as an example of a finite-interval IBVP. The advected heat equation is treated in Section~\ref{sec:advheat}. The half-line Dirichlet IBVPs for the linear KdV equations are treated in Sections~\ref{sec:kdv1} (one boundary condition) and \ref{sec:kdv2} (two boundary conditions), respectively. 
Sections~\ref{sec:sdadv}, \ref{sec:sdheat} and the Appendix deal with spatially discrete IBVPs for the advection equation and the heat equation. It should be noted that the analytic continuation of solutions of spatially discrete IBVPs can be done in many ways, as the sole requirement is that the analytic continuation interpolates the solution at the fixed grid points, with a well-defined continuous limit as the grid spacing vanishes. Our approach is to consider functions as analytically depending on a discrete variable, a popular approach in this context, see \cite{JC_Joshi}. We conclude with a summary of the extension formulae obtained, collected in one convenient location.

\section{The heat equation} \label{sec:heat}

For problems with a continuous spatial variable, we may start by considering second-order problems. Indeed, for IBVPs involving the first-order transport equation, 
\begin{subequations}
\begin{align}
u_t+ c u_x&=0, &&\hspace*{-1.0in} x>0,\, t>0,\hspace*{1.0in}\\
u(x,0)&=u_0(x), &&\hspace*{-1.0in} x>0,\\
u(0,t)&=f_0(t), &&\hspace*{-1.0in} t>0,
\end{align}
\end{subequations}
\noindent (the boundary condition in the last line is omitted if $c<0$), the solution is analytically continued trivially using the d'Alembert form \cite{JC_haberman} of the solution, 
\begin{equation} \label{eqn:soln_adv}
    u(x,t)=\left\{
    \begin{array}{rcl}
    u_0(x-ct), && t<x/c, \, x>0,\\
    f_0(t-x/c), && t>x/c, \, x>0,
    \end{array}
    \right.
\end{equation}
for $c>0$, and 
\begin{equation}
    u(x,t)=u_0(x-ct), ~~~x>0, \, t>0, 
\end{equation}
for $c<0$. If $f_0(t)$ ($u_0(x)$) is analytic for $c>0$ ($c<0$), then the solution is trivially extended for negative values of $x$. If these functions are not analytic, then no analytic extension exists. 

\subsection{Dirichlet boundary conditions}\la{sec:DBC}


\def \contour {{\partial\Omega}}

Consider the heat equation on the half line with Dirichlet boundary conditions \rf{heatdirichletibvp}. Using Fokas' Unified Transform Method (UTM) \cite{JC_fokas_book}, its solution is written as 
\begin{equation}\label{utmsol1}
    u(x,t) = I_0(x,t) + I_{f_0}(x,t),
\end{equation}
\no with
\begin{align}\la{I0heat1}
I_0(x,t) &= \frac{1}{2\pi} \int_{-\infty}^\infty e^{ikx-k^2t} \hat u_0(k) \,dk - \frac{1}{2\pi} \int_\contour e^{ikx - k^2t}\hat u_0(-k) \,dk,\\\la{Ifheat1}
I_{f_0}(x,t) &=  \frac{1}{i\pi} \int_\contour ke^{ikx - k^2t} F_0(k^2,t) \,dk,
\end{align}
\no and 
\beq
\hat u_0(k) = \int_0^\infty e^{-iky} u_0(y) \,dy, \qquad
F_0(k^2,t) = \int_0^t e^{k^2s} f_0(s) \,ds.
\eeq 
\no The region $\Omega = \{
k\in \C: |k|>r, $ and $\pi/4 < \Arg(k) < 3\pi/4 \}$ for some $r>0$ is shown in Figure~\ref{fig:Omega}.

\figcl{\imsize}{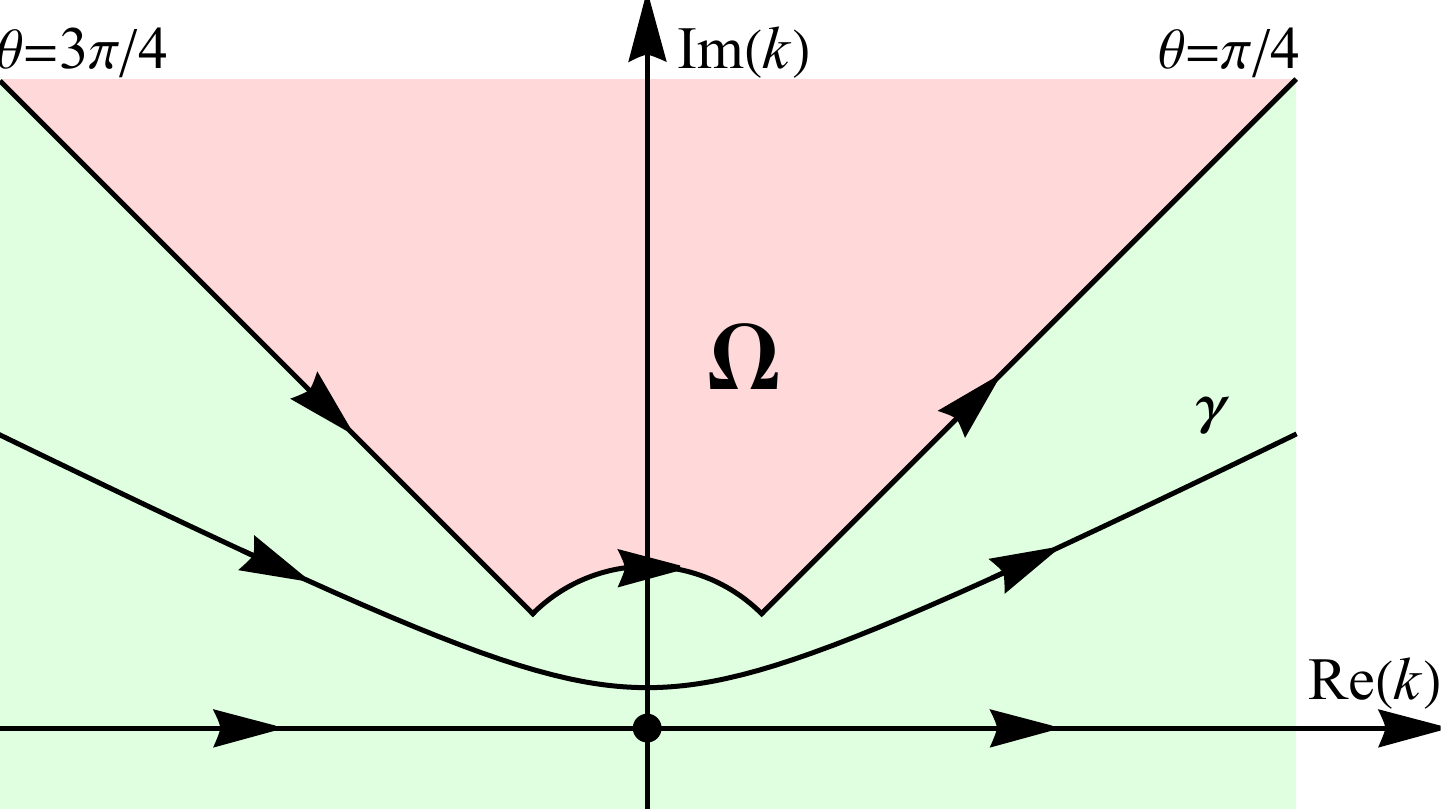}{The region $\Omega 
$ for the heat equation on the half-line.}{Omega}

A simple contour deformation followed by a combination of terms allows for the rewriting of the solution using the classical Fourier sine transform
\begin{equation}\la{classicalheat}
u(x,t)= \frac{2}{\pi} \int_0^\infty dk \, e^{-k^2t} \sin(kx) \left[\int_0^\infty u_0(y)\sin(ky)\,dy + k \int_0^t e^{k^2s}f_0(s)\,ds\right],
\end{equation}
\no from which it is immediately clear that $u(-x,t)=-u(x,t)$. Thus, an analytic continuation of the solution for $x<0$ is not obtained from considering negative values of $x$ in \rf{classicalheat}, unless $f_0(t)\equiv 0$. We pursue a different approach. 




The part $I_0(x,t)$ of the solution containing the initial condition, is entire in $x$ by Theorem~\ref{thm:IC}.

\begin{theorem} \label{thm:IC}

If $u_0\in L^1(\mathbb R^+)$, then $I_0(x,t)$ is entire in $x$ for $t>0$. If $f_0\in L^1(0,T)$ for some $T>0$, then $I_{f_0}(x,t)$ is analytic for $|\Im(x)|<\Re(x)$ for $0\leq t \leq T$. 

\end{theorem}

\begin{proof} 
Let $\Gamma$ be any closed, piecewise smooth contour in the complex $x$-plane. Then, 
\[
\oint_\Gamma dx \int_{-\infty}^\infty e^{ikx-k^2t} \hat u_0(k) \, dk = \int_{-\infty}^\infty dk \, e^{-k^2t}\hat u_0(k) \oint_\Gamma e^{ikx} \, dx = 0,
\]
\no and, since $\hat u_0(-k)$ is analytic for $\Im(k)>0$,
\[ 
\oint_\Gamma dx \int_\contour e^{ikx-k^2t} \hat u_0(-k) \, dk = \oint_\Gamma dx \int_{-\infty}^\infty e^{ikx-k^2t} \hat u_0(-k) \, dk = \int_{-\infty}^\infty dk \, e^{-k^2t}\hat u_0(-k) \oint_\Gamma e^{ikx} \, dx = 0, 
\]
\no by Cauchy's theorem \cite{af}. The order of integration may be switched using Fubini's theorem \cite{folland}, since
$$ \oint_\Gamma |dx| \, \int_{-\infty}^\infty  \left|e^{ikx-k^2t}\hat u_0(\pm k)\right| \, dk \leq \ell(\Gamma)\,\|u_0\|_1 \sqrt{\frac{\pi}{t}} \max_{x\in\Gamma} e^{|x|^2/4t} < \infty, $$
\no where 
%
%
$\ell(\Gamma)$ is the arclength of $\Gamma$. By Morera's theorem \cite{af}, $I_0(x,t)$ is entire in $x$. 

We proceed similarly for $I_{f_0}(x,t)$. We have
\[ 
\oint_\Gamma dx \int_\contour ke^{ikx - k^2t} F_0(k^2,t) \,dk =  \int_\contour  dk \, ke^{- k^2t} F_0(k^2,t) \oint_\Gamma e^{ikx} \,dx = 0,
\]
\no so that $I_{f_0}(x,t)$ is analytic in $x$ wherever Fubini's theorem can be applied to switch the order of integration. Since for $k\in \partial \Omega$, 
\[
\left|e^{-k^2t}F_0(k^2,t)\right| = \left| \int_0^t e^{-k^2(t-s)}f_0(s) \, ds \right| \leq \int_0^T |f_0(s)| \, ds = \|f_0\|_1,
\]
\no it follows that, with $\Gamma$ in the region $|\Im(x)|<\Re(x)$,  
\begin{align*}
\oint_\Gamma |dx| \int_\contour \left|ke^{ikx - k^2t} F_0(k^2,t)\right| \, |dk| &\leq \ell(\Gamma) \|f_0\|_1 \max_{x\in\Gamma} \int_\contour \left|ke^{ikx} \right| \, |dk| \\
&\leq \ell(\Gamma) \|f_0\|_1 \max_{x\in\Gamma} \frac{4|x|^2}{(\Re(x)^2-\Im(x)^2)^2} < \infty,
\end{align*}
\no so that $I_{f_0}(x,t)$ is an analytic function of $x$ for $|\Im(x)|<\Re(x)$.
\end{proof}

Returning to \rf{utmsol1}, for $I_{f_0}(x,t)$ we can switch the order of integration and integrate over $k$ to find
\beq\la{Ifheat2}
I_{f_0}(x,t) = \frac{1}{i\pi} \int_0^t ds \,  f_0(s)\int_\contour ke^{ikx - k^2(t-s)}\,dk   = \frac{x}{2\sqrt{\pi}} \int_0^t \frac{f_0(s)}{(t-s)^\frac32} e^{-\frac{x^2}{4(t-s)}}\,ds,
\eeq
\no which is analytic for $x>0$ by Theorem~\ref{thm:IC}, but is discontinuous at $x=0$, unless $f_0(t)\equiv 0$, which is confirmed below, see \rf{heat_cont_full}. In fact \cite{JC_fokas_book}, 
\beq
\lim_{x\ra 0^+} I_{f_0}(x,t)=f_0(t), ~~~~I_{f_0}(0,t)=0, ~~~~\lim_{x\ra 0^-} I_{f_0}(x,t)=-f_0(t), 
\eeq
as is easily seen from \rf{classicalheat} or \rf{Ifheat2}. 

We extend $I_{f_0}(x,t)$ for $x>0$ to an entire function by constructing its Taylor series about $x=0$. If we do so,
\beq
I_{f_0}(x,t) = \sum_{n=0}^\infty a_{2n}(t) x^{2n} + \sum_{n=0}^\infty a_{2n+1}(t) x^{2n+1}, 
\eeq 
\no which is valid for $x>0$, but trivially extended to negative $x$. With $x\rightarrow -x$, we find
\beq
I_{f_0}(-x,t) = \sum_{n=0}^\infty a_{2n}(t) x^{2n} - \sum_{n=0}^\infty a_{2n+1}(t) x^{2n+1},
\eeq 
\no so that
\beq
I_{f_0}(x,t) = 2\sum_{n=0}^\infty a_{2n}(t) x^{2n} - I_{f_0}(-x,t), 
\eeq    
\no which relates the values for $x<0$ to the values for $x>0$. We have chosen to extend $I_{f_0}(x,t)$ to negative values of $x$ using the even terms in the Taylor series because $I_{f_0}(x,t)$ in \rf{Ifheat2} is odd and the Dirichlet condition immediately gives $a_0(t) = f_0(t)$.

We cannot find the Taylor series for $I_{f_0}(x,t)$ by finding the Taylor series for the kernel in \rf{Ifheat2}, since the resulting integrals diverge. It is possible to transform the $s$ variable and find the Taylor series, but it is easier and more generalizable to use the unevaluated contour integrals (since their evaluation is not possible, in general):
\beq 
\frac{\partial^{2n}I_{f_0}}{\partial x^{2n}}\bigg|_{x=0} = \frac{(-1)^n}{i\pi} \int_0^t ds \, f_0(s) \int_{\contour} k^{2n+1} e^{-k^2(t-s)}\,dk.
\eeq 
\no Using \cite{JC_fokas_book},
\beq 
\frac{\partial^n}{\partial s^n}\delta(s-t) = \frac{1}{i\pi} \int_{\contour} k^{2n+1} e^{-k^2(t-s)} dk,
\eeq 
we obtain 
\beq \label{eqn:nonrigorouseven}
\frac{\partial^{2n}I_{f_0}}{\partial x^{2n}}\bigg|_{x=0} = (-1)^n\int_0^t f_0(s) \delta^{(n)}(s-t)\,ds 
= \int_0^t f_0^{(n)}(s)  \delta(s-t)\,ds = f_0^{(n)}(t). 
\eeq 
\no We obtain this same result rigorously using small $x$ asymptotics \cite{benderorszag}, see \rf{eqn:rigorouseven}, but it is convenient to use the delta function when possible. 

Then, for $x<0$, 
\beq
I_{f_0}(x,t) 
= 2\sum_{n=0}^\infty \frac{x^{2n}}{(2n)!}f_0^{(n)}(t) - I_{f_0}(-x,t) = \tilde f_0(x,t)-I_{f_0}(-x,t),
\eeq
\no where we have defined 
\beq 
\label{f0tilde}
\tilde f_0(x,t) = 2 \sum_{n=0}^\infty \frac{x^{2n}}{(2n)!} f_0^{(n)}(t).
\eeq 
\no This equation has been derived by Burggrag, see \cite{burggraf}. It follows that 
\beq
I_{f_0}^{\text{ext}}(x,t) = \case{1.2}{I_{f_0}(x,t), & x> 0, \\ f_0(t), &x=0, \\ \tilde f_0(x,t) - I_{f_0}(-x,t), & x< 0,}
\label{heat_cont_full}
\eeq 
\no is an analytic function for $x\in \R$ (see Theorem~\ref{thm2} below). This illustrates that for $I_{f_0}^{\text{ext}}(x,t)$ to equal the odd function $I_{f_0}(x,t)$ in \rf{Ifheat2}, we need $f_0(t)\equiv 0$. 

In order to obtain all coefficients of the Taylor series for $I_{f_0}(x,t)$, we start by deforming $\Omega$ down to $\gamma$, see Figure~\ref{fig:Omega}. Using repeated integration by parts, 
\begin{align}\nonumber 
I_{f_0}(x,t) &= \frac{1}{i\pi} \int_{\gamma} dk \, ke^{ikx-k^2t} \int_0^t e^{k^2s} f_0(s)\,ds\\\nonumber
&= \frac{1}{i\pi} \int_{\gamma} dk \, ke^{ikx-k^2t} \left[ \sum_{m=1}^n\frac{ (-1)^m f_0^{(m-1)}(0)}{k^{2m}} + \frac{(-1)^n}{k^{2n}} \int_0^t  e^{k^2s}f_0^{(n)}(s)\,ds \right] \\
\label{If0IBP}
&= \sum_{m=1}^n \phi_{m}(x,t)f_0^{(m-1)}(0) + \int_0^t f_0^{(n)}(s)\phi_n(x,t-s)\,ds,
\end{align}
\no where the integral around $\gamma$ of the terms involving $f_0^{(m-1)}(t) e^{k^2t}$ is zero. We define 
\beq  \label{eqn:hdbc_phi}
\phi_m(x,t) = \frac{(-1)^m}{i\pi} \int_{\gamma} \frac{k e^{ikx-k^2t}}{k^{2m}}\,dk.
\eeq 
\no Switching the order of integration above is allowed since, assuming analyticity of $f_0(t)$, for $x>0$, 
\beq
\int_0^t ds \, \left|f_0^{(n)}(s)\right| \int_{\gamma} \left|\frac{ke^{ikx-k^2(t-s)}}{k^{2m}}\right| \,|dk| 
< \infty.  
\eeq 



\no Differentiating \rf{eqn:hdbc_phi} $(2n-q)$-times with respect to $x$ ($q=0,1$), setting $x=0$, and evaluating the integral, we have
\beq
\phi_m^{(2n-q)}(0,t) = \frac{(-1)^m}{i\pi} \int_{\gamma}  \frac{(ik)^{2n-q}ke^{-k^2t}}{k^{2m}}\,dk
= -\delta_{q,1}\frac{(-1)^{n-m}\Gamma\left(n-m+\frac12\right)}{\pi t^{n-m+\frac12}},
\eeq
\no where $\delta_{q,1}=1$ if $q=1$ and $0$ otherwise.
%
%
We can switch the order of differentiation, integration, and taking limits because of absolute integrability and the boundedness in $x$.  This is why we deform to $\gamma$. Therefore,  

\beq\la{kernel}
\frac{\partial^{2n-q}I_{f_0}}{\partial x^{2n-q}}\bigg|_{x=0} = \sum_{m=1}^n \phi_m^{(2n-q)}(x,t) + \int_0^t 
f_0^{(n)}(s)\phi_n^{(2n-q)}(x,t-s)\,ds.
\eeq 

It follows that for the Taylor coefficients of the even powers, we have
\beq
\frac{\partial^{2n}I_{f_0}}{\partial x^{2n}} \bigg|_{x=0} 
= \lim_{x\to0^+} \int_0^t ds \, f_0^{(n)}(s) \frac{1}{i\pi} \int_{\partial \Omega} k e^{ikx-k^2t}\,dk= \lim_{x\to0^+} \frac{x}{2\sqrt{\pi}}\int_0^t \frac{f_0^{(n)}(s)}{(t-s)^{\frac32}} e^{-\frac{x^2}{4(t-s)}}\,ds.  
\eeq
\no Note that the non-integral terms from \rf{kernel} vanish. Substituting $z=x^2/(4(t-s))$, 
\beq \label{eqn:rigorouseven}
\frac{\partial^{2n}I_{f_0}}{\partial x^{2n}} \bigg|_{x=0}= \lim_{x\to0^+} \frac{1}{\sqrt{\pi}}\int_{\frac{x^2}{4t}}^\infty  f_0^{(n)}\left(t-\frac{x^2}{4z^2}\right) \frac{e^{-z}}{\sqrt{z}}\,dz = \frac{f_0^{(n)}(t)}{\sqrt{\pi}} \int_0^\infty \frac{e^{-z}}{\sqrt{z}}\,dz = f_0^{(n)}(t). 
\eeq
\no We can take the limit inside the integral by the dominated convergence theorem \cite{folland}. This more rigorous derivation confirms the result  \rf{eqn:nonrigorouseven} obtained above using delta functions. Due to the simplicity of the delta function approach, we prefer it in what follows. 

Repeating this line of thought for the coefficients of the odd powers in the Taylor series, we get
\beq 
a_{2n-1}(t) = -\frac{\Gamma\left(\frac12\right)}{\pi(2n-1)!} \left[\sum_{m=1}^n \frac{(-1)^{n-m}\Gamma\left(n-m+\frac12\right)}{t^{n-m+\frac12}\Gamma\left(\frac12\right)} f_0^{(m-1)}(0) + \int_0^t   \frac{f_0^{(n)}(s)}{(t-s)^{\frac12}}\,ds \right].
\eeq 
\no These results may be combined to write 
\beq\la{seriesrep}
I_{f_0}(x,t) = \sum_{n=0}^\infty \frac{(-1)^n x^n}{n!}f_0^{\left(\frac{n}{2}\right)}(t),
\eeq 
where the coefficients are interpreted as Riemann--Liouville fractional derivatives \cite{fractional_calculus}. For non-monomial dispersion relations $W(k)$, a concise notation using fractional derivatives is not possible. This representation of $I_{f_0}(x,t)$ is trivially analytically continued for $x\in \C$. Thus the right-hand side of \rf{seriesrep} is a representation for $I_{f_0}^\text{ext}(x,t)$ for all $x\in \C$, which we prove below in Theorem~\ref{thm2}, provided $f_0(t)$ is analytic in a neighborhood of the positive real $t$-axis. 

If an extension is required close to the boundary only (for instance to set up a numerical scheme), the series representation \rf{seriesrep} may be more convenient than the more global extension provided by \rf{heat_cont_full}. In this section, we have been careful to distinguish between a function and its representation in a part of the complex plane. For the sake of brevity, less care is used below. We expect the distinction to be clear in context. 

\begin{theorem} \label{thm2}
Define $\mathcal D=\{t\in \C: \text{dist}(t,[0,T])\leq r\}$, for some $r>0$, $T>0$, a domain in the complex $t$-plane containing the interval $[0,T]$.  If $f_0(t)$ is analytic in $\mathcal D$, then the series representation \rf{seriesrep} is entire in $x$ for $0 \leq t \leq T$.
\end{theorem}

\begin{proof}
Since $f_0(t)$ is analytic in $\mathcal D$, then for all  $0\leq s\leq t\leq T$, $f_0(\tau)$ is analytic in $|\tau-s|\leq r$. By Cauchy's integral formula \cite{af}, we have
\[
\left|f_0^{(n)}(s)\right| \leq \frac{n!}{r^n} \sup_{|\tau-s|=r}|f_0(\tau)| \leq \frac{n!}{r^n} \|f_0\|_{\infty},
\]
\no where $\|f_0\|_\infty = \sup_{t\in \mathcal D} |f_0(t)|$,
%
%
\no which exists, since $f_0(t)$ is analytic in $\mathcal D$. Then, 
$$
|a_{2n}(t)| = \frac{|f_0^{(n)}(t)|}{(2n)!} \leq \frac{n!}{r^n(2n)!}\|f_0\|_\infty, 
$$
\no so that the series of the even terms converges absolutely for all $x$. We also have
\begin{align*}  
|a_{2n-1}(t)| &= \frac{1}{\sqrt{\pi}(2n-1)!}\left|\sum_{m=1}^{n} \frac{(-1)^{n-m}\Gamma\left(n-m+\frac12\right)}{\sqrt{\pi}t^{n-m+\frac12}}f_0^{(m-1)}(0) + \int_0^t \frac{f_0^{(n)}(s)}{\sqrt{t-s}} \, ds\right| \\
&\leq \frac{\|f_0\|_\infty}{\sqrt{\pi}(2n-1)!}\left[\sum_{m=1}^{n} \frac{\Gamma\left(n-m+\frac12\right)(m-1)!}{\sqrt{\pi} t^{n-m+\frac12}r^{m-1}} + \frac{n!}{r^n} \int_0^t  \frac{ds}{\sqrt{t-s}}\right] \\
&\leq \frac{2n!\|f_0\|_\infty\sqrt{t}}{\sqrt{\pi}(2n-1)!r^n}\left[1+\sum_{m=1}^{n} \frac{\Gamma\left(n-m+\frac12\right)(m-1)!}{2\sqrt{\pi} n!} \left(\frac{r}{t}\right)^{n-m+1} \right], 
\end{align*}
\no so that for $t>r$, 
\[
|a_{2n-1}(t)| \leq \frac{3n!\|f_0\|_t\sqrt{t}}{\sqrt{\pi}(2n-1)!r^n}, 
\]
\no while for $t<r$,
\[
|a_{2n-1}(t)| \leq \frac{3n!\|f_0\|_t}{\sqrt{\pi}(2n-1)!t^{n-\frac12}},
\]
\no where we used 
\[
\sum_{m=1}^{n} \frac{\Gamma(n-m+\frac12)(m-1)!}{2\sqrt{\pi}n!} \leq \frac12,
\]
\no which is proved by showing the left-hand side is a decreasing function of $n$. Thus, the series of the odd terms also converges absolutely for all $x\in \C$ and for any $0 \leq t \leq T$.
\end{proof}

\subsubsection{Boundary-to-Initial Map}

Consider \rf{f0tilde}. 
%
%
\no It follows that $\tilde f_{0t}=\tilde f_{0xx}$, so that  
\beq\la{heatacw} 
w(x,t) =u_\text{ac}(x,t) = \case{1.2}{ u(x,t), & x\geq 0, \\ \tilde f_0(x,t)-u(-x,t), & x<0,}
\eeq 

\no is an analytic solution to the heat equation on the whole line with initial condition
\beq \la{heatacw0}
w_0(x) = \case{1.2}{ u_0(x), & x\geq 0, \\ \tilde f_0(x,0) - u_0(-x), & x<0. }
\eeq 
\no This initial condition is analytic when the compatibility conditions $u_0^{(2n)}(0) = f_0^{(n)}(0)$, for $n\geq 0$ are satisfied \cite{trogdon}. If $u_0(x)$ is analytic and these compatibility conditions are satisfied, then $w_0(x)$ equals the analytic extension of $u_0(x)$ for negative values of $x$. Note that for the homogeneous boundary condition $f_0(t)=0$, we recover the method of images' odd extension of the initial condition.

It is noteworthy that $w_0(x)$ may be unbounded for $x<0$. For example, if $f_0(t) = te^{-t}$, then $f_0^{(n)}(0) = -(-1)^{n} n$, and
\beq
\tilde f_0(x,0) = -\sum_{n=1}^\infty \frac{(-1)^n x^{2n}}{(2n-1)!} =x\sin(x).
\eeq 
\no Thus, the corresponding whole-line problem may not be of physical interest. Nevertheless, the expressions obtained can be used to estimate the solution outside but near the physical domain to be used in numerical schemes requiring information at so-called ghost points.

\subsubsection{Examples}

We demonstrate our results using two examples. Our first example starts from a whole-line solution, 
\beq\la{wholeline1}
u_{\R}(x,t) = \frac{e^{-\frac{(x-1)^2}{4t+1}}}{\sqrt{4t+1}},
\eeq
\no from which we construct the initial and boundary conditions for a half-line problem: $u_0(x)=u_{\R}(x,0)$ for $x>0$, and $f_0(t)=u_{\R}(0,t)$. Next, we use our analytic continuation result $u_{\text{ac}}(x,t)$ to reconstruct the solution $u_{\R}(x,t)$. The results are shown in Figure~\ref{fig:DHL_u1}.

\begin{figure}[tb]
\begin{center}
\def \sc {0.475}
\begin{tabular}{cc}
\includegraphics[scale=\sc]{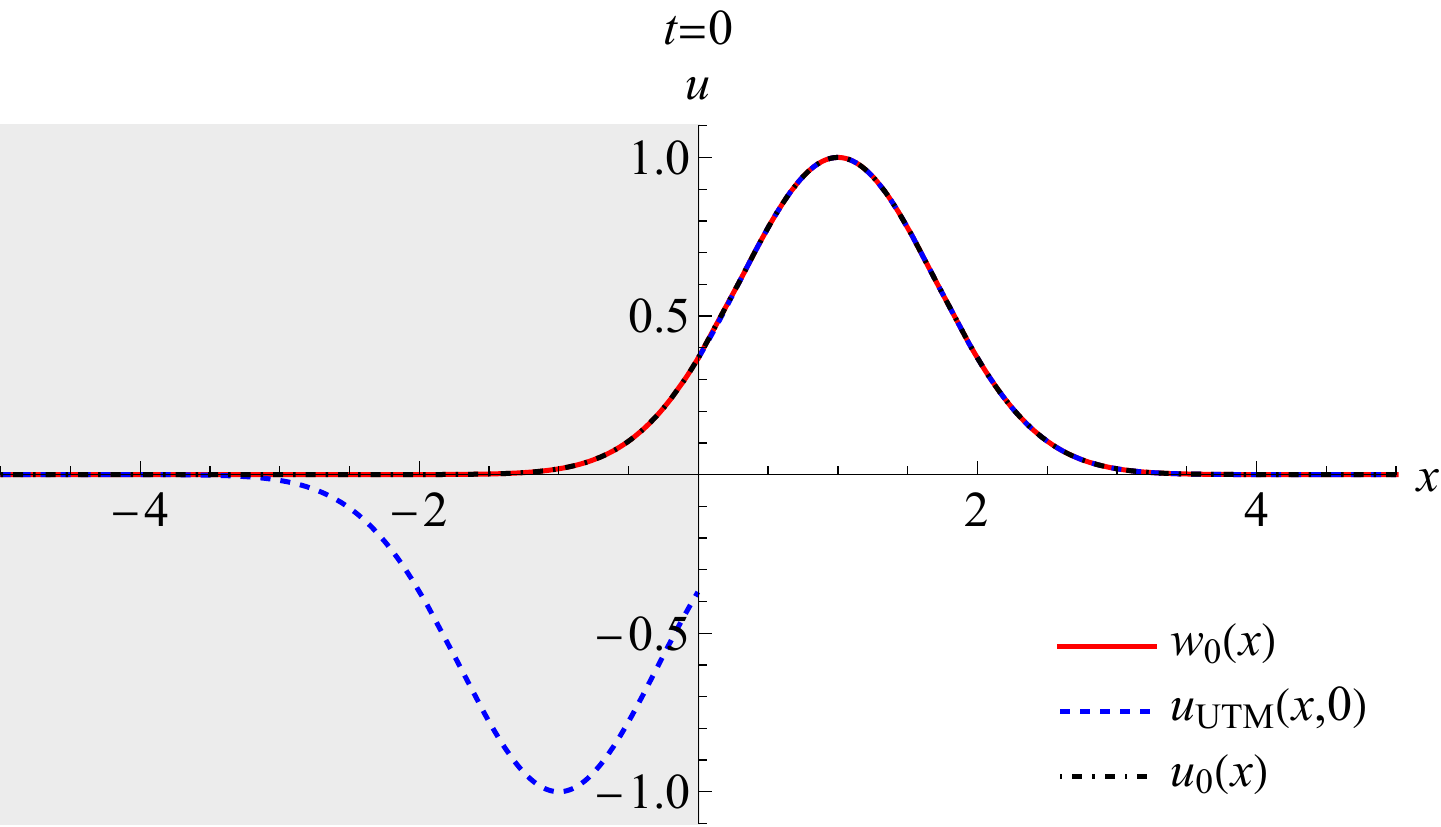} & \includegraphics[scale=\sc]{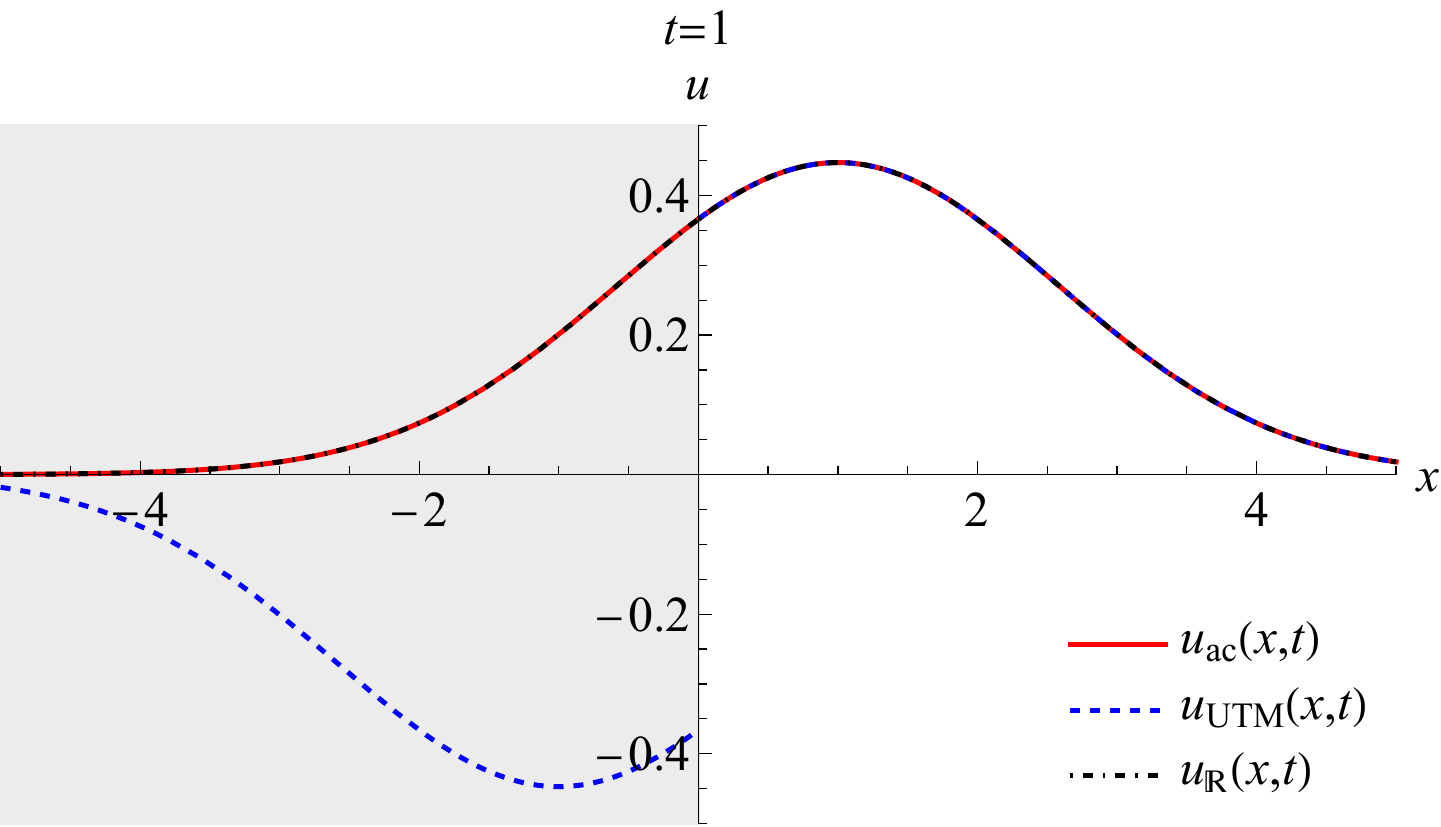}\\
(a) & (b)
\end{tabular}
\end{center}
\caption{(a) The initial condition $w_0(x)$, given by \rf{heatacw0}, leading to the (analytic) solution on the right, shown with the UTM solution at $t=0$ and the whole-line initial condition, $u_0(x) = u_\R(x,0)$. (b) The solution $u_{\text{ac}}(x,t)$ at $t=1$ obtained through analytic continuation, overlaid with the solution $u_{\R}(x,t)$ \rf{wholeline1} and the discontinuous extension $u_{\text{UTM}}(x,t)$ resulting from evaluating the UTM solution for negative $x$ values.}
\la{fig:DHL_u1}
\end{figure}

Next, we consider \rf{heatdirichletibvp} with $f_0(t) = te^{-t}$ and $u_0(x)=u_{\R}(x,0)$ for $x>0$. Since $f_0^{(n)}(t) = (-1)^n e^{-t}(t-n),$ we obtain 
\beq
\tilde f_0(x,t) = 2e^{-t}\sum_{n=0}^\infty \frac{(-1)^n x^{2n}(t-n)}{(2n)!} = e^{-t}\big(2t\cos(x)+x\sin(x)\big).
\eeq 
\no Using \rf{heatacw} and \rf{heatacw0}, we find the analytic continuation for $x\in \R$. In this case, $w_0(x)$ is discontinuous at $x=0$, (since the first compatibility condition is not satisfied), see Figure~\ref{heatacex2}. 

\begin{figure}[tb]
\begin{center}
\def \sc {0.475}
\begin{tabular}{cc}
\includegraphics[scale=\sc]{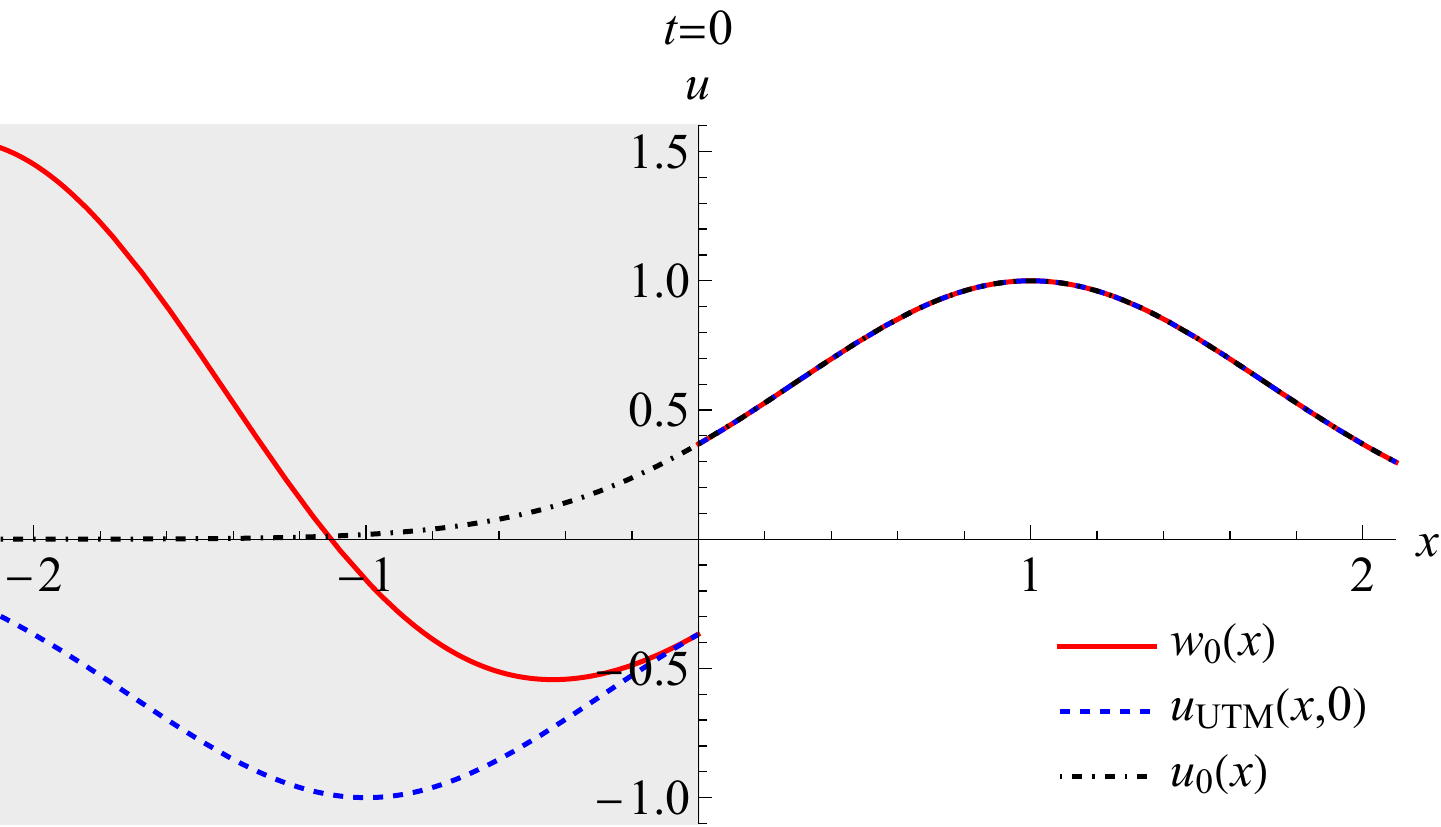} & \includegraphics[scale=\sc]{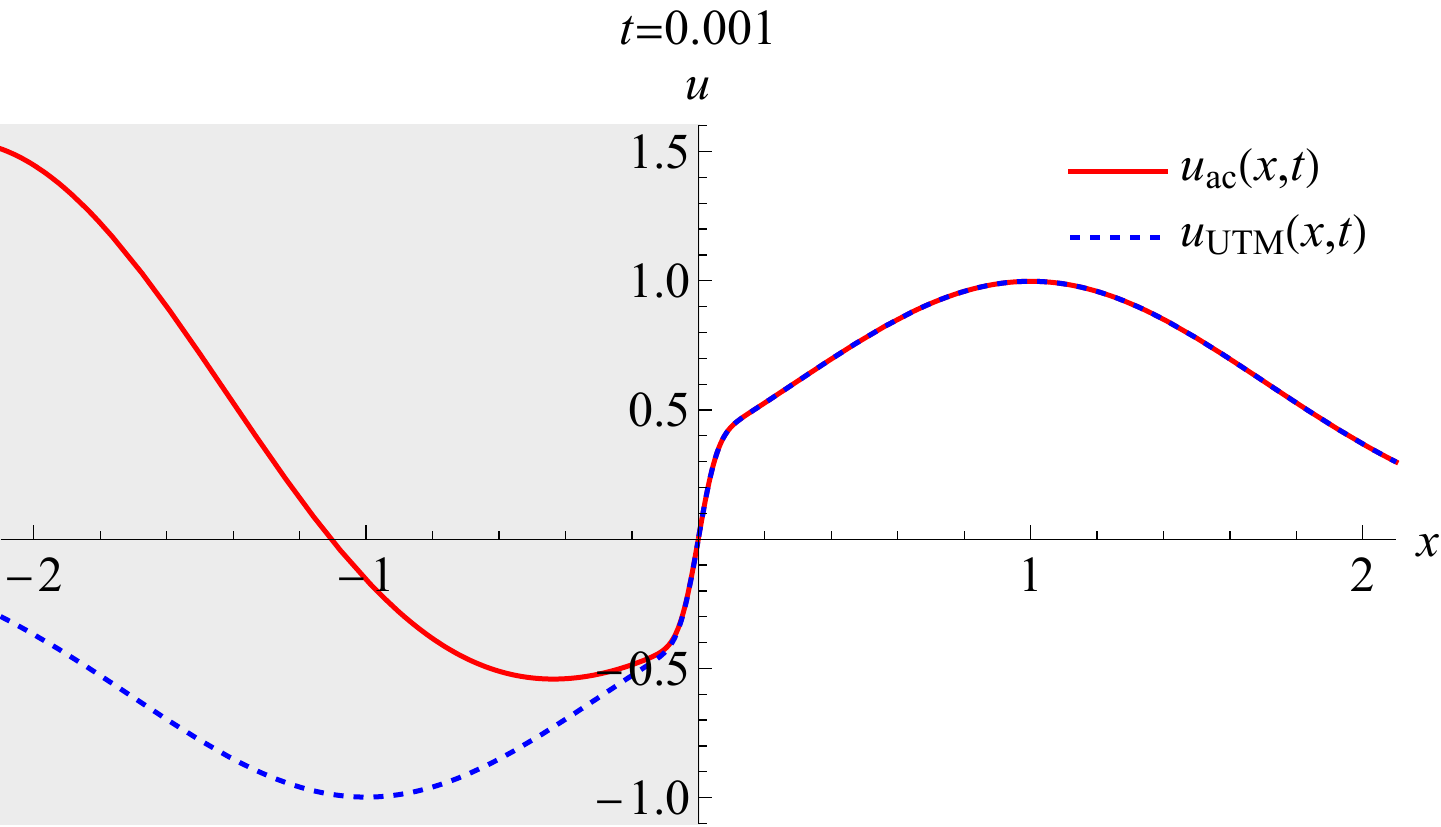}\\
(a) & (b)\\
\includegraphics[scale=\sc]{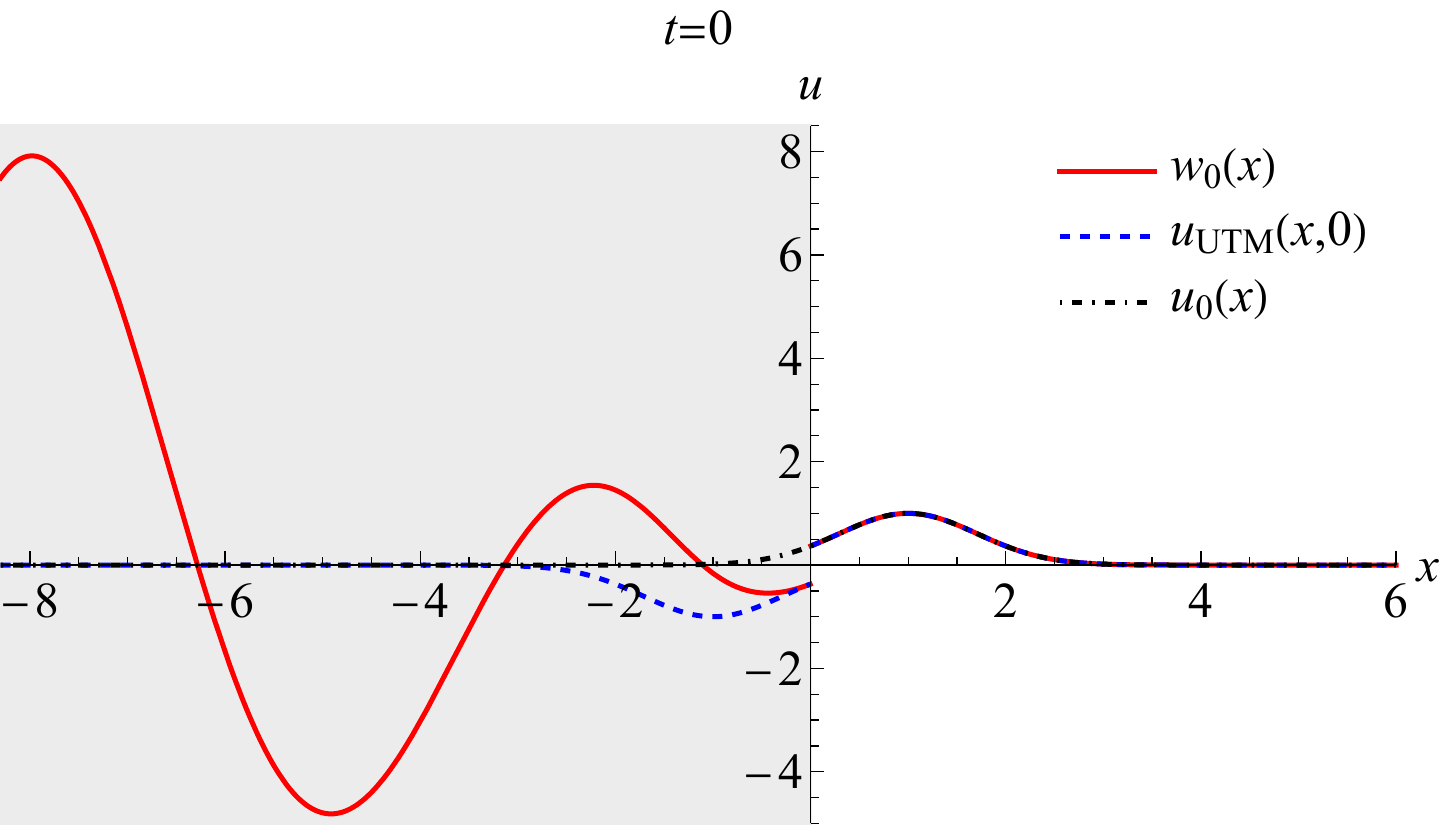} & \includegraphics[scale=\sc]{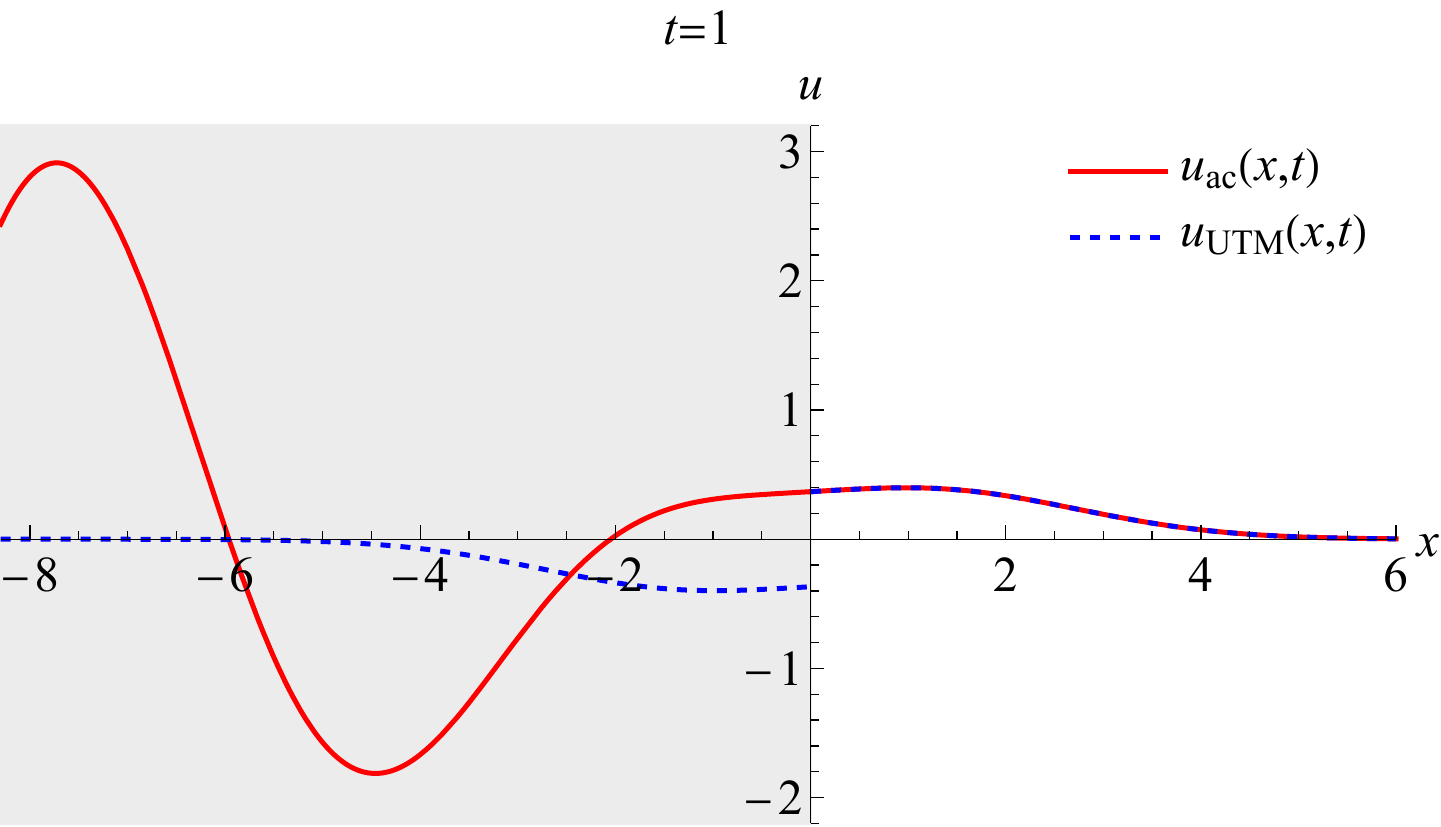}\\
(c) & (d)
\end{tabular}
\end{center}
\caption{(a, c) The (discontinuous) initial condition $w_0(x)$ \rf{heatacw0} leading to the (analytic) solution on the right, shown with $u_\text{UTM}(x,0)$ and the analytic continuation of $u_0(x)$. Note that (a) is a close-up of (c). (b) The solution $u_{\text{ac}}(x,t)$ \rf{heatacw} at $t=0.001$ obtained through analytic continuation, shown with the discontinuous UTM solution $u_{\text{UTM}}(x,t)$. Note that the jump in $u_\text{UTM}(x,t)$, $u_\text{UTM}(0^+,t)-u_\text{UTM}(0^-,t)= 2f_0(t)\approx 0.002$. 
(d) The solution $u_{\text{ac}}(x,t)$ at $t=1$ obtained through analytic continuation, shown with the discontinuous UTM solution $u_{\text{UTM}}(x,t)$. }
\la{heatacex2}
\end{figure}

\subsection{Neumann boundary conditions}\label{heat_N_cont_sec}

\def \contour {{\partial\Omega}}

The heat equation on the half line with Neumann boundary conditions, 
\begin{subequations}\la{heatneumannibvp}
\begin{align}
u_t &= u_{xx}, &&\hspace*{-1.0in} x>0, \, t>0,\hspace*{1.0in}\\
u(x,0) &= u_0(x), &&\hspace*{-1.0in} x > 0, \\ \la{heatneumannbc}
u_x(0,t) &= f_1(t),&&\hspace*{-1.0in} t>0,
\end{align}
\end{subequations}  
\no has the UTM solution \cite{JC_fokas_book}
\beq 
u(x,t) = I_0(x,t) + I_{f_1}(x,t), 
\eeq 
\no with
\begin{align}\la{heatneumanni0}
I_0(x,t) &= \frac{1}{2\pi} \int_{-\infty}^\infty e^{ikx-k^2t} \hat u_0(k)\,dk  + \frac{1}{2\pi} \int_\contour e^{ikx - k^2t} \hat u_0(-k)\,dk, \\
I_{f_1}(x,t) &=   - \frac{1}{\pi} \int_\contour e^{ikx - k^2t}F_1(k^2,t)\,dk ,
\end{align}
\no and 
\beq 
\hat u_0(k) = \int_0^\infty e^{-iky} u_0(y) \, dy , \qquad
F_1(k^2,t) = \int_0^t e^{k^2s} f_1(s) \, ds.
\eeq 
\no The region $\Omega$ remains as in Section~\ref{sec:DBC}, Figure~\ref{fig:Omega}.

%
%
As before, $I_0(x,t)$ 
%
%
%
is an entire function of $x\in \C$ for $t>0$ by Theorem~\ref{thm:IC}. Again as before, the integration with respect to $k$ in $I_{f_1}(x,t)$ can be evaluated so that
\beq 
\label{If1}
I_{f_1}(x,t) = -\frac{1}{\pi} \int_{0}^t ds \, f_1(s) \int_{\partial \Omega} e^{ikx-k^2(t-s)}\,dk = -\frac{1}{\sqrt{\pi}}\int_0^t  \frac{f_1(s)}{\sqrt{t-s}}e^{-\frac{x^2}{4(t-s)}}\,ds,
\eeq
\no which is analytic for $x>0$ (using a result analogous to Theorem~\ref{thm:IC}), but not at $x=0$, although it is continuous there. This is an even function of $x$ (immediately seen using the classical cosine transform representation or from \rf{If1}) and not analytic at $x=0$ unless $f_1(t)\equiv 0$, as we see below, \rf{eqn:If1_ext}. We extend $I_{f_1}(x,t)$ to an entire function by finding a Taylor series representation about $x=0$: 
\beq 
I_{f_1}(x,t) = 2\sum_{n=0}^\infty b_{2n+1}(t) x^{2n+1} + I_{f_1}(-x,t),
\eeq 
\no where we use the odd terms in the Taylor series to extend the function to $x<0$. This choice is convenient because the Neumann condition 
\rf{heatneumannbc} gives $b_1(t) = f_1(t)$. Proceeding as before, we have
\beq
\frac{\partial^{2n+1}I_{f_1}}{\partial x^{2n+1}}\bigg|_{x=0} =  \frac{(-1)^{n}}{i\pi} \int_0^t ds \, f_1(s) \int_{\contour} k^{2n+1} e^{-k^2(t-s)}\,dk  =   f_1^{(n)}(t).
\eeq 

Thus, for $x<0$, 
\beq
I_{f_1}(x,t) = \tilde f_1(x,t)+I_{f_1}(-x,t),
\eeq 
\no where we define 
\beq \label{eqn:f1tilde}
\tilde f_1(x,t) = 2\sum_{n=0}^\infty \frac{x^{2n+1}}{(2n+1)!} f_1^{(n)}(t),
\eeq 
\no so that 
\beq \label{eqn:If1_ext}
I_{f_1}^\text{ext}(x,t) = \case{1.2}{I_{f_1}(x,t), & x\geq 0, \\ \tilde f_1(x,t) + I_{f_1}(-x,t), & x< 0,}
\eeq 
\no is an analytic function for $x\in \R$. At this point, it is clear that in order for this to equal the even function $I_{f_1}(x,t)$, 
we need $f_1(t)\equiv 0$. If we want the full series for $I_{f_1}(x,t)$, we can use integration by parts in the same way as in Section~\ref{sec:DBC}, to get the coefficients of the even terms in the Taylor series as
\beq 
b_{2n}(t) = - \frac{\Gamma\left(\frac12\right)}{\pi(2n)!} \left[\sum_{m=1}^n \frac{(-1)^{n-m}\Gamma\left(n-m+\frac12\right)}{t^{n-m+\frac12}\Gamma\left(\frac12\right)} + \int_0^t \frac{f_1^{(n)}(s)}{(t-s)^{\frac12}}\,ds\right],
\eeq 
\no allowing us to rewrite $I_{f_1}(x,t)$ as 
\beq 
I_{f_1}(x,t) = -\sum_{n=0}^\infty \frac{(-1)^n x^n}{n!}f_1^{\left(\frac{n-1}2\right)}(t).
\eeq 
\no As before, the radius of convergence of this series is infinite under the assumptions of Theorem~\ref{thm2}. 

\subsubsection{Boundary-to-Initial Map}



Considering \rf{eqn:f1tilde}, we see that
$\tilde f_{1t}=\tilde f_{1xx}$. 
Therefore 

\beq \label{heatnbcac}
w(x,t) = u_\text{ac}(x,t) = \case{1.2}{ u(x,t), & x\geq 0, \\ \tilde f_1(x,t)+u(-x,t), & x<0,}
\eeq 

\no is an analytic solution to the heat equation on the whole line with initial condition 

\beq
\label{eqn:hlnbcw0}
w_0(x) = \case{1.2}{ u_0(x), & x\geq 0, \\ \tilde f_1(x,0) + u_0(-x), & x<0. }
\eeq 

\no This initial condition is analytic when the compatibility conditions $u_0^{(2n+1)}(0) = f_1^{(n)}(0)$, $n\geq 0$ are satisfied \cite{trogdon}. If $u_0(x)$ is analytic and these corner compatibility conditions are met, then $w_0(x) = u_0(x)$. As before, $w_0(x)$ is not necessarily bounded for $x<0$. Note that for the homogeneous boundary condition $f_1(t)=0$, we recover the method of images' even extension of the initial condition.

\subsubsection{Examples}

As an example, we construct a half-line problem using \rf{wholeline1} but with the appropriate Neumann boundary condition. Since the solution of the Neumann IBVP, as given by the UTM, is even if we allow $x<0$, we obtain a continuous, but non-differentiable function.  Using the analytic continuation described above, $u_{\R}(x,t)$ is recovered, see Figure~\ref{fig:heatneumannfig1}.

\begin{figure}[tb]
\begin{center}
\def \sc {0.475}
\begin{tabular}{cc}
\includegraphics[scale=\sc]{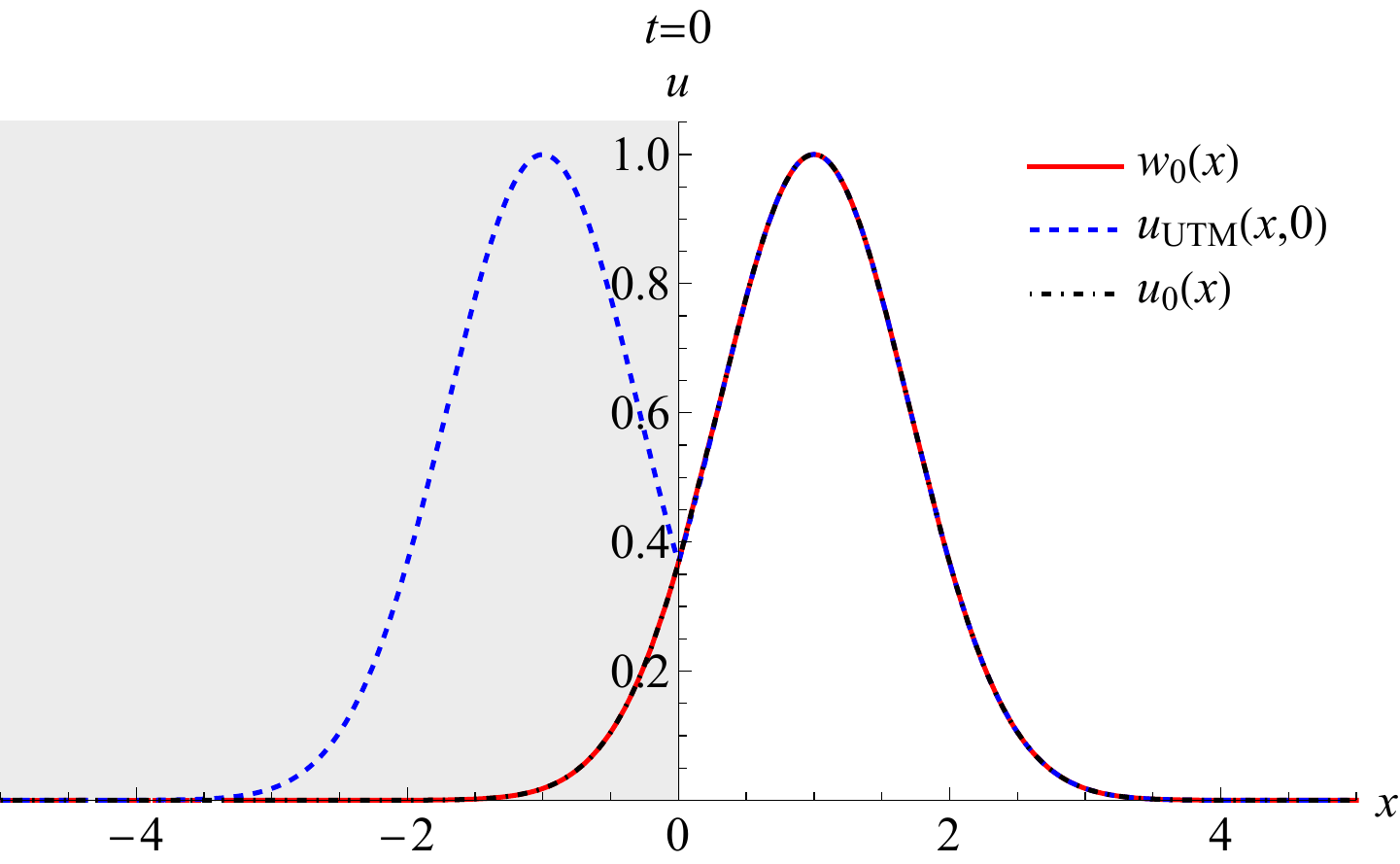} & \includegraphics[scale=\sc]{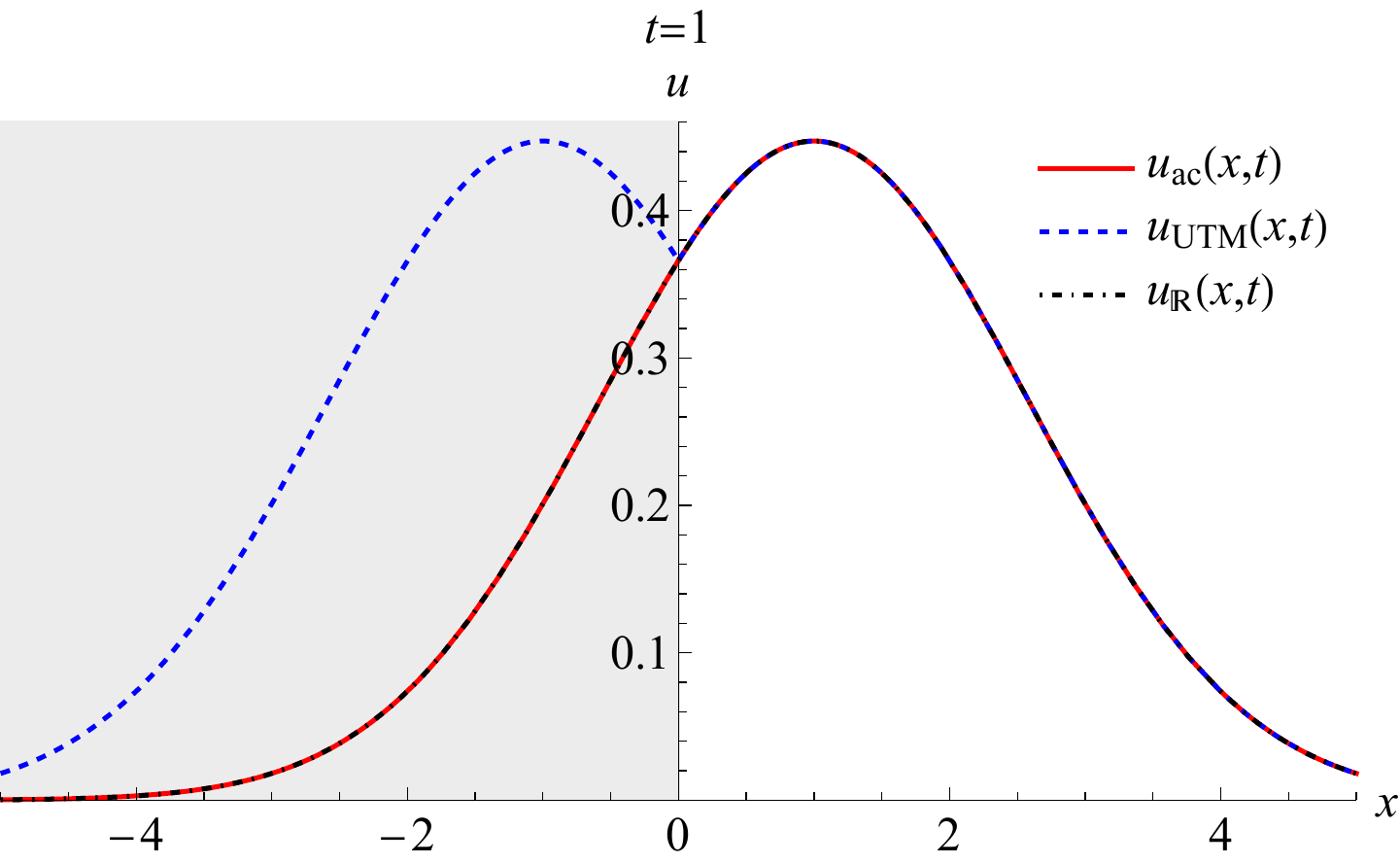}\\
(a) & (b)
\end{tabular}
\end{center}
\caption{(a) The initial condition $w_0(x)$ given by \rf{eqn:hlnbcw0}, leading to the (analytic) solution on the right, shown with the UTM solution at $t=0$ and the whole-line initial condition $u_0(x)=u_\R(x,0)$. (b) The solution $u_{\text{ac}}(x,t)$  \rf{heatnbcac} at $t=1$ obtained through analytic continuation, overlaid with the solution $u_{\R}(x,t)$ \rf{wholeline1} and the nondifferentiable extension $u_{\text{UTM}}(x,t)$ resulting from evaluating the UTM solution for negative $x$ values.}
\la{fig:heatneumannfig1}
\end{figure}

Using instead the boundary function $f_1(t) = te^{-t}$, 

\beq
\tilde f_1(x,t) = e^{-t}\big((2t-1)\sin(x)-x\cos(x)\big). 
\eeq  

\no In this case, the resulting analytic continuation is a smooth function, but the corresponding initial condition is not differentiable at $x=0$, see Figure~\ref{fig:heatneumannfig2}

\begin{figure}[tb]
\begin{center}
\def \sc {0.475}
\begin{tabular}{cc} 
\includegraphics[scale=\sc]{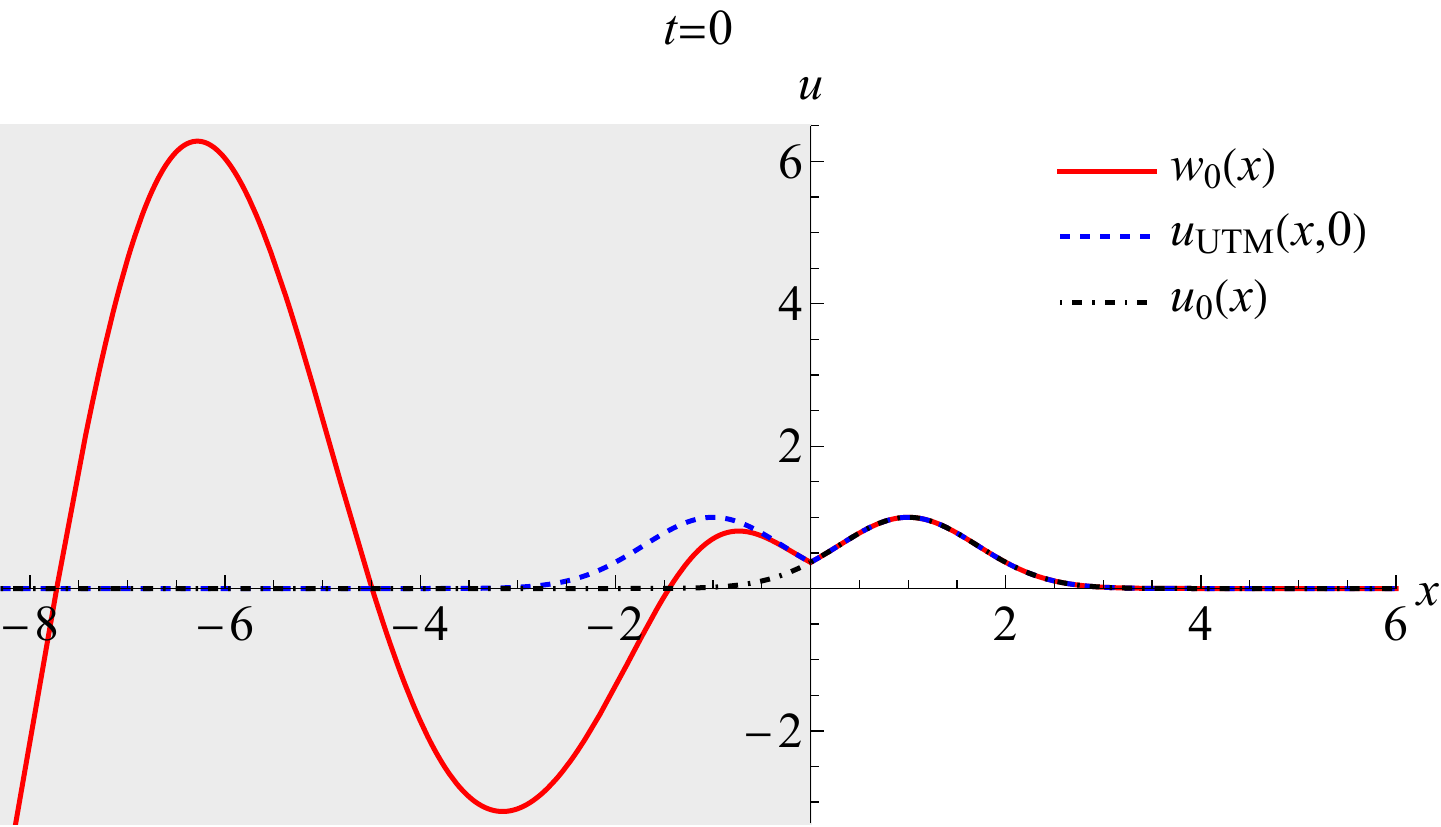} &\includegraphics[scale=\sc]{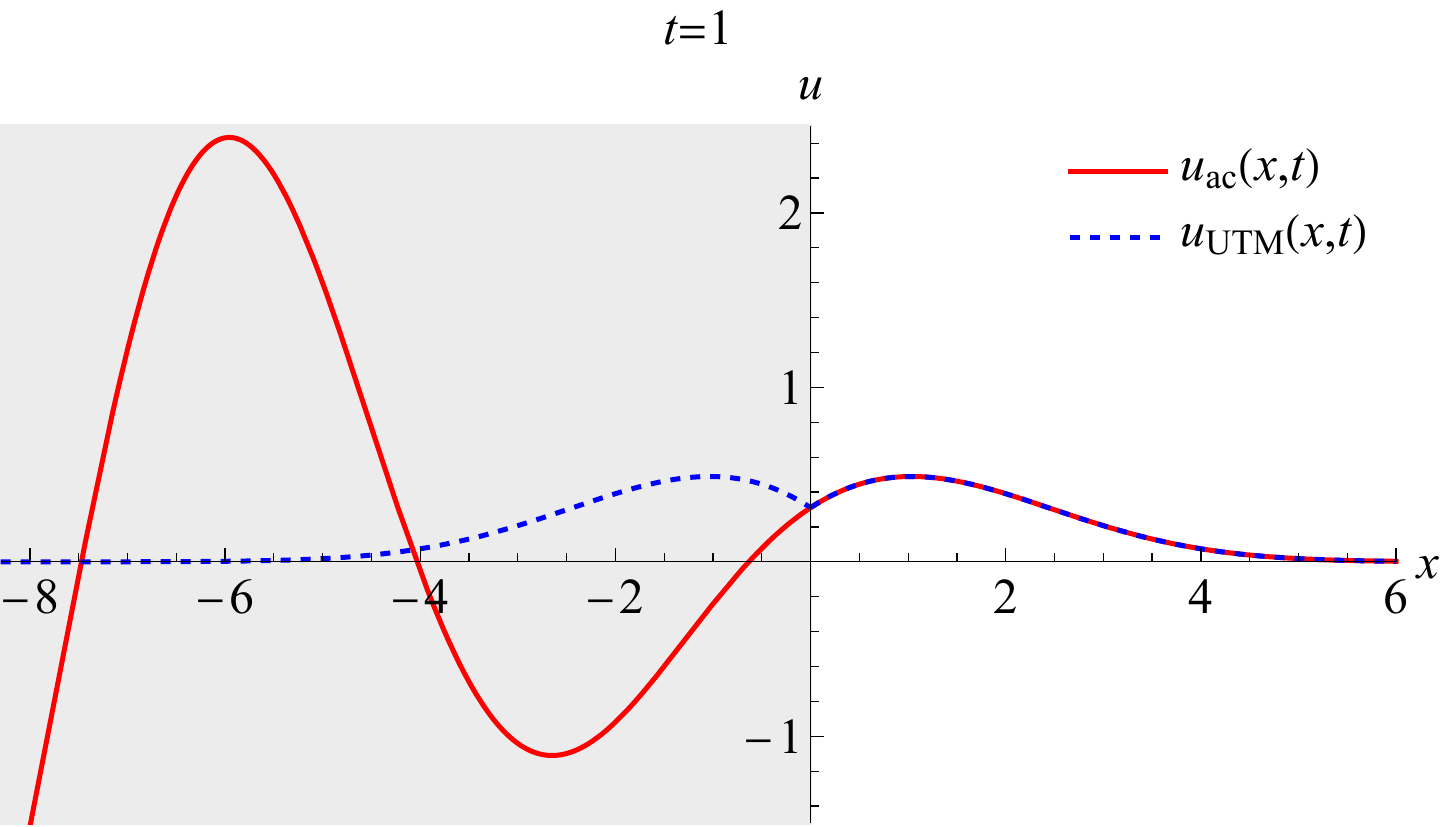}\\
(a) & (b)
\end{tabular}
\end{center}
\caption{(a) The (non-differentiable) initial condition $w_0(x)$ \rf{eqn:hlnbcw0} leading to the (analytic) solution on the right, shown with the analytic continuation of $u_0(x)$. (b) The solution $u_{\text{ac}}(x,t)$  \rf{heatnbcac} at $t=1$ obtained through analytic continuation, shown with the even extension of $u_\text{UTM}(x,t)$.}
\la{fig:heatneumannfig2}
\end{figure}

\subsection{Finite interval with Dirichlet boundary conditions}

\def \contour {{\partial\Omega}}

The solution to the heat equation on a finite interval with Dirichlet boundary conditions, 
\begin{subequations}\la{heatfiibvp}
\begin{align}\la{heatfi}
u_t&=u_{xx}, &&\hspace*{-1.0in} x\in(0,L),\, t>0,\hspace*{1.0in}\\
u(x,0)&=u_0(x), &&\hspace*{-1.0in} x\in (0,L),\\
u(0,t)&=f_0(t), &&\hspace*{-1.0in} t>0,\\
u(L,t)&=g_0(t), &&\hspace*{-1in} t>0,
\end{align}
\end{subequations}
\no can be written as \cite{JC_fokas_book} 
\beq 
u(x,t) = I_0(x,t) + I_{f_0}(x,t) + I_{g_0}(x,t),
\eeq 
\no with 
\begin{align} \nonumber
I_0(x,t) &=  \frac{1}{2\pi} \int_{-\infty}^\infty e^{ikx-k^2t} \hat u_0(k)\,dk 
- \frac{1}{2\pi} \int_{\partial\Omega} e^{-k^2t} \frac{e^{ik(L+x)} -e^{ik(L-x)}}{\Delta(k)}\hat u_0(k)\,dk\\\la{heatfii0}
& -\frac{1}{2\pi} \int_{\partial\Omega} e^{-k^2t} \frac{e^{ik(L-x)} - e^{-ik(L-x)}}{\Delta(k)} \hat u_0(-k)\,dk, \\
I_{f_0}(x,t) &= \frac{1}{i\pi} \int_{\contour} ke^{-k^2t} \frac{e^{ik(L-x)} - e^{-ik(L-x)}}{\Delta(k)}F_0(k^2,t)\,dk, \\
I_{g_0}(x,t) &= \frac{1}{i\pi} \int_{\partial  \Omega} ke^{-k^2t} \frac{e^{ikx}-e^{-ikx}}{\Delta(k)}G_0(k^2,t)\,dk. 
\end{align}
\no Here
\begin{align}
\hat u_0(k) = \int_0^L e^{-iky} u_0(y)\,dy, &\qquad \Delta(k) = 2i\sin(kL), \\
F_0(k^2,t) = \int_0^t e^{k^2s} f_0(s)\,ds, &\qquad G_0(k^2,t) = \int_0^t e^{k^2s} g_0(s)\,ds,
\end{align}
\no and $ \Omega=\{$
\!\!$k\in \C: |k|>r$, and $\pi/4 < \Arg(k) < 3\pi/4 \}$, for some $r>0$, shown in Figure~\ref{fig:Omega}.

The integral $I_0(x,t)$ 
%
%
is an entire function of $x$ by Theorem~\ref{thm:IC_FI} below. 

\begin{theorem}\label{thm:IC_FI}
If $u_0\in L^1(0,L)$, then $I_0(x,t)$ is entire in $x$ for $t>0$. If $f_0\in L^1(0,T)$ for some $T>0$, then for $0\leq t \leq T$, $I_{f_0}(x,t)$ is analytic for $x\in \{x\in C: |\Im(x)|<\Re(x) \text{ and } |\Im(2L-x)|<\Re(2L-x)\}$. If $g_0 \in L^1(0,T)$, then for $0\leq t \leq T$, $I_{g_0}(x,t)$ is analytic for $x\in \{x\in C: |\Im(L+x)|<\Re(L+x) \text{ and } |\Im(L-x)|<\Re(L-x)\}$.
\end{theorem}

\begin{proof} The proof is similar to the proof of Theorem~\ref{thm:IC}. 

\end{proof}

The function $I_{f_0}(x,t)$ is defined and analytic only for $x\in(0,2L)$, as otherwise its exponential kernel is growing on $\contour$. Deforming $\contour$ to the real axis, picking up residue contributions from the singularities, we obtain the classical Fourier series:
\beq
I_{f_0}(x,t) = \sum_{n=1}^\infty \frac{2n\pi}{L^2}e^{-\frac{n^2\pi^2t}{L^2}}F_0\left(\frac{n^2\pi^2}{L^2},t\right)\sin\left(\frac{n\pi x}{L}\right).
\eeq
\no This Fourier series is defined outside $x\in(0,2L)$, but this periodic extension is not in general analytic, see Figure~\ref{heatacfiex1}. Upon integrating by parts,
\beq
F_0(k^2,t) =  \int_0^t e^{k^2s} f_0(s)\,ds = \frac{f_0(t)e^{k^2t}- f_0(0)}{k^{2}} - \frac{1}{k^{2}} \int_0^t  e^{k^2s}f_0'(s)\,ds \sim \mathcal O\left(k^{-2}e^{k^2t}\right),
\eeq
\no and, due to a lack of convergence when taking $x$-derivatives, a Taylor series cannot be obtained from this representation.  

Using the contour integral representation, we can analytically continue $I_{f_0}(x,t)$ for $x<0$, by finding a Taylor series about $x=0$. As before, 
\beq 
I_{f_0}(x,t) = 2\sum_{n=0}^\infty a_{2n}(t)x^{2n} - I_{f_0}(-x,t)=\tilde f_0(x,t) - I_{f_0}(-x,t),  
\eeq 
\no which defines $\tilde f_0(x,t)$. Using
\beq 
I_{f_0}(x,t) = \frac{1}{i\pi} \int_0^t ds \, f_0(s) \int_{\partial\Omega} k\frac{\sin(k(L-x))}{\sin(kL)}e^{-k^2(t-s)}\,dk,
\eeq 
\no if we take $2n$ $x$-derivatives, we find
\beq 
\frac{\partial^{2n}I_{f_0}}{\partial x^{2n}}\bigg|_{x=0} = \frac{(-1)^n}{i\pi} \int_0^t ds \, f_0(s)\int_{\partial\Omega} k^{2n+1} e^{-k^2(t-s) }\,dk= f_0^{(n)}(t). 
\eeq  
\no It follows that, for $x\in (-2L,0)$, as before, 
\beq
I_{f_0}(x,t) = \tilde f_0(x,t)-I_{f_0}(-x,t).
\eeq 
\no Thus, 
\beq 
I_{f_0}^\text{ext}(x,t) = \case{1}{-I_{f_0}(-x,t)+\tilde f_0(x,t), & -2L<x<0, \\ f_0(t), & x=0, \\ I_{f_0}(x,t), & 0< x < 2L,}
\eeq
\no is an analytic function for $-2L<x<2L$. For the remainder of this section, we define $I_{f_0}(0,t) = f_0(t)$ and $I_{g_0}(L,t) = g_0(t)$, so as to not write the boundary terms separately. Since for $0<x<2L$, $I_{f_0}(2L-x,t) = -I_{f_0}(x,t)$,
\beq 
I_{f_0}^\text{ext}(x-2L,t) = \case{1}{I_{f_0}(x,t)+\tilde f_0(x-2L,t), & 0<x<2L, \\ I_{f_0}(x-2L,t), & 2L\leq x < 4L,}
\eeq 
\no is an analytic function for $0<x<4L$, so that
\beq 
I_{f_0}^\text{ext}(x,t) = \case{1}{-I_{f_0}(-x,t)+\tilde f_0(x,t), & -2L<x<0, \\ I_{f_0}(x,t), & 0\leq x < 2L,\\ I_{f_0}(x-2L,t)-\tilde f_0(x-2L,t), & 2L\leq x<4L,} 
\eeq
\no is an analytic function for $-2L<x<4L$. Continuing this, we find
\beq\la{heatif0ext}
I_{f_0}^\text{ext}(x,t) = \case{1}{\vdots \\ I_{f_0}(x+4L,t)+\tilde f_0(x,t)+\tilde f_0(x+2L,t), & -4L\leq x<-2L, \\I_{f_0}(x+2L,t)+\tilde f_0(x,t), & -2L\leq x<0, \\ I_{f_0}(x,t), & 0 \leq x < 2L,\\ I_{f_0}(x-2L,t)-\tilde f_0(x-2L,t), & 2L\leq x<4L, \\I_{f_0}(x-4L,t)-\tilde f_0(x-2L,t)-\tilde f_0(x-4L,t), & 4L\leq x<6L, \\ \vdots }
\eeq



We can repeat the analysis similar to what we did for the whole-line problem, by finding the odd coefficients. Alternatively, we can expand about a point other than $x=0$. For instance, denoting the Taylor series about $x=L$ as
\beq 
\label{heat_fi_fullseries}
I_{f_0}(x,t) = \sum_{n=1}^\infty A_{2n-1}(t) (x-L)^{2n-1},
\eeq 
\no (the coefficients of the even terms vanish), we have 
\beq 
A_{2n-1}(t)= \frac{1}{(2n-1)!} \frac{\partial^{2n-1}I_{f_0}}{\partial x^{2n-1}} \bigg|_{x=L} = \frac{2(-1)^n}{\pi(2n-1)!} \int_0^t ds \, f_0(s) \int_{\partial \Omega} k^{2n}\frac{e^{ikL}}{e^{2ikL}-1} e^{-k^2(t-s)}\,dk. 
\eeq 
\no which provides an analytic continuation of $I_{f_0}(x,t)$ for all $x\in \C$, as shown below.

\begin{theorem}

If $f_0(t)$ is analytic in a strip in the complex $t$-plane, containing part of the positive real axis, $\mathcal D = \{t\in \C: \text{dist}(t,[0,T])\leq r\}$, for some $r>0$ and $T>0$, then the series representation \rf{heat_fi_fullseries} is entire in $x$ for $0\leq t\leq T$.

\end{theorem}

\begin{proof}

The proof is similar to the proof of Theorem~\ref{thm2}.

\end{proof}

We proceed for $I_{g_0}(x,t)$ in exactly the same way as for $I_{f_0}(x,t)$. First, note that $I_{g_0}(x,t)$ is analytic for 
$x\in (-L,L).$ We write
\beq
I_{g_0}(x,t) = \frac{1}{i\pi} \int_0^t ds \, g_0(s)\int_{\partial \Omega} k \frac{\sin(kx)}{\sin(kL)}e^{-k^2(t-s)}\,dk. 
\eeq 
\no Taking $2n$ $x$-derivatives, 
\beq 
\frac{\partial^{2n}I_{g_0}}{\partial x^{2n}}\bigg|_{x=L} = \frac{(-1)^n}{i\pi} \int_0^t ds \, g_0(s)\int_{\partial\Omega} k^{2n+1} e^{-k^2(t-s) }\,dk = g_0^{(n)}(t). 
\eeq 
\no Since 
\beq 
I_{g_0}(x,t) = 2\sum_n B_n(t)(x-L)^{2n} - I_{g_0}(2L-x,t)=\tilde g_0(x,t)- I_{g_0}(2L-x,t),
\eeq 
\no which defines $\tilde g_0(x,t)$. This implies that for $x\in(L,3L)$, 
\beq
I_{g_0}(x,t) = 2\sum_{n=0}^\infty \frac{(x-L)^{2n}}{(2n)!}g_0^{(n)}(t) - I_{g_0}(2L-x,t) = \tilde g_0(x,t)-I_{g_0}(2L-x,t).
\eeq 
\no This process can be continued as before to get
\beq\la{heatig0ext}
I_{g_0}^\text{ext}(x,t) = \case{1}{\vdots  \\ I_{g_0}(x+4L,t) - \tilde g_0(x+2L,t)-\tilde g_0(x+4L,t), & -5L< x\leq-3L, \\ I_{g_0}(x+2L,t)-\tilde g_0(x+2L,t), & -3L< x\leq-L, \\ I_{g_0}(x,t),&-L< x\leq L, \\ I_{g_0}(x-2L,t)+\tilde g_0(x,t), & L< x\leq 3L, \\ I_{g_0}(x-4L,t)+\tilde g_0(x,t) + \tilde g_0(x-2L,t), & 3L < x\leq 5L, \\ \vdots}
\eeq 
\no As for $I_{f_0}(x,t)$, we can also find the odd terms in the Taylor series or write the Taylor series about $x=0$ instead. 

\subsubsection{Boundary-to-Initial Map}

Using the definitions for $\tilde f_0(x,t)$ and $\tilde g_0(x,t)$ above, 
\beq \label{heatfiac}
w(x,t) = u_\text{ac}(x,t) = I_0(x,t) + I_{f_0}^\text{ext}(x,t)+ I_{g_0}^\text{ext}(x,t),
\eeq 
\no is a solution to the heat equation on the full-line with initial condition 
\beq 
\label{eqn:heatfiw0}
w_0(x) = I_0(x,0) + I_{f_0}^\text{ext}(x,0)+ I_{g_0}^\text{ext}(x,0),
\eeq 
\no where, for integer $n$, 
\begin{align} 
I_0(x,0) &= \case{0.5}{ u_0(x-2nL), &~~~ 2nL\leq x<(2n+1)L, \\ -u_0(2(n+1)L-x),&~~~ (2n+1)L \leq x<2(n+1)L,} \\
I_{f_0}^\text{ext}(x,0) &= \case{0.75}{\sum_{j=0}^{|n|-1}\tilde f_0(x+2jL,0), & 2nL\leq x<2(n+1)L, \qquad n<0, \\  -\sum_{j=1}^{n} \tilde f_0(x-2jL,0), & 2nL\leq x<2(n+1)L, \qquad n\geq 0,} \\ 
I_{g_0}^\text{ext}(x,0) &= \case{0.75}{-\sum_{j=1}^{|n|}\tilde g_0(x+2jL,t), & (2n-1)L< x\leq (2n+1)L, \qquad n<0, \\ \sum_{j=0}^{n-1} \tilde g_0(x-2jL,t), & (2n-1)L < x \leq (2n+1)L, \qquad n\geq 0.} 
\end{align} 

\subsubsection{Examples}

Using \rf{wholeline1} again, we consider the IBVP with $u_0(x)=u_\R(x,0)$, for $x\in (0,1)$, and boundary data $f_0(t)=u_{\R}(0,t)$, $g_0(t)=u_{\R}(1,t)$ for $t>0$. We recover $u_\R(x,t)$, as shown in Figure~\ref{heatacfiex1}.

\begin{figure}[tb]
\begin{center}
\def \sc {0.475}
\begin{tabular}{cc}
\includegraphics[scale=\sc]{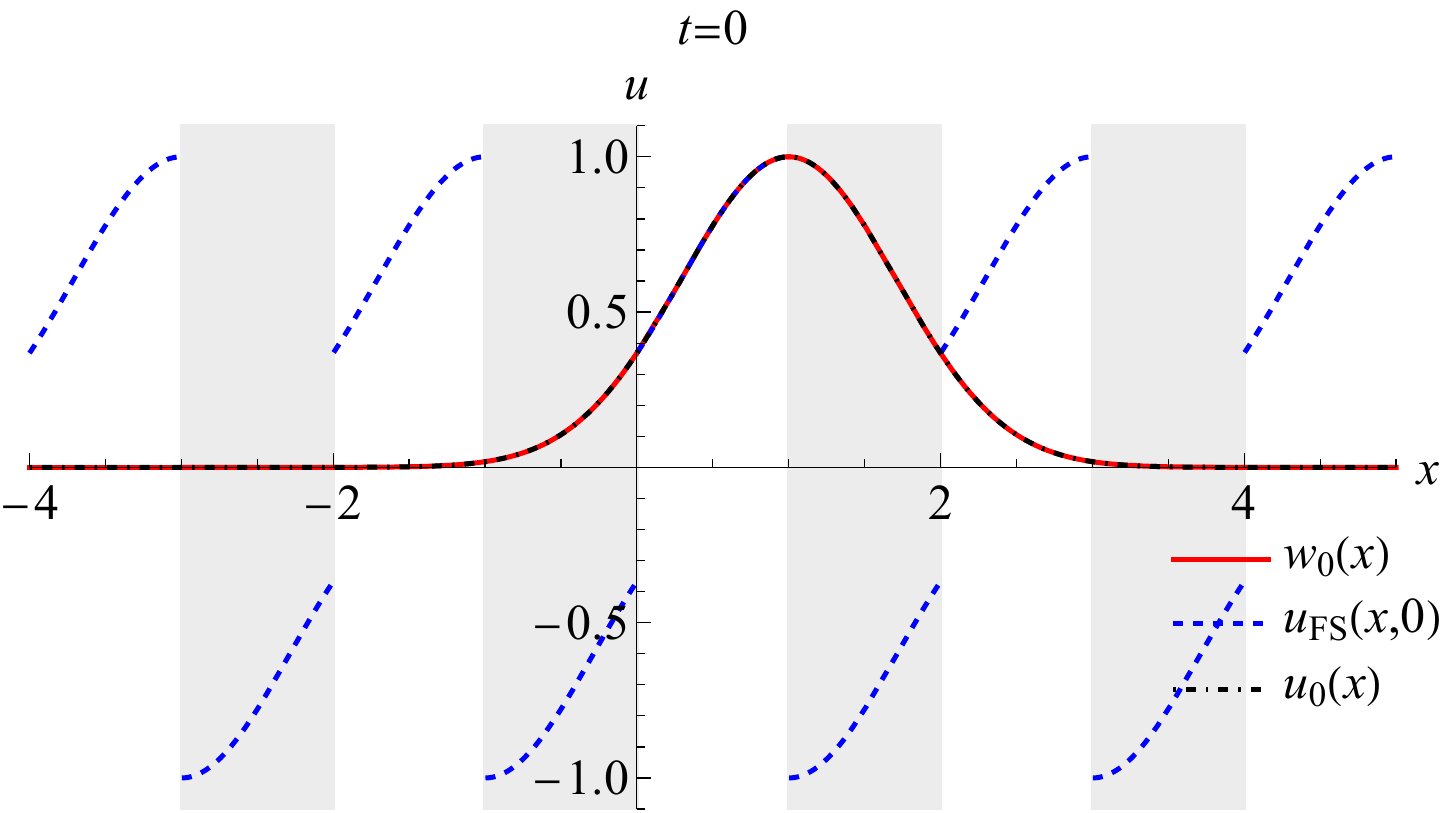} & \includegraphics[scale=\sc]{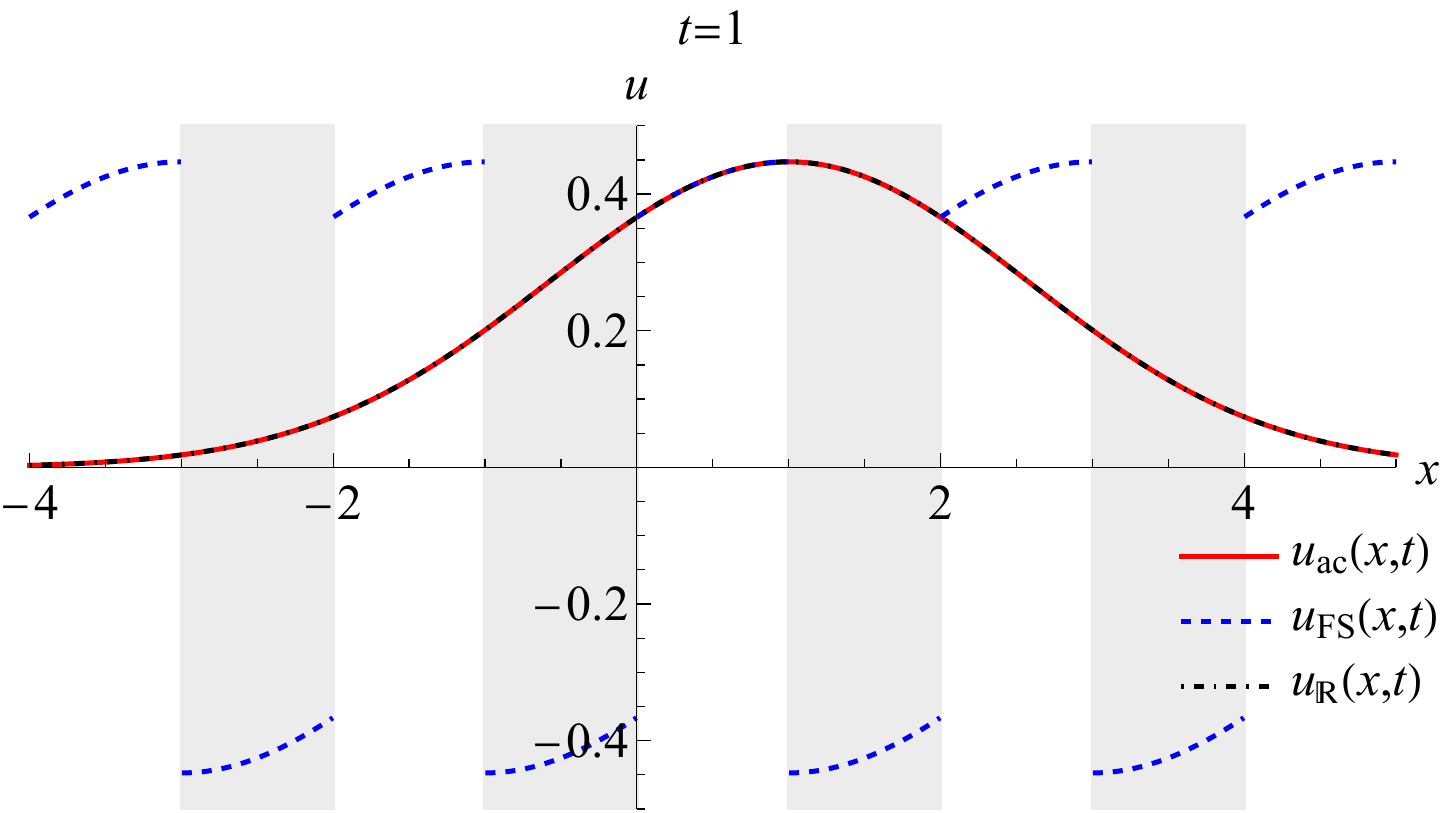}\\
(a) & (b)
\end{tabular}
\end{center}
\caption{(a) The initial condition $w_0(x)$ given by \rf{eqn:heatfiw0}, leading to the (analytic) solution on the right, shown with the Fourier Series solution at $t=0$  $u_\text{FS}(x,0)=u_\text{UTM}(x,0)$ and the whole-line initial condition $u_0(x)=u_\mathbb{R}(x,0)$. (b) The solution $u_{\text{ac}}(x,t)$ \rf{heatfiac} at $t=1$ obtained through analytic continuation, overlaid with the solution $u_{\R}(x,t)$ \rf{wholeline1} and the periodic odd Fourier series solution, $u_{\text{FS}}(x,t)=u_\text{UTM}(x,t)$.}
\la{heatacfiex1}
\end{figure}

As our second example, we change the boundary functions to $f_0(t) = te^{-t}$ and $g_0(t) = 1/(1+t)$. We show the analytic continuation 
of the solution on the full line in Figure~\ref{heatacfiex2}. As in our previous examples, it should be noted that the solution is not bounded on $\R$, nor is the extended initial condition continuous. 

\begin{figure}[tb]
\begin{center}
\def \sc {0.475}
\begin{tabular}{cc}
\includegraphics[scale=\sc]{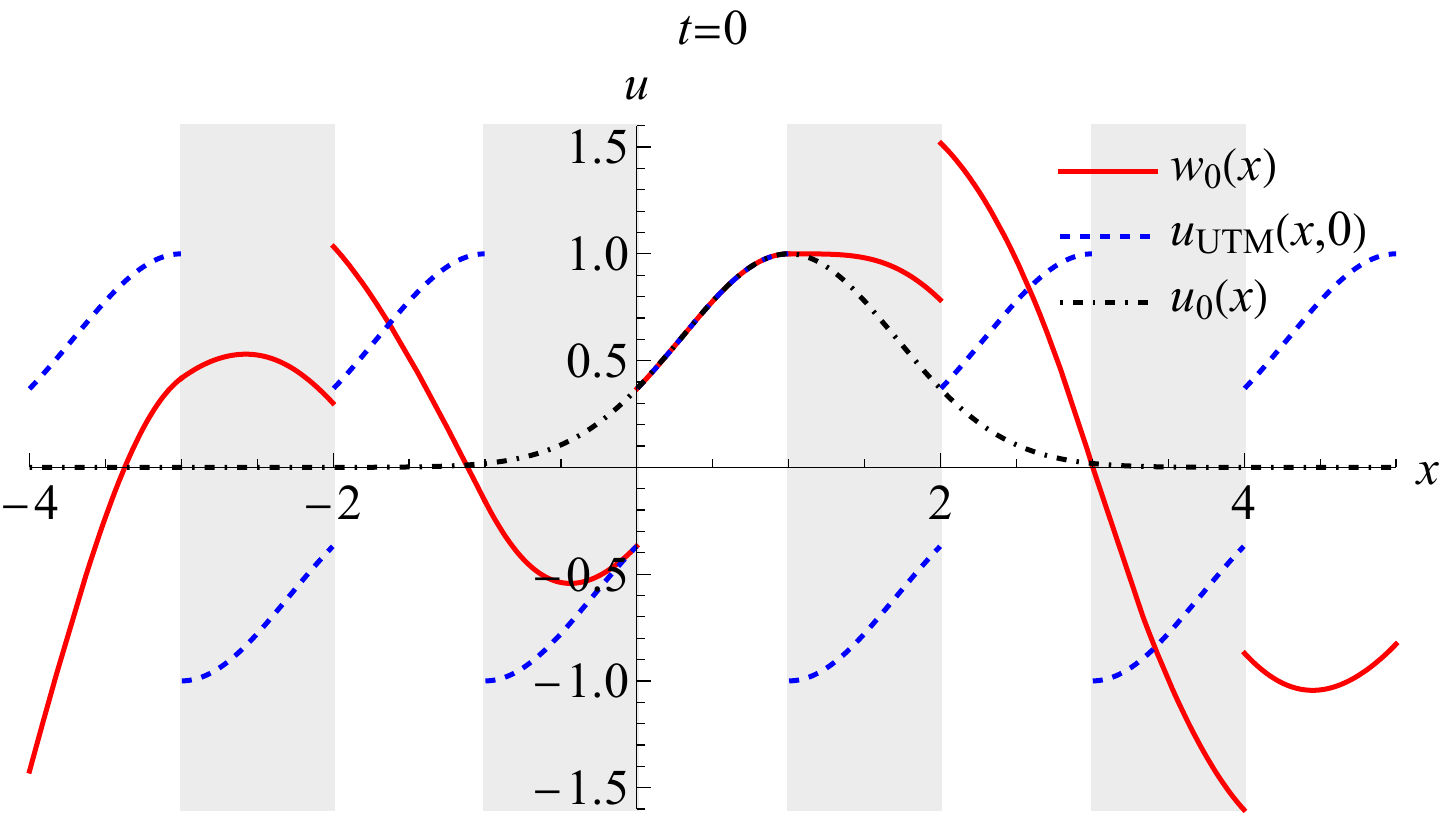} & \includegraphics[scale=\sc]{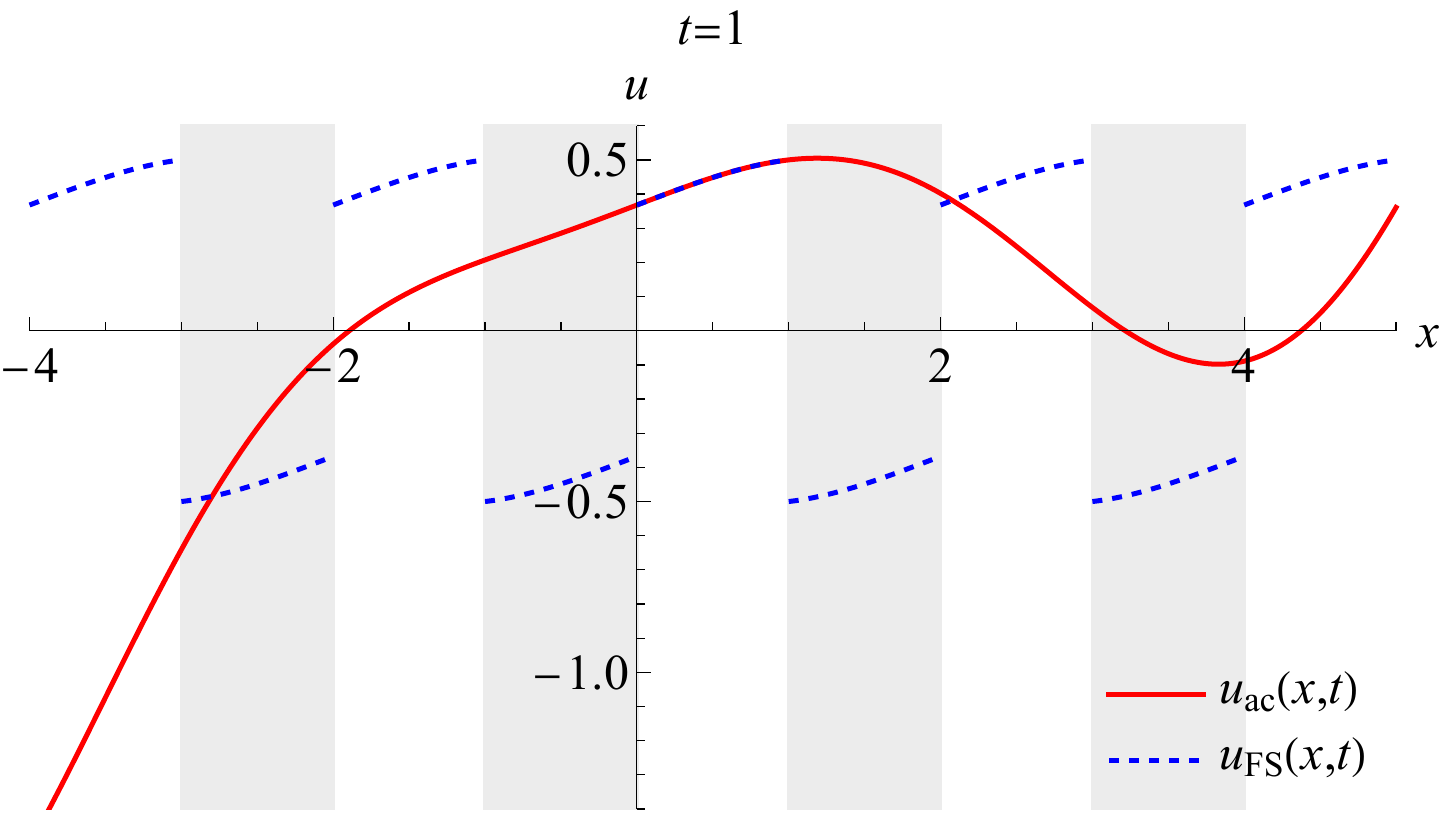}\\
(a) & (b)
\end{tabular}
\end{center}
\caption{(a) The (discontinuous) initial condition $w_0(x)$ \rf{eqn:heatfiw0} leading to the (analytic) solution on the right, shown with the analytic continuation of $u_0(x)$. (b) The solution $u_{\text{ac}}(x,t)$ \rf{heatfiac} at $t=1$ obtained through analytic continuation. }
\la{heatacfiex2}
\end{figure}

\section{The advected heat equation} \label{sec:advheat}

The advected heat equation on the half-line, 
\begin{subequations}\la{advheatibvp}
\begin{align}\la{advheat}
u_t&=u_{xx}+cu_x, &&\hspace*{-1.0in} x>0,\, t>0,\hspace*{1.0in}\\
u(x,0)&=u_0(x), &&\hspace*{-1.0in} x>0,\\
u(0,t)&=f_0(t), &&\hspace*{-1.0in} t>0,
\end{align}
\end{subequations}
\no is an interesting next example since its underlying operator is not self adjoint. It is also a first example whose solution is more conveniently obtained using Fokas' UTM than using classical techniques \cite{JC_fokas_book}. The solution is represented as
\beq
u(x,t) = I_0(x,t) + I_{f_0}(x,t),
\eeq
\no with
\begin{align}\la{I0heatadv1}
I_0(x,t) &=  \frac{1}{2\pi} \int_{-\infty}^\infty  e^{ikx-Wt}\hat u_0(k)\, dk-\frac{1}{2\pi} \int_{\partial \Omega}  e^{ikx-Wt} \hat u_0(-k+ic)\, dk,\\\la{Ifheatadv1}
I_{f_0}(x,t) &=  -\frac{1}{2\pi} \int_{\partial \Omega} (2ik+c)e^{ikx-Wt}  F_0(W,t) \, dk,
\end{align}
\no where $W=k^2-ikc$, and 
\beq
\hat u_0(k) = \int_0^\infty e^{-iky} u_0(y) \, dy, \qquad F_0(W,t) = \int_0^t e^{Ws} f_0(s) \, ds.
\eeq
\no Here $\Omega=\{k\in \C: \pi/4 < \Arg(k) < 3\pi/4 $ and $|k|>r\}$ for some $r>|c|$ (due to $W(ic)=0$), shown in Figure~\ref{fig:Omega}.


The initial condition contribution $I_0(x,t)$
%
%
%
is entire, as in Theorem~\ref{thm:IC}. The boundary contribution can be written as
\beq \la{advheat27} 
I_{f_0}(x,t) = -\frac{1}{2\pi}  \int_{\partial \Omega} dk \, (2ik+c)e^{ikx-Wt}\int_0^t  e^{Ws}f_0(s) \, ds= \frac{x}{2\sqrt{\pi}} \int_0^t  \frac{f_0(s)}{(t-s)^\frac32} e^{-\frac{(x+c(t-s))^2}{4(t-s)}}\, ds,
\eeq
\no which is analytic for $x>0$ (again, as in Theorem~\ref{thm:IC}), but discontinuous at $x=0$ since $I_{f_0}(0,t)=0$ but $\lim_{x\rightarrow 0^+}I_{f_0}(x,t)=f_0(t)$. Instead, we deform $\partial \Omega$ to $\gamma$ (see Figure~\ref{fig:Omega}), passing above $k=i|c|$ and integrate by parts $n+1$ times so that
\begin{align}\nonumber
I_{f_0} (x,t) &= -\frac{1}{2\pi} \int_{\gamma} dk \, (2ik+c)e^{ikx-Wt} \left[ \sum_{m=1}^{n+1}(-1)^m\frac{f_0^{(m-1)}(0)}{W^m} - \frac{(-1)^n}{W^{n+1}} \int_0^t  e^{Ws} f_0^{(n+1)}(s) \, ds \right] \\
&= \sum_{m=0}^{n+1} f_0^{(m-1)}(0)\phi_{m}(x,t) +  \int_0^t f_0^{(n+1)}(s) \phi_{n+1}(x,t-s) \, ds,
\end{align}
\no with
\beq
\phi_{m}(x,t) = - \frac{(-1)^m}{2\pi} \int_{\gamma}  \frac{(2ik+c)e^{ikx-Wt}}{W^m} \, dk.
\eeq
\no Switching the order of integration is allowed, due to absolute integrability. Taking $2n+q$ derivatives ($q=0$ or $1$), 
\beq
\frac{\partial^{2n+q} I_{f_0}}{\partial x^{2n+q}} = \sum_{m=1}^{n+1} f_0^{(m-1)}(0)\phi_{m}^{(2n+q)}(x,t) +  \int_0^t f_0^{(n+1)}(s) \phi_{n+1}^{(2n+q)}(x,t-s) \, ds .
\eeq
\no Using $k\mapsto k/\sqrt{t}$, 
\begin{align}\nonumber
\phi_{n+1}^{(2n+q)}(0,t) &= \frac{(-1)^n}{2\pi}  \int_{\gamma}  \frac{(ik)^{2n+q}(2ik+c)e^{-Wt}}{W^{n+1}} \, dk = \frac{i^q}{2\pi}  \int_{\gamma} \frac{k^{2n+q}(2ik+c)e^{-Wt}}{W^{n+1}} \, dk \\
&= \frac{i^q}{2\pi t^{\frac q2}}  \int_{\gamma} \frac{k^{2n+q}(2ik+c\sqrt{t})e^{-k^2+ikc\sqrt{t}}}{(k^2-ikc\sqrt{t})^{n+1}} \, dk 
= \mathcal{O}\left(t^{-\frac q2}\right).
\end{align}
\no Using the dominated convergence theorem, we get the coefficients of the Taylor series for $I_{f_0}(x,t)$ as ($n\in \mathbb{N}$),
\beq
\label{heat_adv_coeff_even}
a_{2n}(t) = \frac{1}{(2n)!} \left[\sum_{m=1}^{n+1}  f_0^{(m-1)}(0) \phi_{m}^{(2n)}(0,t) + \int_0^t f_0^{(n+1)}(s) \phi_{n+1}^{(2n)}(0,t-s)\, ds\right],
\eeq
\no and
\beq
\label{heat_adv_coeff_odd}
a_{2n+1}(t) = \frac{1}{(2n+1)!} \left[ \sum_{m=1}^{n+1} f_0^{(m-1)}(0)\phi_{m}^{(2n+1)}(0,t) + \int_0^t f_0^{(n+1)}(s) \phi_{n+1}^{(2n+1)}(0,t-s)\, ds \right].
\eeq
\no As before, we can use the UTM integral representation \rf{advheat27} to reduce the Taylor series to even or odd terms. With $a_0(t) = f_0(t)$, even is preferred, so that
\beq \la{advheatbext}
I_{f_0}^\text{ext}(x,t) = \case{1}{I_{f_0}(x,t), & x\geq 0, \\ \tilde f_0(x,t) - I_{f_0}(-x,t), & x<0,} ~~\text{where} ~~ \tilde f_0(x,t) = 2\sum_{n=0}^\infty a_{2n}(t) x^{2n}.
\eeq

\subsection{Boundary-to-Initial Map}

Define 
\beq
I_{f_0}^\text{ext}(x,t) = \sum_{n=0}^\infty a_{n}(t)x^{n},
\eeq
\no where $a_{n}(t)$ is defined above in \rf{heat_adv_coeff_even} and \rf{heat_adv_coeff_odd}. Then
\beq \label{heatadvac}
w(x,t) = u_\text{ac}(x,t) = I_0(x,t) + I_{f_0}^{\text{ext}}(x,t),
\eeq
\no is a whole-line solution to the heat equation with advection with initial condition 
\beq
w_0(x) = I_0(x,0) + I_{f_0}^\text{ext}(x,0).
\eeq

Since
\begin{align}
I_0(x,t) 
&=  \frac{1}{2\pi} \int_0^\infty dy\, u_0(y) \int_{-\infty}^\infty  e^{ik(x-y)-Wt}\, dk-\frac{1}{2\pi} \int_{0}^\infty dy \, u_0(y)e^{cy} \int_{\partial \Omega}  e^{ik(x+y)-Wt} \, dk, 
\end{align}
\no we have
\beq
I_0(x,0) = \case{1}{u_0(x), & x\geq 0, \\ -e^{-cx}u_0(-x), & x< 0.}
\eeq
\no For small $t>0$,
\begin{align}\nonumber
\tilde f_0(x,t) 
&= 2\sum_{n=0}^\infty \frac{x^{2n}}{(2n)!}  \left[\sum_{m=1}^{n+1}  f_0^{(m-1)}(0) \phi_{m}^{(2n)}(0,t) + \int_0^t f_0^{(n+1)}(s) \phi_{n+1}^{(2n)}(0,t-s)\, ds\right] \\
\nonumber
&\sim 2\sum_{m=0}^{\infty}f_0^{(m)}(0)\sum_{n=m}^\infty \frac{x^{2n}}{(2n)!}  \phi_{m+1}^{(2n)}(0,t) \\
&= \frac{1}{\pi}\sum_{m=0}^{\infty}f_0^{(m)}(0)\sum_{n=m}^\infty \frac{(-1)^{m+n}x^{2n}}{(2n)!}  \int_{\gamma}  \frac{(2ik+c)k^{2n-m-1}e^{-(k^2-ikc)t}}{(k-ic)^{m+1}} \, dk.
\end{align}
\no Using $k=\kappa+ic/2$, 
\begin{align}
\tilde f_0(x,t) &\sim \frac{1}{\pi}\sum_{m=0}^{\infty}f_0^{(m)}(0)\sum_{n=m}^\infty \frac{(-1)^{m+n}x^{2n}}{(2n)!}  \int_{\gamma'}  \frac{2i\kappa \left(\kappa+\frac{ic}{2}\right)^{2n-m-1}e^{-\left(\kappa^2+\frac{c^2}{4}\right)t}}{\left(\kappa+\frac{ic}{2}\right)^{m+1}} \, d\kappa,
\end{align}
\no where $\gamma'$ is $\gamma$ shifted by $ic/2$. If $n=m=0$, the integrand has a simple pole in the UHP at $k=i|c|/2$. Otherwise it has a pole in the upper half plane if $c>0$ at $k=ic/2$ of order $m+1$. Deforming down to the real axis, around these poles, and using the residue theorem gives 
\beq
\tilde f_0(x,t) \sim \tilde f_0^{(i)}(x,t) + \tilde f_0^{(r)}(x,t),
\eeq
\no with integral part
\begin{align}\nonumber
\tilde f_0^{(i)}(x,t) &= \frac{1}{\pi}\sum_{m=0}^{\infty}f_0^{(m)}(0)\sum_{n=m}^\infty \frac{(-1)^{m+n}x^{2n}}{(2n)!}  \int_{-\infty}^\infty  \frac{2i\kappa \left(\kappa+\frac{ic}{2}\right)^{2n-m-1}e^{-\left(\kappa^2+\frac{c^2}{4}\right)t}}{\left(\kappa+\frac{ic}{2}\right)^{m+1}} \, d\kappa\\
&\to\frac{1}{\pi}\sum_{m=0}^{\infty} \frac{x^{2m}}{(2m)!}f_0^{(m)}(0) \int_{-\infty}^\infty  \frac{2i\kappa \left(\kappa+\frac{ic}{2}\right)^{m-1}
}{\left(\kappa+\frac{ic}{2}\right)^{m+1}} {}_1F_{2}\left(\begin{array}{c} 1 \\ m+\frac12, \, m+1 \end{array};\, -\frac{\left(\kappa + \frac{ic}{2}\right)^2x^2}{4}\right) \, d\kappa, 
\end{align} 
\no after switching the sum and the integral. Here ${}_1F_{2}$ denotes the generalized hypergeometric function~\cite{dlmf}. The residue part equals
\beq
\tilde f_0^{(r)}(x,t) = 2\sum_{m=0}^{\infty}f_0^{(m)}(0)\sum_{n=m}^\infty \frac{(-1)^{m+n}x^{2n}}{(2n)!} \underset{\kappa=i|c|/2}{\text{Res}}\left( \frac{2\kappa \left(\kappa+\frac{ic}{2}\right)^{2n-m-1}e^{-\left(\kappa^2+\frac{c^2}{4}\right)t}}{\left(\kappa+\frac{ic}{2}\right)^{m+1}}\right).
\eeq
\no For $(n,m)=(0,0)$, 
\begin{align}
\underset{\kappa=i|c|/2}{\text{Res}}\left( \frac{2\kappa e^{-\left(\kappa^2+\frac{c^2}{4}\right)t}}{\kappa^2+\frac{c^2}{4}}\right) &= 1,
\end{align}
\no while for $(n,m)\neq (0,0)$,
\begin{align}\nonumber 
\frac{2\kappa\left(\kappa+\frac{ic}{2}\right)^{2n-m-1}}{\left(\kappa-\frac{ic}{2}\right)^{m+1}} &= \frac{2\kappa(\kappa-\frac{ic}{2}+ic)^{2n-m-1}}{\left(\kappa-\frac{ic}2\right)^{m+1}} = \sum_{j=0}^{2n-m-1} \binom{2n-m-1}{j} \frac{2\kappa(ic)^{2n-m-j-1}}{\left(\kappa-\frac{ic}2\right)^{m-j+1}}\\
&= \sum_{j=0}^{2n-m-1} \binom{2n-m-1}{j} \left[ \frac{2(ic)^{2n-m-j-1}}{\left(\kappa-\frac{ic}2\right)^{m-j}} +  \frac{(ic)^{2n-m-j}}{\left(\kappa-\frac{ic}2\right)^{m-j+1}}\right],
\end{align}
\no so that
\begin{align}\nonumber 
\underset{\kappa=i|c|/2}{\text{Res}}\left( \frac{2\kappa \left(\kappa+\frac{ic}{2}\right)^{2n-m-1}e^{-\left(\kappa^2+\frac{c^2}{4}\right)t}}{\left(\kappa+\frac{ic}{2}\right)^{m+1}}\right) &\sim  2\binom{2n-m-1}{m-1} (ic)^{2n-2m} + \binom{2n-m-1}{m}(ic)^{2n-2m}\\
&= \frac{(-1)^{n-m}(2n)(2n-m-1)!c^{2n-2m}}{m!(2n-2m)!}.
\end{align}
\no Then,
\beq
\tilde f_0^{(r)}(x,t) \sim \case{2}{2f_0(0), & c<0, \\ \displaystyle 2f_0(0)+2f_0(0)\sum_{n=1}^\infty \frac{c^{2n}x^{2n}}{(2n)!}+2\sum_{m=1}^{\infty}f_0^{(m)}(0)\sum_{n=m}^\infty \frac{(2n-m-1)!c^{2n-2m}x^{2n}}{(2n-1)!m!(2n-2m)!}, & c>0,}
\eeq
\no and
\beq
\tilde f_0^{(r)}(x,0) = \case{2}{2f_0(0), & c<0, \\ \displaystyle 2f_0(0)\cosh(cx)+ 2\sqrt{\pi} \cosh\left(\frac{cx}{2}\right)\sum_{m=1}^{\infty}\frac{|x|^{m+\frac12}}{c^{m-\frac12}m!}I_{m-\frac12}\left(\frac{cx}{2}\right)f_0^{(m)}(0), & c>0,}\\
\eeq
\no which gives the boundary-to-initial map,
\beq
\label{eqn:w0adv}
w_0(x) = \case{1}{u_0(x), & x\geq 0, \\ -e^{-cx}u_0(-x) + \tilde{f}_0(x,0), & x<0.}
\eeq
\no It can be shown that these integrals and sums are convergent and that in the limit as $c\to0$, $\tilde f_0(x,t)$ limits to \rf{f0tilde}.





\subsection{Examples}

We start with the whole line solution $u_\R(x+ct+1,t)$ \rf{wholeline1} with the corresponding boundary and initial conditions, equating $c=\pm 1$. The boundary integral is discontinuous, but the analytic continuation recovers the exact whole-line solution as shown in Figures~\ref{fig:advheat1} and~\ref{fig:advheat2}. Considering the boundary condition $f_0(t) = te^{-t}$, we find the analytic continuation of the solution shown in Figure~\ref{fig:advheatalt}.

\begin{figure}[tb]
\begin{center}
\def \sc {0.475}
\begin{tabular}{cc}
\includegraphics[scale=\sc]{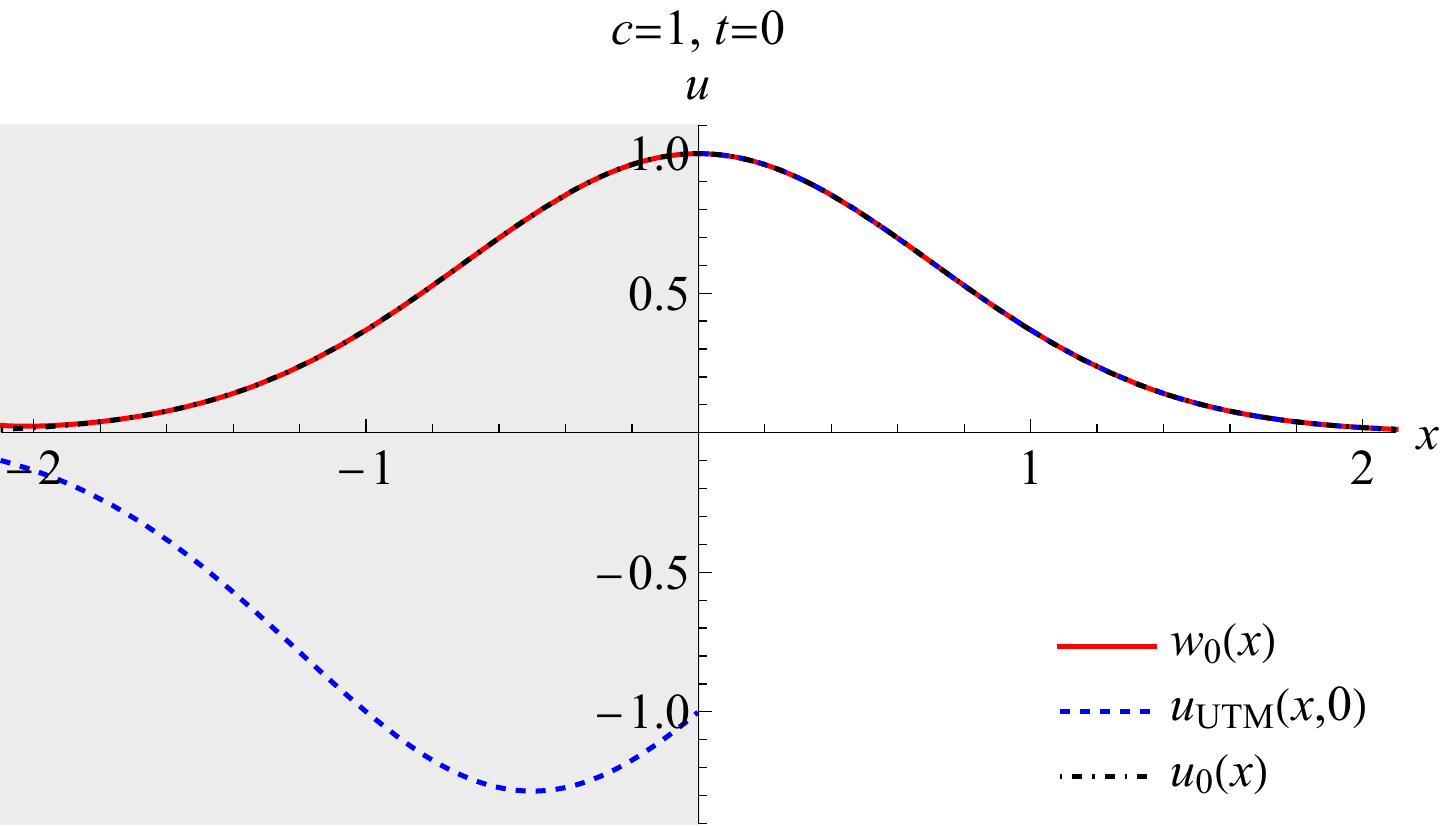} & \includegraphics[scale=\sc]{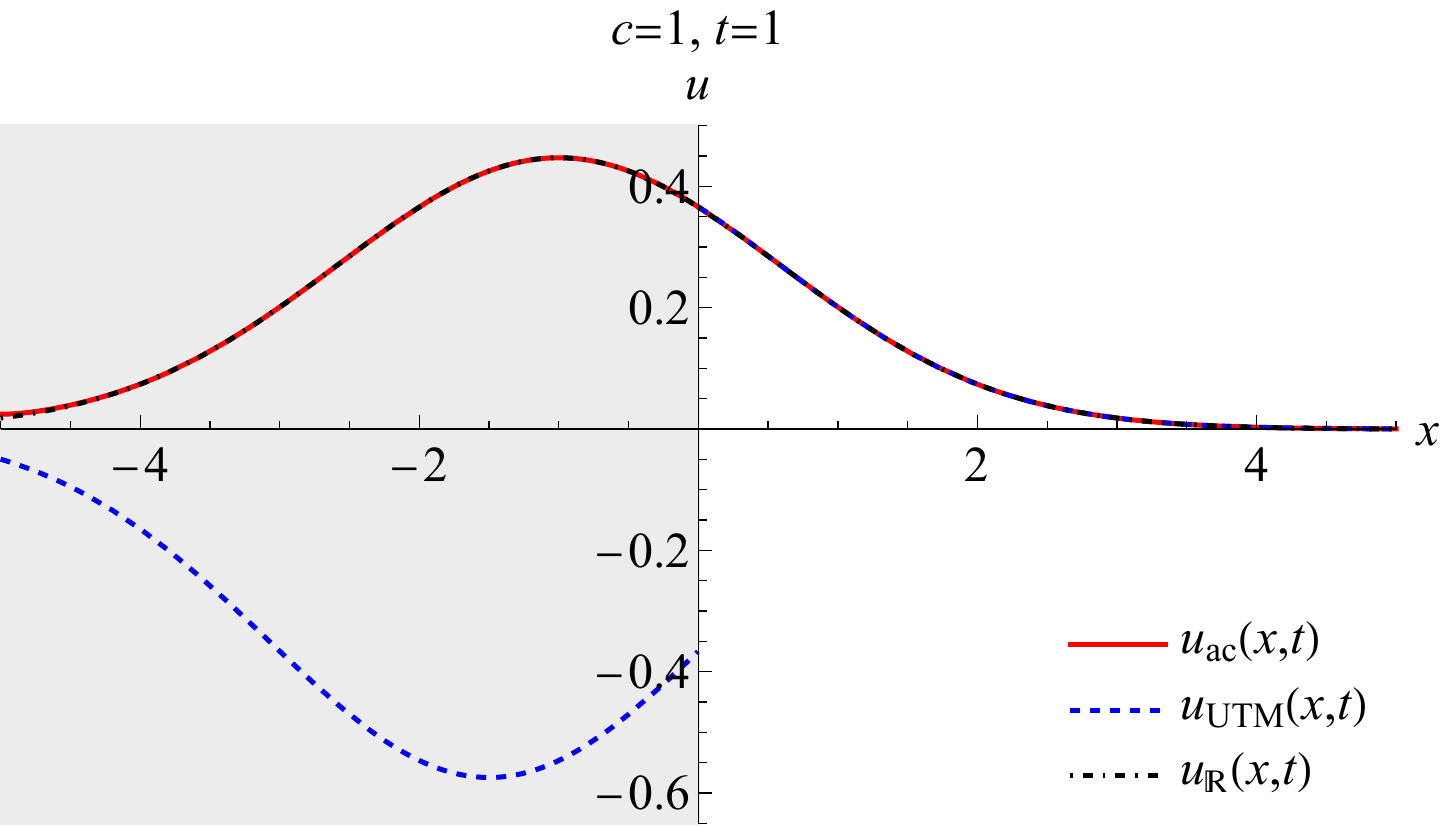}\\
(a) & (b)
\end{tabular}
\end{center}
\caption{With $c=1$: (a) The initial condition $w_0(x)$ \rf{eqn:w0adv}, leading to the (analytic) solution on the right, shown together with the UTM solution at $t=0$ and the whole-line initial condition $u_0(x) = u_\R(x+1,0)$. (b) The solution $u_{\text{ac}}(x,t)$ \rf{heatadvac} at $t=1$ obtained through analytic continuation.}
\la{fig:advheat1}
\end{figure}

\begin{figure}[tb]
\begin{center}
\def \sc {0.475}
\begin{tabular}{cc}
\includegraphics[scale=\sc]{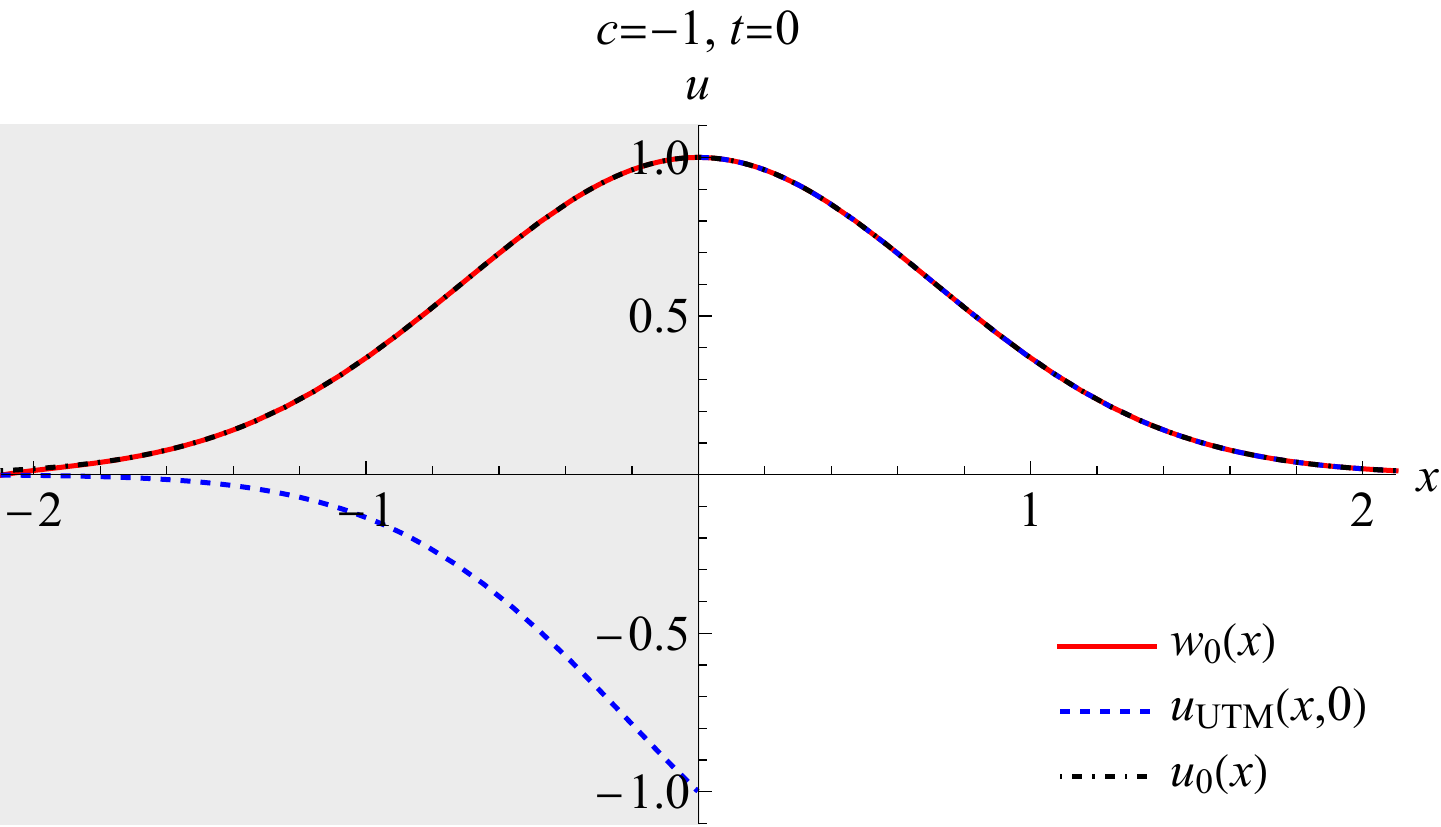} & \includegraphics[scale=\sc]{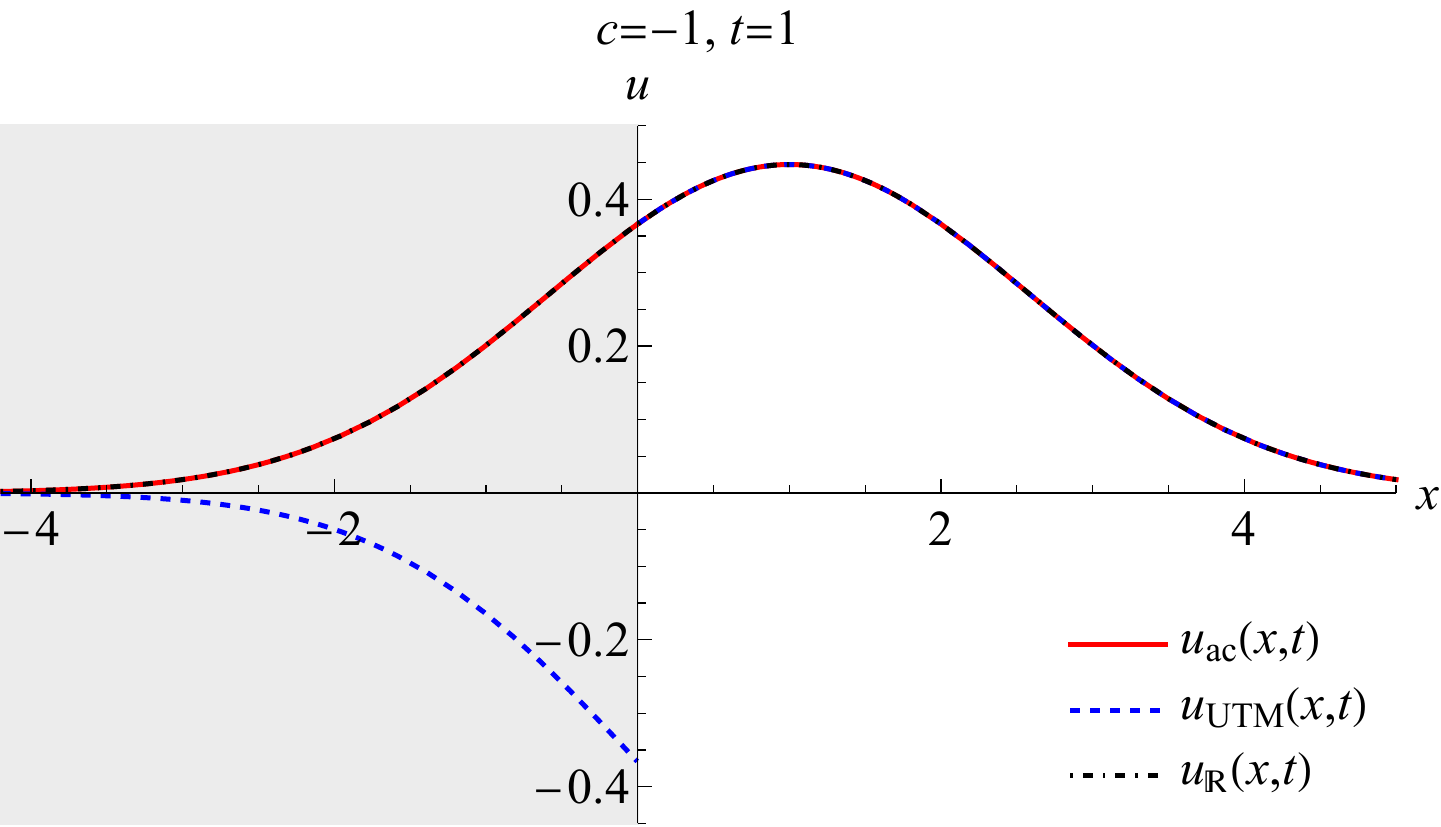}\\
(a) & (b)
\end{tabular}
\end{center}
\caption{with $c=-1$: (a) The initial condition, $w_0(x)$, \rf{eqn:w0adv}, leading to the (analytic) solution on the left, shown with the the UTM solution at $t=0$ and the whole-line initial condition $u_0(x)=u_\R(x+1,0)$. (b) The solution $u_{\text{ac}}(x,t)$ \rf{heatadvac} at $t=1$ obtained using analytic continuation.}
\la{fig:advheat2}
\end{figure}

\begin{figure}[tb]
\begin{center}
\def \sc {0.475}
\begin{tabular}{cc}
\includegraphics[scale=\sc]{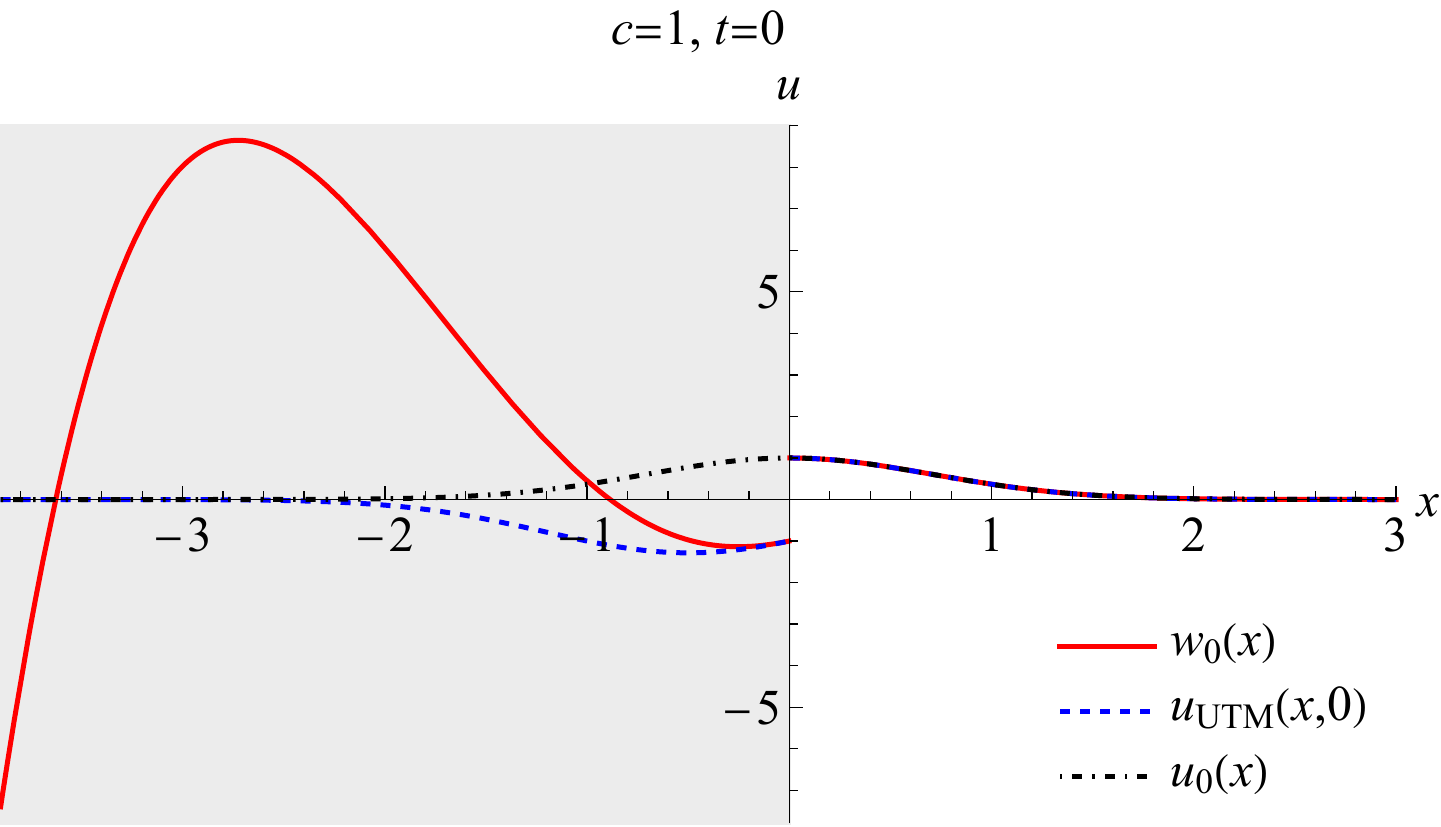} & \includegraphics[scale=\sc]{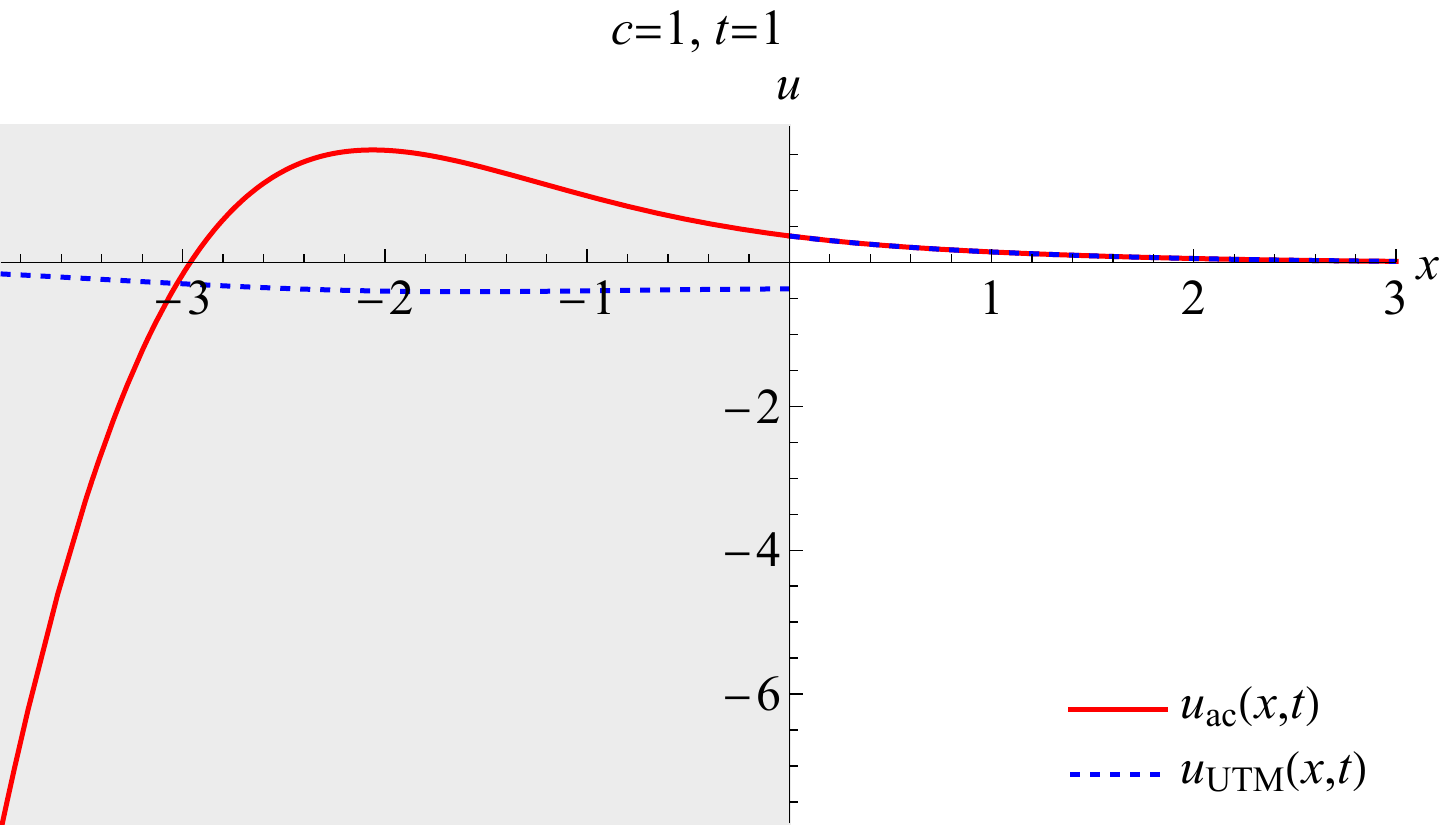}\\
(a) & (b)
\end{tabular}
\end{center}
\caption{(a) The (discontinuous) initial condition $w_0(x)$ \rf{eqn:w0adv} leading to the (analytic) solution on the right, shown with the analytic continuation of $u_0(x)=u_\R(x+1,0)$. (b) The solution $u_{\text{ac}}(x,t)$ \rf{heatadvac} at $t=1$ obtained through analytic continuation. }
\la{fig:advheatalt}
\end{figure}

\section{The linear KdV equation, 1 boundary condition} \label{sec:kdv1}

\def \contour {{\partial\Omega}}

Consider the linear Korteweg--de Vries (KdV) equation on the half-line with Dirichlet boundary conditions,
\begin{subequations}\la{kdvprob1}
\begin{align}\la{kdv1ibvp}
u_t+u_{xxx} &= 0, &&\hspace*{-1.0in} x>0, \, t>0,\hspace*{1.0in}\\
u(x,0)&=u_0(x), &&\hspace*{-1.0in} x>0,\\
u(0,t)&=f_0(t), &&\hspace*{-1.0in} t>0,
\end{align}
\end{subequations}
\no which requires only one boundary condition \cite{JC_fokas_book}, unlike the similar equation in the next section. Its solution is written as
\beq 
u(x,t) = I_0(x,t) + I_{f_0}(x,t),
\eeq
\no with
\begin{align}\la{kdv1I0}
I_0(x,t) &= \frac{1}{2\pi} \int_{-\infty}^\infty  e^{ikx+ik^3t}\hat u_0(k) \, dk +\frac{1}{2\pi} \int_{\contour}  e^{ikx+ik^3t}\left(\alpha \hat  u_0(\alpha k) + \alpha^2 \hat u_0(\alpha^2k)\right)  dk, \\
I_{f_0}(x,t) &= -\frac{1}{2\pi} \int_{\contour}  3k^2e^{ikx+ik^3t}F_0(-ik^3,t)\, dk, 
\end{align}
\no and 
\beq
\hat u_0(k) = \int_0^\infty e^{-iky} u_0(y) \,dy, \qquad
F_0(-ik^3,t) = \int_0^t e^{-ik^3s} f_0(s) \,ds,
\eeq 
\no and $\alpha = \exp({2i\pi}/{3}).$ The region $\Omega =\{k\in \C:\pi/3<\Arg(k)<2\pi/3\}$ is shown in Figure~\ref{fig:Omega_kdv}.
Theorem~\ref{thm:IC_KdV} shows that after some contour deformations, $I_0(x,t)$ defines an entire function of~$x$.  

\figcl{\imsize}{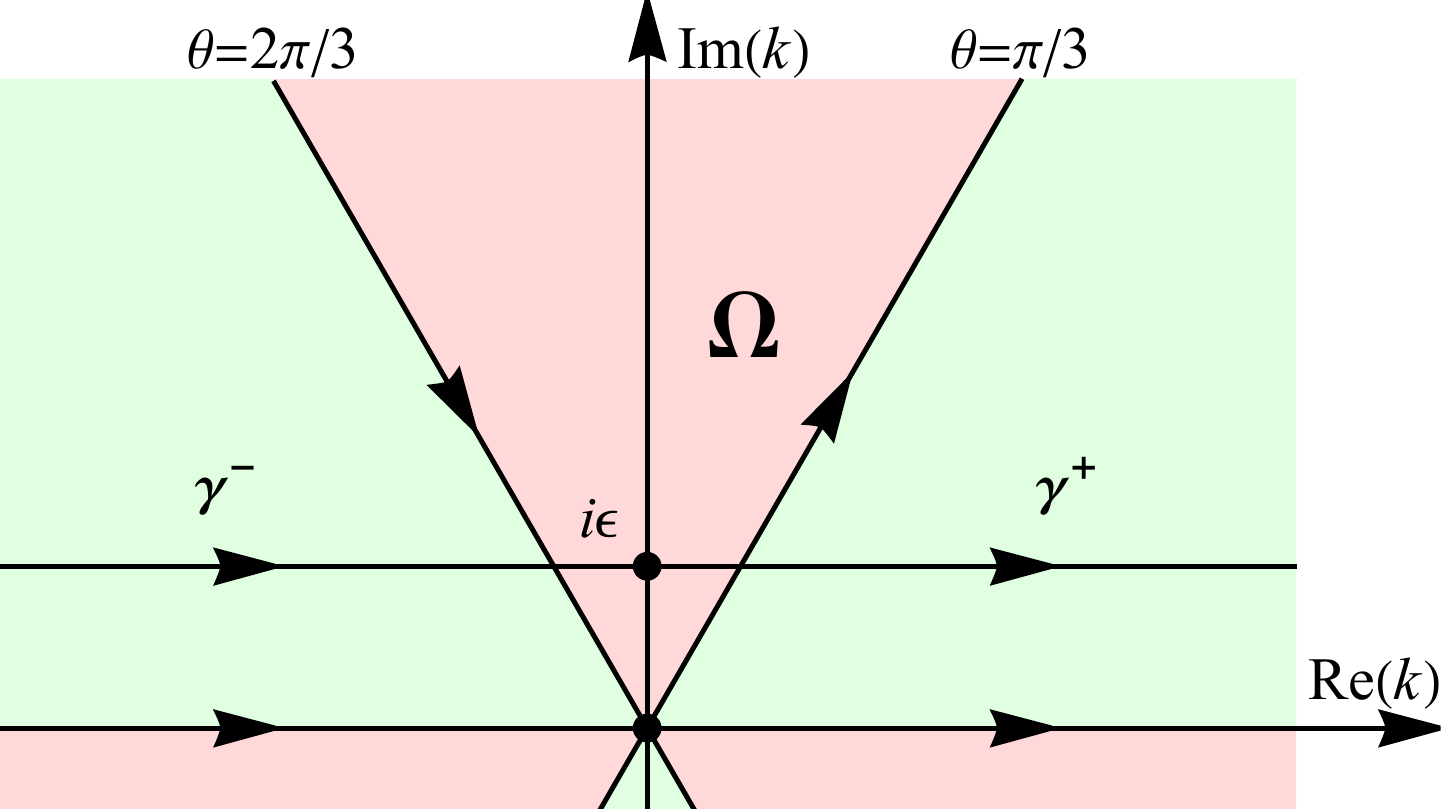}{The region $\Omega$ for KdV and the contour $\gamma=\gamma^-\cup \gamma^+$.}{Omega_kdv}

\begin{theorem}\label{thm:IC_KdV}

If there exists an $\epsilon>0$ such that 
\[
\left\|e^{\epsilon y}u_0(y)\right\|_1 = \int_0^\infty |u_0(y)|e^{\epsilon y} \, dy < \infty,
\]
\no then $I_0(x,t)$, is entire in $x$ for each $t>0$.\footnote{In fact, less decay is needed. If $u_0(y) =\mathcal O( e^{-z^p})$ for $p>1/2$, then $I_0(x,t)$ is analytic for $x\in \R$.} If $f_0\in L^1(0,T)$ for some $T>0$, then $I_{f_0}(x,t)$ is analytic for $|\Im(x)|<\sqrt{3}\Re(x)$ for $0<t<T$.

\end{theorem}

\begin{proof}

We deform the path for the first integral in $I_0(x,t)$ up to $\gamma=\gamma^-\cup \gamma^+$, see Figure~\ref{fig:Omega_kdv}, using the assumption on $u_0(y)$, still assuming $x>0$. Next, an integral over any closed contour $\Gamma$ in the complex $x$-plane is zero by Cauchy's theorem, where we can use Fubini's theorem to swap the order of integration since
\begin{align*}
\oint_\Gamma |dx| \int_{\gamma} \left|e^{ikx+ik^3t}\hat u_0(k)\right| \, |dk| &\leq \ell(\Gamma) \left\|e^{\epsilon y} u_0(y)\right\|_1 \max_{x\in\Gamma} \int_{-\infty}^\infty e^{\epsilon |x|+\kappa |x|+(\epsilon^3-3\epsilon\kappa^2)t} \, d\kappa \\
&= \ell(\Gamma) \left\|e^{\epsilon y} u_0(y)\right\|_1  \sqrt{\frac{\pi}{3\epsilon t}} \max_{x\in \Gamma} e^{\epsilon|x|+\epsilon^3t+\frac{|x|^2}{12\epsilon t}} < \infty,
\end{align*}
\no where we parametrized the contour $\gamma$ with $k=\kappa+i\epsilon$. It follows that the first integral in $I_0(x,t)$, after the contour deformation, is entire by Morera's theorem. 

For the second integral, we break the contour $\partial\Omega$ into two parts: the left part $\partial\Omega^-$ and the right part $\partial \Omega^+$. We deform $\partial\Omega^-$ down to $\gamma^-$, so that the corresponding integral over $\Gamma$ is zero by Cauchy's theorem, since
\begin{align*}
\oint_\Gamma |dx| \int_{\gamma^-} \left| e^{ikx+ik^3t} \hat u_0(\alpha k) \right|\, |dk| &\leq \ell(\Gamma)\|u_0\|_1 \max_{x\in \Gamma} \int_0^\infty e^{\epsilon|x|+\kappa |x|+(\epsilon^3 - 3 \epsilon \kappa^2)t} \, d\kappa \\
&\leq \ell(\Gamma)\|u_0\|_1 \sqrt{\frac{\pi}{3\epsilon t}} \max_{x\in \Gamma}  e^{\epsilon|x|+\epsilon^3t+\frac{|x|^2}{3\epsilon t}} < \infty,
\end{align*}
\no where we used $1+\erf(y)\leq 2$ for $y\in\R$.

For the $\partial\Omega^+$ part, we use the transformation $\kappa = \alpha k$ and deform up to $\gamma^-$, so that
\[
\int_{\partial\Omega^+} e^{ikx+ik^3t} \alpha \hat u_0(\alpha k) \, dk = \int_{\R^-} e^{i\alpha^2\kappa x+i\kappa^3t} \hat u_0(\kappa) \, dk= \int_{\gamma^-} e^{i\alpha^2\kappa x+i\kappa^3t} \hat u_0(\kappa) \, dk.
\]
\no If we integrate this over $\Gamma$, we get zero by Cauchy's theorem, where we can use Fubini's theorem since
\begin{align*}
\oint_\Gamma |dx| \int_{\gamma^-} \left|e^{i\alpha^2\kappa x+i\kappa^3t} \hat u_0(\kappa)\right| \, |dk| &\leq \ell(\Gamma) \left\|e^{\epsilon y}u_0(y)\right\|_1 \max_{x\in\Gamma} \int_0^\infty e^{\sqrt{3}\epsilon|x|+\sqrt{3}\kappa |x|+\left(\epsilon^3-3\epsilon \kappa^2\right) t} \, d\kappa \\
&\leq \ell(\Gamma) \left\|e^{\epsilon y} u_0(y) \right\|_1 \sqrt{\frac{\pi}{3\epsilon t}} \max_{x\in\Gamma} e^{\sqrt{3}\epsilon|x| +\epsilon^2t+ \frac{|x|^2}{4\epsilon t}} < \infty.
\end{align*}
\no The same holds for the third term in $I_0(x,t)$, thus $I_0(x,t)$ is entire.
For $I_{f_0}(x,t)$, an integral over a closed contour $\Gamma$ in the region $|\Im(x)|<\sqrt{3}\Re(x)$ is zero if
\begin{align*}
\oint |dx| \int_{\partial\Omega} \left|3k^2 e^{ikx+ik^3t}F_0(-ik^3,t)\right| \, |dk| &\leq \ell(\Gamma) \|f_0\|_1 \max_{x\in\Gamma} \int_{\partial\Omega} 3|k|^2 |e^{ikx}| \, |dk| \\
&\leq \ell(\Gamma) \|f_0\|_1 \max_{x\in\Gamma} \frac{288 \sqrt{3}|x|^3}{(3\Re(x)^2-\Im(x)^2)^3} < \infty,
\end{align*}
\no so that $I_{f_0}(x,t)$ is an analytic function of $x$ for $|\Im(x)| < \sqrt{3}\Re(x)$.

\end{proof}

If we swap the order of integration for $I_{f_0}(x,t)$ and integrate the $k$-integral, we find
\beq
I_{f_0}(x,t) = \frac{x}{\sqrt[3]{3}} \int_0^t \frac{f_0(s)}{(t-s)^{\frac43}} \Ai\left(\frac{x}{\sqrt[3]{3(t-s)}}\right) \, ds,
\eeq
\no where $\Ai(z)$ denotes the Airy function, see \cite{dlmf}. From this, it is clear this is not defined for $x<0$. Instead, we deform to a contour $\gamma$ lying under $\Omega$ (ensuring $\gamma$ passes above the origin) and above $\R$, see Figure~\ref{fig:Omega_kdv}:
\beq 
I_{f_0}(x,t) = -\frac{1}{2\pi} \int_{\gamma} dk \, 3k^2 e^{ikx+ik^3t} \int_0^t  e^{-ik^3s} f_0(s)\, ds.
\eeq
\no For the $x$-derivatives of order $3n$, we use
\begin{align} 
\frac{\partial^n}{\partial s^n}\delta(s-t) &= -\frac{(-i)^n}{2\pi} \int_{\partial\Omega}  3k^{3n+2} e^{ik^3(t-s)}\, dk,
\end{align}
\no leading to
\beq 
a_{3n}(t) = \frac{1}{(3n)!} \int_0^t f_0(s) \delta^{(n)}(s-t) \, ds = \frac{(-1)^nf_0^{(n)}(t)}{(3n)!} . 
\eeq

For the other derivatives, we deform to $\gamma$ and integrate by parts $n$ times to find
\begin{align}\nonumber
I_{f_0}(x,t) 
&= -\frac{1}{2\pi} \int_{\gamma}dk \, 3k^2 e^{ikx+ik^3t} \left[\sum_{m=1}^n \frac{f_0^{(m-1)}(0)}{(ik^3)^m} + \frac{1}{(ik^3)^n}\int_0^t e^{-ik^3s} f_0^{(n)}(s) ds\right] \\
&= \sum_{m=1}^n f_0^{(m-1)}(0) \phi_m(x,t) + \int_0^t f_0^{(n)}(s) \phi_n(x,t-s) ds,
\end{align}
\no with 
\beq 
\phi_m(x,t) = -\frac{1}{2\pi} \int_{\gamma} \frac{3k^2 e^{ikx+ik^3t}}{(ik^3)^m} dk.
\eeq
\no Taking $3n-q$ derivatives ($q=1,2$),
\beq 
\frac{\partial^{3n-q}I_{f_0}}{\partial x^{3n-q}}  = \sum_{m=1}^n f_0^{(m-1)}(0) \phi_m^{(3n-q)}(x,t) + \int_0^t f_0^{(n)}(s) \phi_n^{(3n-q)}(x,t-s) ds.
\eeq
\no The contour deformation to $\gamma$ enables the differentiating under the integral sign and the use of the dominated convergence theorem. With
\begin{align}\nonumber 
\phi_m^{(3n-q)}(0,t) &= -\frac{1}{2\pi} \int_{\gamma} dk \, \frac{3(ik)^{3n-q}k^2e^{ik^3t}}{(ik^3)^m} 
= \frac{3(-1)^{m-q}}{\pi}\sin\left(\frac{q\pi}{3} \right)\int_{0}^\infty d\rho \, \rho^{3(n-m)-q+2}e^{-\rho^3t}  \\
&= \frac{(-1)^{m-q}}{\pi t^{n-m+1-\frac{q}{3}}}\sin\left(\frac{q\pi}{3} \right) \Gamma\left(n-m+1-\frac{q}{3}\right).
\end{align} 
\no It follows that the coefficients of the Taylor series are 
\beq
a_{3n-2}(t) = \frac{\sqrt{3}\Gamma\left(\frac{1}{3}\right)}{2\pi (3n-2)!} \left[\sum_{m=1}^n \frac{(-1)^{m}\Gamma\left(n-m+\frac{1}{3}\right)}{\Gamma\left(\frac{1}{3}\right) t^{n-m+\frac{1}{3}}} f_0^{(m-1)}(0) + (-1)^{n}\int_0^t ds \, \frac{f_0^{(n)}(s)}{(t-s)^{\frac{1}{3}}}  \right], 
\eeq
\no and
\beq 
a_{3n-1}(t) = -\frac{\sqrt{3}\Gamma\left(\frac{2}{3}\right)}{2\pi(3n-1)!} \left[ \sum_{m=1}^n  \frac{(-1)^{m}\Gamma\left(n-m+\frac{2}{3}\right)}{ \Gamma\left(\frac{2}{3}\right)t^{n-m+\frac{2}{3}}} f_0^{(m-1)}(0)+ (-1)^n\int_0^t ds \, \frac{f_0^{(n)}(s)}{(t-s)^{\frac{2}{3}}}  \right].
\eeq
\no For $x<0$, we may write
\beq
\label{eqn:f0tildekdv1}
I_{f_0}(x,t) = 2\sum_{n=0}^\infty a_{2n}(t)x^{2n} - I_{f_0}(-x,t) = \tilde f_0(x,t) - I_{f_0}(-x,t),
\eeq
\no so that
\beq\la{advheatIf0ext}
I_{f_0}^\text{ext} (x,t) = \case{1}{I_{f_0}(x,t), &  x>0, \\ \tilde f_0(x,t) - I_{f_0}(-x,t), & x< 0, }
\eeq
\no and we may reduce the number of coefficients to be computed. As for the other monomial dispersion relations, the Taylor series may be written more compactly as
\beq \label{If0series}
I_{f_0}^\text{ext}(x,t) =\sum_{n=0}^\infty a_n(t) x^n = \sum_{n=0}^\infty \frac{(-1)^n}{n!} f_0^{\left(\frac n3\right)}\!(t)\,x^{n} ,
\eeq
\no using Riemann--Liouville fractional derivatives. 

\subsection{Boundary-to-Initial Map}

The function $I_{f_0}^{\text{ext}}(x,t)$, defined by \rf{If0series},
%
%
is easily seen to be a solution to the linear KdV equation, \rf{kdv1ibvp}.
We have
\beq
\label{eqn:If0extkdv1}
I_{f_0}^\text{ext}(x,t) = \sum_{n=0}^\infty a_{3n}(t)x^{3n} +\sum_{n=1}^\infty a_{3n-1}(t)x^{3n-1}+\sum_{n=1}^\infty a_{3n-2}(t)x^{3n-2},
\eeq
\no and
\beq
\sum_{n=0}^\infty a_{3n}(t)x^{3n} = \sum_{n=0}^\infty \frac{(-1)^n x^{3n}}{(3n)!}f_0^{(n)}(t) \to \sum_{n=0}^\infty \frac{(-1)^nx^{3n}}{(3n)!} f_0^{(n)}(0),
\eeq
\no as $t\to0^+$. For the other terms, 
\beq
\sum_{n=1}^\infty a_{3n-2}(t)x^{3n-2} \sim \frac{\sqrt{3}}{2\pi}\sum_{n=1}^\infty\frac{x^{3n-2}}{(3n-2)!} \sum_{m=1}^{n}  \frac{(-1)^{m}\Gamma\left(n-m+\frac13\right)}{t^{n-m+\frac13}}f_0^{(m-1)}(0),
\eeq
\no since the integral terms approach zero as $t\to 0^+$. We switch sums, so that
\begin{align}\nonumber
\sum_{n=1}^\infty a_{3n-2}(t)x^{3n-2} 
&\sim \frac{\sqrt{3}}{2\pi} \sum_{m=1}^{\infty}(-1)^{m}f_0^{(m-1)}(0) \sum_{n=m}^\infty\frac{x^{3n-2}}{(3n-2)!} \frac{\Gamma\left(n-m+\frac13\right)}{t^{n-m+\frac13}}\\
&= \frac{\sqrt{3}\Gamma\left(\frac13\right)}{2\pi t^{\frac13}} \sum_{m=1}^{\infty}\frac{(-1)^{m}x^{3m-2}}{(3m-2)!}f_0^{(m-1)}(0){}_2F_3\left(\begin{array}{c} \frac13, \,  1 \\ m-\frac13,\, m,\, m+\frac13\end{array}; \, \frac{x^3}{27t}\right).
\end{align}
\no Since \cite{dlmf}
\begin{align}
\frac{\sqrt{3}\Gamma\left(\frac13\right)}{2\pi t^{\frac13}}{}_2F_3\left(\begin{array}{c} \frac13, \,1 \\ m-\frac13,\,m,\,m+\frac13 \end{array};\,\frac{x^3}{27t}\right)
&\sim -\frac{\left(3m-2\right)}{x}, ~~~\mbox{as}~t\to 0^+,~~\mbox{for}~x<0,  
\end{align}
\no we obtain  
\beq
\sum_{n=1}^\infty a_{3n-2}(t)x^{3n-2} \to \sum_{m=0}^{\infty}\frac{(-1)^{m}x^{3m}}{(3m)!}f_0^{(m)}(0), ~~~\mbox{as}~t\to 0^+,~~\mbox{for}~x<0. 
\eeq
\no Similarly, 
\beq
\sum_{n=1}^\infty a_{3n-1}(t) x^{3n-1} \to  \sum_{m=0}^{\infty}\frac{(-1)^{m}x^{3m}}{(3m)!}f_0^{(m)}(0),~~~\mbox{as}~t\to 0^+,~~\mbox{for}~x<0.
\eeq
\no Combining \rf{eqn:f0tildekdv1} and \rf{eqn:If0extkdv1}, 
\beq
I_{f_0}(x,t) \to \tilde f(x,0) = 3\sum_{n=0}^\infty \frac{(-1)^n x^{3n}}{(3n)!} f_0^{(n)}(0), ~~~\mbox{as}~t\to 0^+,~~\mbox{for}~x<0.
\eeq

For $I_0(x,t)$, we have
\beq 
\frac{1}{2\pi} \int_{-\infty}^\infty  e^{ikx+ik^3t}\hat u_0(k) \, dk \to \case{1}{u_0(x), & x>0, \\ 0, & x<0,}
\eeq
\no and for $m=1,2$,
\begin{align}\nonumber
\frac{\alpha^m}{2\pi} \int_0^\infty dy \, u_0(y) \int_{\gamma} dk \, e^{ik(x-\alpha^m y)+ik^3t} &= \frac{1}{2\pi} \int_{\alpha^{m}\R^+} dy \, u_0(\alpha^{-m}y) \int_{\gamma} dk \, e^{ik(x- y)+ik^3t} \\
&\to \case{1}{0, & x>0, \\ -u_0(\alpha^{-m}x), & x<0,}
\end{align}
\no as $t\to 0^+$, so that
\beq 
\label{eqn:w0kdv1}
w_0(x) = \case{1}{ u_0(x), & x>0, \\ \tilde f_0(x,0) -u_0(\alpha x) - u_0(\alpha^2x), & x<0,}
\eeq
\no provides the boundary-to-initial map. \\

\subsection{Examples}

Our first example uses the whole-line solution,
\beq
\label{eqn:uRkdv}
u_{\R}(x,t) = 2e^{-(x+2t)}\cos(x-2t),
\eeq
\no restricted to $x>0$, with $u_0(x)=u_{\R}(x,0)$ and $f_0(t)=u_{\R}(0,t)$. The solution obtained using UTM is no longer defined for $x<0$, but the analytic continuation recovers the exact solution on the whole line. The results are shown in Figure~\ref{fig:kdv1}.





\begin{figure}[tb]
\begin{center}
\def \sc {0.475}
\begin{tabular}{cc}
\includegraphics[scale=\sc]{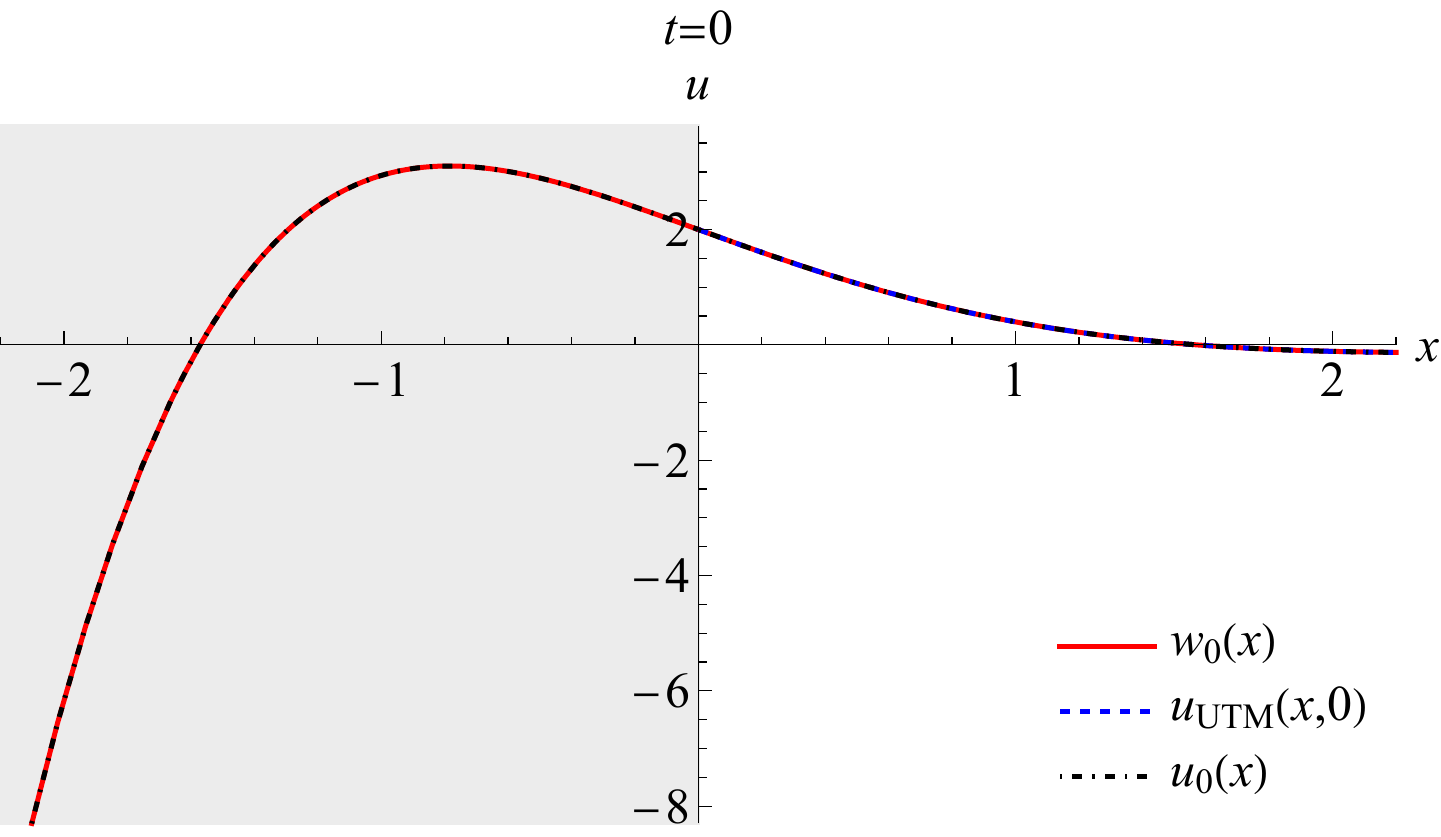} & \includegraphics[scale=\sc]{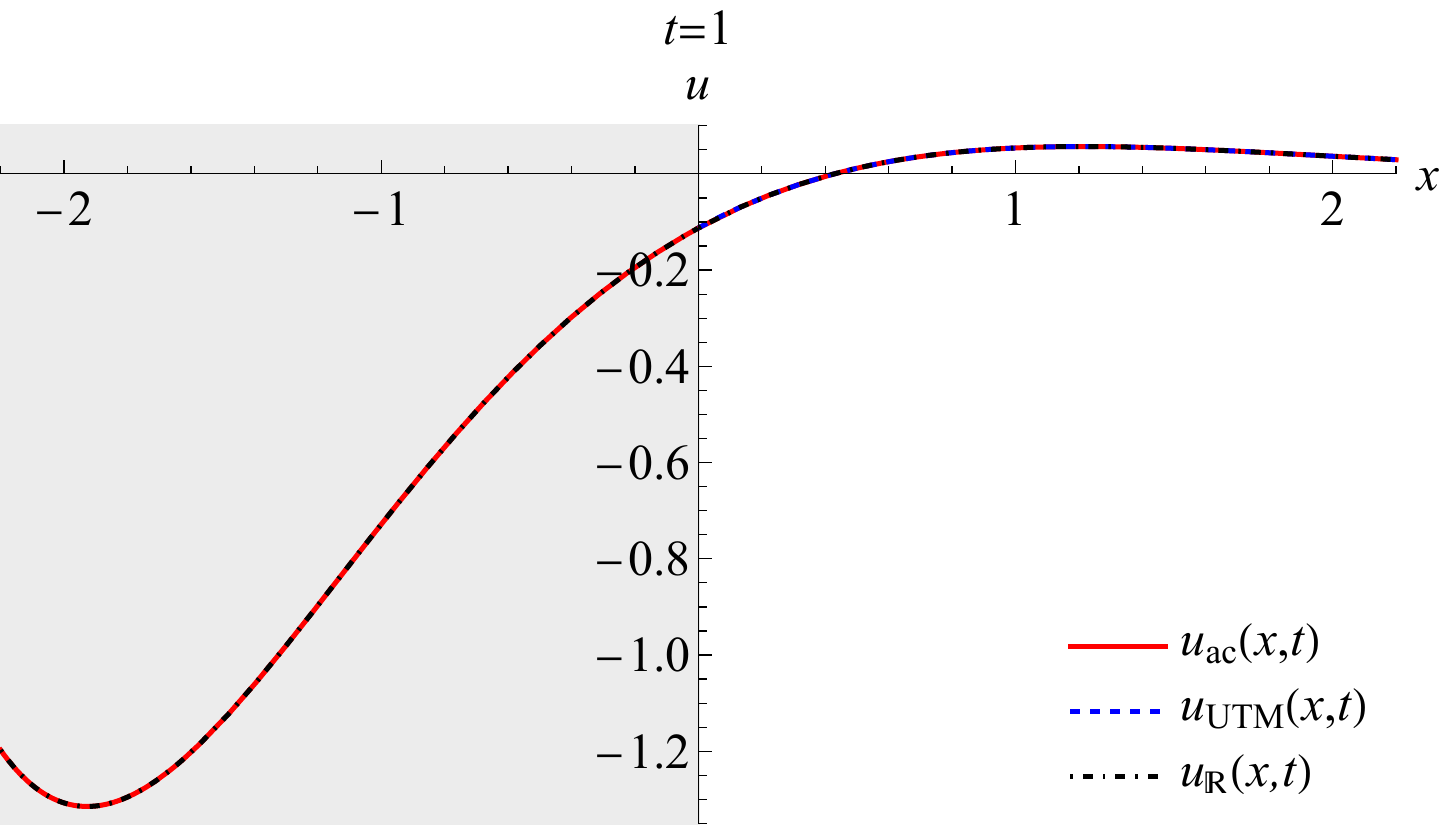}\\
(a) & (b)
\end{tabular}
\end{center}
\caption{(a) The initial condition $w_0(x)$ \rf{eqn:w0kdv1}, leading to the (analytic) solution on the right, shown with the UTM solution at $t=0$ for $x>0$ (it is not defined for $x<0$) and the whole-line initial condition $u_0(x) = u_\R(x,0)$. (b) The solution $u_{\text{ac}}(x,t)$ at $t=1$ obtained through analytic continuation, shown with the whole line solution, $u_\R(x,t)$, and the UTM solution for $x>0$ (it is not defined for $x<0$).}
\la{fig:kdv1}
\end{figure}

Next, we consider \rf{kdv1ibvp} with $f_0(t) = te^{-t}$ and $u_0(x) = u_\R(x,0)$ for $x>0$. Since $f_0^{(n)}(0) = -(-1)^n n,$ 
\beq
\tilde f_0(x,0) = -\frac{1}{3}xe^{x} + \frac{2}{3}xe^{-\frac{x}{2}}\sin\left(\frac{\sqrt{3}}{2} x+ \frac{\pi}{6}\right),
\eeq
\no and we find the analytic continuation of the solution and the corresponding initial condition shown in Figure~\ref{fig:kdv1alt}.

\begin{figure}[tb]
\begin{center}
\def \sc {0.475}
\begin{tabular}{cc}
\includegraphics[scale=\sc]{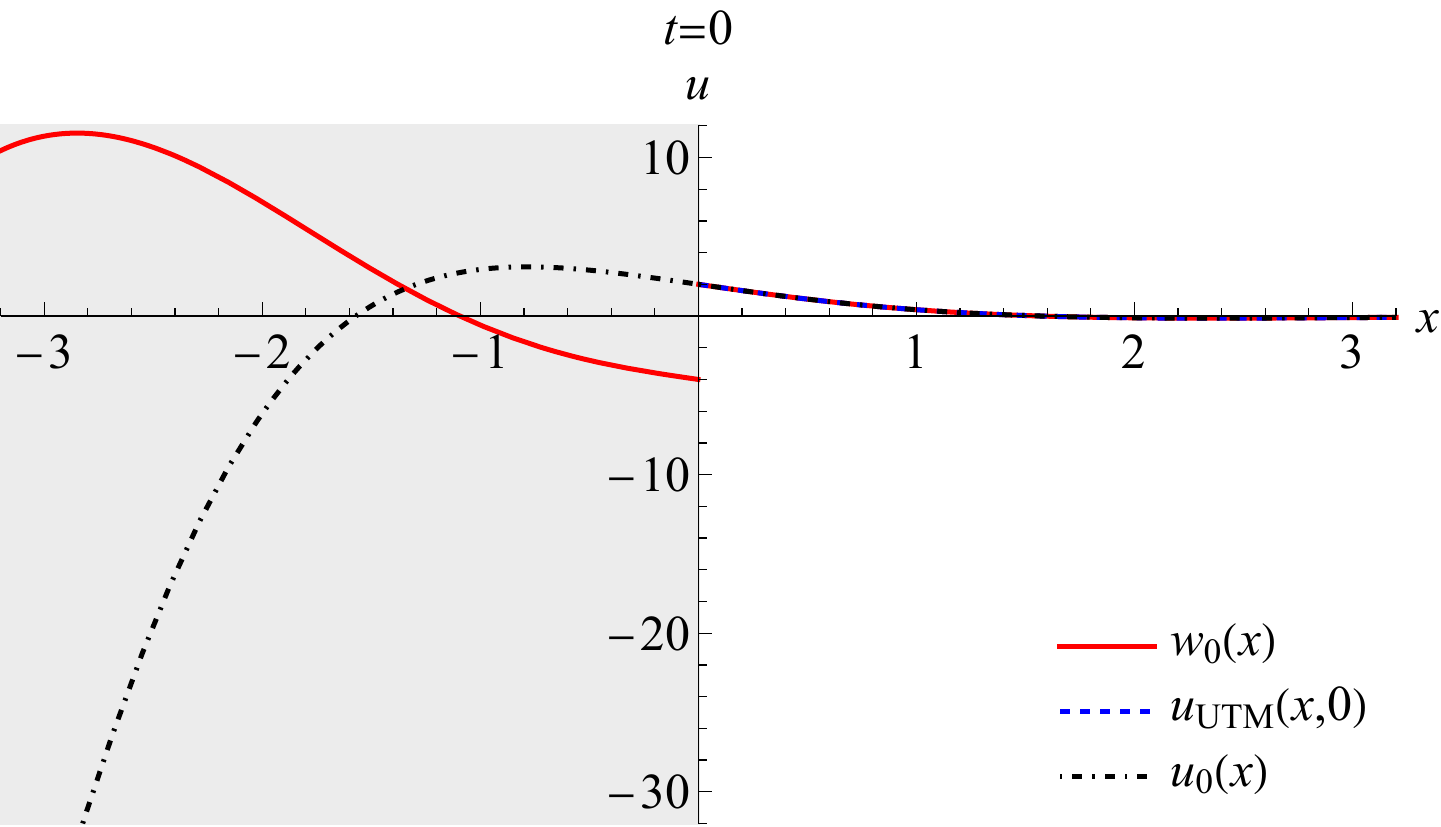} & \includegraphics[scale=\sc]{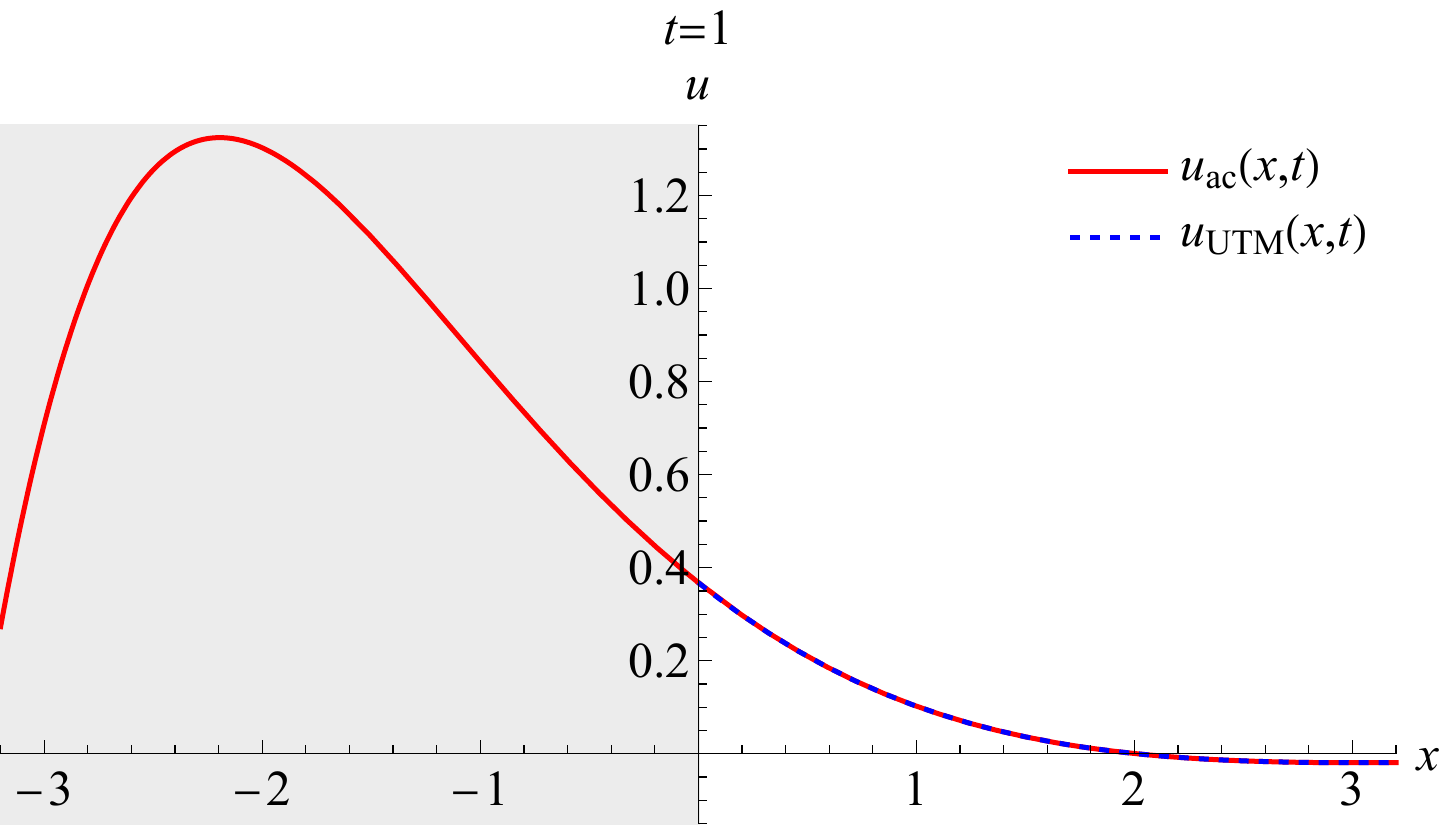}\\
(a) & (b)
\end{tabular}
\end{center}
\caption{(a) The (discontinuous) initial condition $w_0(x)$ \rf{eqn:w0kdv1} leading to the (analytic) solution on the right, shown with the analytic continuation of $u_0(x)=u_\R(x,0)$. (b) The solution $u_{\text{ac}}(x,t)$ at $t=1$ obtained through analytic continuation, and the UTM solution for $x>0$ (it is not defined for $x<0$).}
\la{fig:kdv1alt}
\end{figure}

\section{The linear KdV equation, 2 boundary conditions} \label{sec:kdv2}

Consider the version of linear KdV with 2 boundary conditions

\begin{subequations}\la{kdvprob2}
\begin{align}\la{kdv2ibvp}
u_t-u_{xxx} &= 0, &&\hspace*{-1.0in} x>0, \, t>0,\hspace*{1.0in}\\
u(x,0)&=u_0(x), &&\hspace*{-1.0in} x>0,\\
u(0,t)&=f_0(t), &&\hspace*{-1.0in} t>0, \\
u_x(0,t)&=f_1(t), && \hspace*{-1.0in} t>0.
\end{align}
\end{subequations}
\no The solution is given by \cite{JC_fokas_book}
\beq
u(x,t) = I_0(x,t) + I_{f_0}(x,t) + I_{f_1}(x,t),
\eeq
\no where
\begin{align}\la{kdv2I0}
I_0(x,t) &\!= \frac{1}{2\pi} \int_{-\infty}^\infty\!\!\!\!  e^{ikx-ik^3t} \hat u_0(k) dk\! -\! \frac{1}{2\pi} \int_{\partial\Omega_1}\!\!\!\!  e^{ikx-ik^3t}   \hat u_0(\alpha^2k) dk  \!-\!\frac{1}{2\pi} \int_{\partial\Omega_2}\!\!\!\!  e^{ikx-ik^3t}\hat  u_0(\alpha k) dk, \\
I_{f_0}(x,t) &\!=  \frac{1-\alpha}{2\pi} \int_{\partial\Omega_1}  e^{ikx-ik^3t}  k^2F_0(ik^3,t) \, dk +\frac{1-\alpha^2}{2\pi} \int_{\partial\Omega_2}  e^{ikx-ik^3t}k^2F_0(ik^3,t) \, dk, \\
I_{f_1}(x,t) &\!= \frac{1-\alpha^2}{2\pi i} \int_{\partial\Omega_1}  e^{ikx-ik^3t} kF_1(ik^3,t) \, dk +\frac{1-\alpha}{2\pi i} \int_{\partial\Omega_2}  e^{ikx-ik^3t} k F_1(ik^3,t) \, dk,
\end{align}
\no and
\beq
\hat u_0(k) = \int_0^\infty e^{-iky} u_0(y) \,dy, \qquad
F_m(ik^3,t) = \int_0^t e^{ik^3s} f_m(s) \,ds, ~~ m=0,1,
\eeq 
\no and $\alpha = \exp(2\pi i/3)$.

\figcl{\imsize}{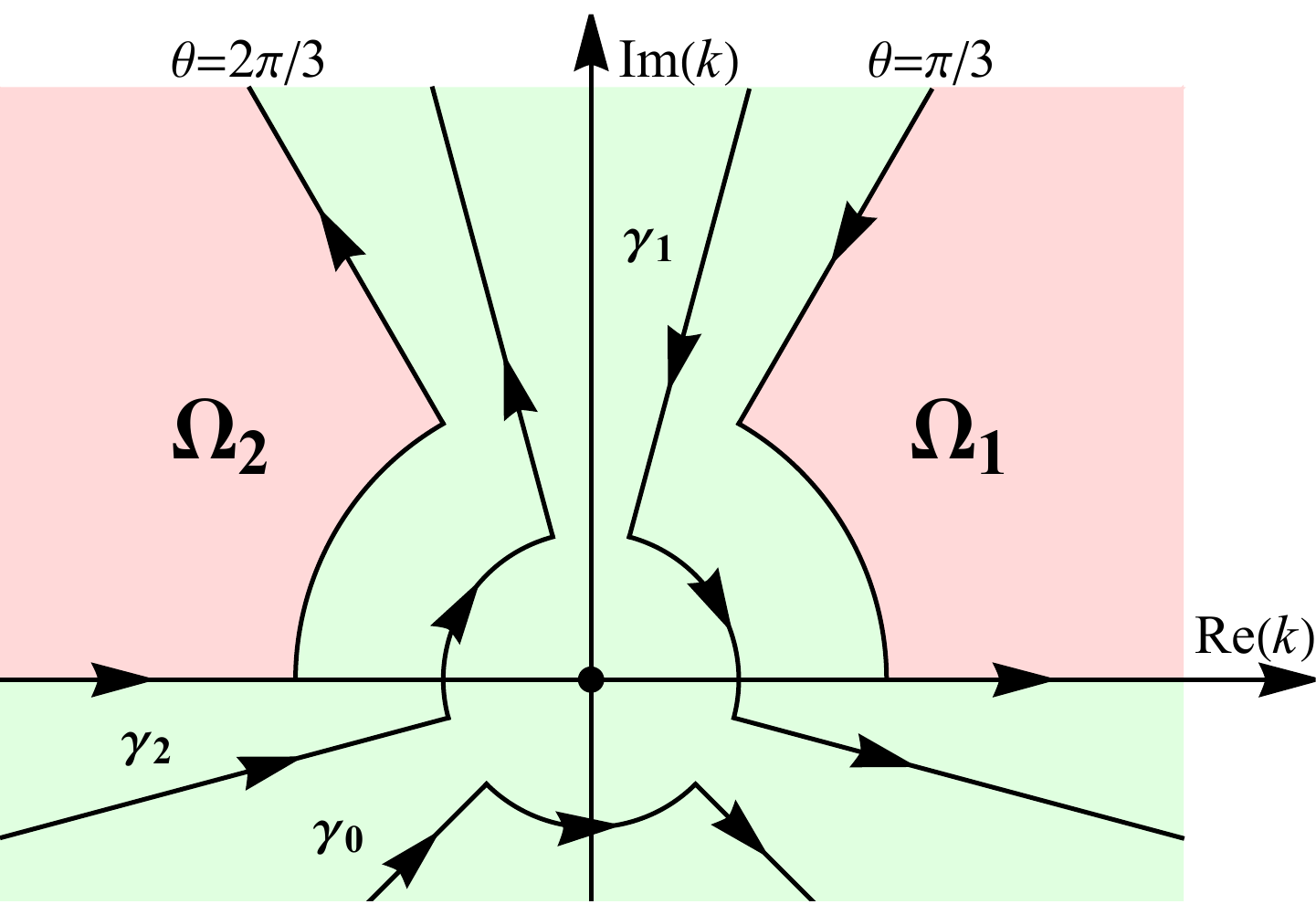}{The regions $\Omega_1$ and $\Omega_2$ for \rf{kdv2ibvp} and the contours $\gamma_0$, $\gamma_1$, and $\gamma_2$.}{kdvOmega}

For the initial condition parts, we deform into the green regions (see Figure~\ref{fig:kdvOmega}) so that
\beq
I_0(x,t) = \frac{1}{2\pi} \int_{\gamma_0}  e^{ikx-ik^3t} \hat u_0(k) \, dk - \frac{1}{2\pi} \int_{\gamma_1}  e^{ikx-ik^3t}   \hat u_0(\alpha^2k) \, dk -\frac{1}{2\pi} \int_{\gamma_2}  e^{ikx-ik^3t}\hat  u_0(\alpha k) \, dk,
\eeq
\no where $\gamma_j$ ($j=1,2$) is the deformation of $\partial\Omega_j$ into the green region (avoiding the origin) defined by $\Re(ik^3)>0$ and $\gamma_0$ is the similar deformation of the real axis into the lower half plane, as shown in Figure~\ref{fig:kdvOmega}. 
%
%

\begin{theorem}
If $u_0\in L^1(\R^+),$ then $I_0(x,t)$ is entire in $x$ for $t>0$.
\end{theorem}

\begin{proof}
If we integrate each of the integrals in $I_0(x,t)$ over a closed contour $\Gamma$ in the complex $x$-plane, we get zero by Cauchy's theorem after switching the order of integration. This is allowed by Fubini's theorem, since
\[
\oint_\Gamma |dx| \int_{\gamma_0} \left|e^{ikx-ik^3t}\hat u_0(k) \right| \, |dk| \leq 2\ell(\Gamma) \|u_0\|_1 \max_{x\in\Gamma} \int_0^\infty e^{2|x|\rho-\rho^3t} \, d\rho <\infty,
\]
\no similarly for the other two terms. Thus, by Morera's theorem, $I_0(x,t)$ is entire in $x$.
\end{proof}

\no It is interesting to compare this to Theorem~\ref{thm:IC_KdV}, where exponential decay is needed for analyticity. We defer the proof for $I_{f_0}(x,t)$ and $I_{f_1}(x,t)$ to Theorem~\ref{thm:convkdv2} below.

Define
\beq
I_{f_0,j}(x,t) = \frac{1-\alpha^j}{2\pi} \int_{\partial\Omega_j}  k^2e^{ikx-ik^3t}F_0(ik^3,t)\, dk,~~j=1,2, 
\eeq
\no so that
\beq
I_{f_0}(x,t) = I_{f_0,1}(x,t)+I_{f_0,2}(x,t).
\eeq
\no We 
integrate by parts $n$ times, resulting in
\begin{align}
I_{f_0,j}(x,t) 
\label{eqn:If0IbP}
&=  \frac{1-\alpha^j}{2\pi} \int_{\partial\Omega_j} dk \, k^2 e^{ikx-ik^3t} \left[ \sum_{m=1}^{n}\frac{f_0^{(m-1)}(0)}{(-ik^3)^{m}} + \frac{1}{(-ik^3)^{n}}\int_0^t e^{ik^3s} f_0^{(n)}(s) \, ds\right].
\end{align}
\no The terms involving $f_0^{(m)}(t)e^{ik^3t}$ integrate to zero around the contour $\partial \Omega_j$ by Jordan's lemma and Cauchy's theorem. We can swap the order of integration by Fubini's Theorem. Technically, we can only do this for $p\geq 2$, but we can also do it for $p=1$, if we do another integration by parts, switch the order, and undo the integration by parts. Then,
\beq
I_{f_0,j}(x,t) = \sum_{m=1}^{n} f_0^{(m-1)}(0)\phi_{m,j}(x,t) + \int_0^t  f_0^{(n)}(s) \phi_{n,j}(x,t-s) \, ds,
\eeq
\no with 
\beq
\phi_{m,j}(x,t) = \frac{1-\alpha^j}{2\pi} \int_{\partial \Omega_j} \frac{k^2 e^{ikx-ik^3t}}{(-ik^3)^{m}} \, dk =\frac{1-\alpha^j}{2\pi} \int_{\gamma_j} \frac{k^2 e^{ikx-ik^3t}}{(-ik^3)^{m}} \, dk,
\eeq
\no where we deform $\partial \Omega_j$ to $ \gamma_j$ in the same way as before, which we can do since $m\geq 1$. It follows that $\phi_{m,j}(x,t)$ is entire in $x$. This expression for $I_{f_0,j}(x,t)$ is smooth up to the $3n$-th $x$-derivative at $x=0$. Then,
\beq
\label{eqn:kdv2IPB}
\frac{\partial^{3n-q} I_{f_0,j}}{\partial x^{3n-q}} \bigg|_{x=0} =\sum_{m=1}^{n} f_0^{(m-1)}(0) \phi_{m,j}^{(3n-q)}(0,t)  +   \int_0^t f_0^{(n)}(s) \phi_{n,j}^{(3n-q)}(0,t-s) \, ds,
\eeq
\no for $q=1,2$. Note that we can swap the limit and integral since
\beq
\phi_{m,j}^{(3n-q)}(0,t) =\frac{1-\alpha^j}{2\pi} \int_{\gamma_j}  \frac{ (ik)^{3n+2-q}e^{-ik^3t}}{(-ik^3)^{m}}\, dk=\frac{1-\alpha^j}{2\pi t^{n-m+1-\frac{q}3}} \int_{\gamma_j}  \frac{ (ik)^{3n+2-q}e^{-ik^3}}{(-ik^3)^{m}}\, dk,
\eeq
\no so that $\phi_{n,j}^{(3n-q)}(0,t) = \mathcal O\left(t^{-1+\frac{q}{3}}\right)$, and the integrand is absolutely integrable. For $j=2$, let $k\mapsto \alpha k$, so that $\gamma_2 \mapsto \gamma_1$, and
\beq
\phi_{m,2}^{(3n-q)}(0,t) =\frac{(1-\alpha^2)\alpha^{-q}}{2\pi} \int_{\gamma_1}  \frac{ (ik)^{3n+2-q} e^{-ik^3t}}{(-ik^3)^{m}}\, dk.
\eeq
\no Therefore, after parametrizing $\gamma_1$ and computing the integral, 
\beq
\label{eqn:kdv2phi}
\sum_{j=1}^2\phi_{m,j}^{(3n-q)}(0,t) = \delta_{q,1} \frac{-\sqrt{3}(-1)^{n-m}\Gamma\left(n-m + \frac{2}{3}\right)}{2\pi t^{n-m+\frac23}},
\eeq
\no where $\delta_{q,1}=1$ if $q=1$ and $0$ otherwise.
%
%
%
%
%
%
%
With
\beq
\label{eqn:kdv2series}
I_{f_0}(x,t) = \sum_{n=0}^\infty a_n(t) x^n, \qquad \text{ and } \qquad I_{f_1}(x,t) = \sum_{n=0}^\infty b_n(t) x^n,
\eeq
\no the above leads to
\begin{align}\nonumber
a_{3n-1}(t) &=-\frac{\sqrt{3}\Gamma\left(\frac23\right)}{2\pi(3n-1)!}\left[\sum_{m=1}^{n}  \frac{(-1)^{n-m}\Gamma\left(n-m+\frac23\right)}{t^{n-m+\frac23}\Gamma\left(\frac23\right)}f_0^{(m-1)}(0)+ \int_0^t ds \, \frac{f_0^{(n)}(s)}{(t-s)^\frac23}\right]\\
\label{eqn:a1kdv2}
&=-\frac{f_0^{\left(n-\frac13\right)}(t)}{(3n-1)!}, ~~n=1,2,\ldots
\end{align}
\no and 
\beq
\label{eqn:a2kdv2}
a_{3n-2}(t)=0, ~~n=1,2,\ldots 
\eeq

If $q=0$, 
%
%
%
%
%
%
%
%
we cannot pass the limit under the integral sign to get \rf{eqn:kdv2IPB}. Using \rf{eqn:kdv2phi},
\begin{align}\nonumber
\frac{\partial^{3n}I_{f_0}}{\partial x^{3n}}\bigg|_{x=0}\!\! &= \int_0^t \!\!ds \, f_0^{(n)}(s) \sum_{j=1}^2\phi_{n,j}^{(3n)}(x,t-s) \\
&= \frac{1}{2\pi} \int_0^t \!\!ds \, f_0^{(n)}(s)\!\left[\!(1-\alpha)\!\! \int_{\gamma_1} \!k^2 e^{ikx-ik^3(t-s)}dk \!+\!(1-\alpha^2) \int_{\gamma_2} \!k^2 e^{ikx-ik^3(t-s)}dk \!\right]\!.
\end{align}
\no Let $z=k^3$ so that $dz = 3k^2\, dk$, and
\begin{align}\nonumber
a_{3n}(t) &= \frac{1}{(3n)!}\frac{\partial^{3n}I_{f_0}}{\partial x^{3n}}\bigg|_{x=0} \\ \nonumber
&= \frac{1}{6\pi(3n)!} \int_0^t ds \, f_0^{(n)}(s)\left[(1-\alpha) \int_{-\infty}^\infty e^{-iz(t-s)} dz +(1-\alpha^2) \int_{-\infty}^\infty e^{-iz(t-s)} dz\right] \\\label{eqn:a0kdv2}
&= \frac{1}{(3n)!}  \int_0^t f_0^{(n)}(s)\delta(t-s) ds=\frac{f_0^{(n)}(t)}{(3n)!},~~n=0, 1, \ldots 
\end{align}
\no which we can confirm rigorously using asymptotics. Similarly, we can find the coefficients for the Taylor series for $I_{f_1}(x,t)$ ($n=1,2$,),
\begin{align}\nonumber 
b_{3n-1}(t) &= -\frac{\sqrt{3}\Gamma\left(\frac13\right)}{2\pi(3n-1)!}\left[\sum_{m=1}^{n}  \frac{(-1)^{n-m}\Gamma\left(n-m+\frac13\right)}{t^{n-m+\frac13}\Gamma\left(\frac13\right)}f_1^{(m-1)}(0)+ \int_0^t ds \, \frac{f_1^{(n)}(s)}{(t-s)^\frac13}\right]\\
\label{eqn:b1kdv2}
&=-\frac{f_1^{\left(n-\frac23\right)}(t)}{(3n-1)!}, ~~n=1,2, \ldots\\
\label{eqn:b0kdv2}
b_{3n}(t)&= 0, ~~~~~~~~~~~~~~~\,~n=0, 1, \ldots\\
\label{eqn:b2kdv2}
b_{3n+1}(t) &= \frac{f_1^{(n)}(t)}{(3n+1)!}, ~~~~~~n=0, 1, \ldots.
\end{align}
\no Then, the two series in \rf{eqn:kdv2series} give analytic extensions for the two functions $I_{f_0}(x,t)$ and $I_{f_1}(x,t)$. To reduce the number of terms in the Taylor series, we may once again write
\begin{align} \label{eqn:kdv2ext0}
I_{f_0}^{\text{ext}} (x,t) &= \case{1}{I_{f_0}(x,t), & x\geq 0, \\ \tilde f_0(x,t) - I_{f_0}(-x,t), & x < 0,} \quad \text{where} \quad \tilde f_0(x,t) = 2\sum_{n=0}^\infty a_{2n}(t) x^{2n}, \\ \label{eqn:kdv2ext1}
I_{f_1}^{\text{ext}} (x,t) &= \case{1}{I_{f_1}(x,t), & x\geq 0, \\ \tilde f_1(x,t) + I_{f_1}(-x,t), & x < 0,}\quad \text{where} \quad \tilde f_1(x,t) = 2\sum_{n=1}^\infty b_{2n-1}(t) x^{2n-1}.
\end{align}

\begin{theorem}
\label{thm:convkdv2} Define $\mathcal D=\{t\in \C: \text{dist}(t,[0,T])\leq r\}$, for some $r>0$, $T>0$, a domain in the complex $t$-plane containing the interval $[0,T]$.  If $f_0(t)$ $(f_1(t))$ is analytic in $\mathcal D$ then $I_{f_0}(x,t)$ $(I_{f_1}(x,t))$ is analytic for $x>0$ and $0\leq t \leq T$, and the series representation \rf{eqn:kdv2series} is entire in $x$ for $0 \leq t \leq T$.
\end{theorem}

\begin{proof}
The proof for the convergence of the series is similar to that of Theorem~\ref{thm2}. Since the series is absolutely convergent, and it solves the corresponding initial and boundary value problem, by uniqueness, it converges to $I_{f_0}(x,t)$ $(I_{f_1}(x,t))$, so that $I_{f_0}(x,t)$ $(I_{f_1}(x,t))$ is an analytic function for $x>0$ and for $0\leq t \leq T$.
\end{proof}

The previously used combination of the Cauchy, Fubini, and Morera theorems is insufficient here. We need the stronger condition of analyticity of $f_0(t)$ ($f_1(t)$) to determine the analyticity of $I_{f_0}(x,t)$ ($I_{f_1}(x,t)$). This is to be compared to the case of the previous section (linear KdV, one boundary condition), where the stronger condition of exponential decay on $u_0(x)$ was needed to determine the analyticity of $I_0(x,t)$.

\subsection{Boundary-to-Initial Map}

The functions
\beq
\tilde f_0(x,t) = 2\sum_{n=0}^\infty a_{2n}(t) x^{2n}, \qquad \text{ and } \qquad \tilde f_1(x,t) = 2\sum_{n=0}^\infty b_{2n}(t) x^{2n},
\eeq
\no are solutions to the linear KdV equation \rf{kdv2ibvp}, which is easily seen from the expressions for the Taylor coefficients in the previous section. Since
\beq
\sum_{n=0}^\infty a_{3n}(t)x^{3n} = \sum_{n=0}^\infty \frac{x^{3n}}{(3n)!} f_0^{(n)}(t) \to \sum_{n=0}^\infty \frac{x^{3n}}{(3n)!} f_0^{(n)}(0), ~~~\mbox{as}~t\to0^+,
\eeq
\no and
\beq
\sum_{n=1}^\infty a_{3n-1}(t)x^{3n-1} 
\sim -\frac{\sqrt{3}}{2\pi}\sum_{n=1}^\infty\frac{x^{3n-1}}{(3n-1)!} \sum_{m=1}^{n}  \frac{(-1)^{n-m}\Gamma\left(n-m+\frac23\right)}{t^{n-m+\frac23}}f_0^{(m-1)}(0),~~~\mbox{as}~t\to0^+,
\eeq
\no and the integral term approaches zero as $t\to 0^+$. 
%
%
%
%
%
%
%
%
Switching the order of summation and evaluating the inner sum, we find for $x<0$,
\beq
\sum_{n=1}^\infty a_{3n-1}(t)x^{3n-1} \sim -\frac{\sqrt{3}\Gamma\left(\frac23\right)}{2\pi t^{\frac23}} \sum_{m=0}^{\infty}\frac{x^{3m+2}}{(3m+2)!}f_0^{(m)}(0){}_2F_3\left(\begin{array}{c}\frac23,\, 1\\ m+1,\, m+\frac43,\, m+\frac53\end{array};\,\frac{|x|^3}{27t}\right),
\eeq
\no as $t\to0^+$, and since
\beq
\frac{\sqrt{3}\Gamma\left(\frac23\right)}{2\pi t^{\frac23}}{}_2F_3\left(\begin{array}{c}\frac23,\, 1\\ m+1,\, m+\frac43,\, m+\frac53\end{array};\,\frac{|x|^3}{27t}\right)
\sim  \frac{3^{\frac{3m}{2}+\frac{3}{4}}(3m+2)!t^{\frac{3m}{2}+\frac{1}{4}}}{\sqrt{\pi}|x|^{\frac{9m}{2}+\frac{11}{4}}}e^{\frac{2|x|^\frac32}{3\sqrt{3t}}}, ~~~\mbox{as}~t\to0^+,
\eeq
\no we find for $x<0$,
\beq
\label{eqn:If0inf}
\sum_{n=1}^\infty a_{3n-1}(t)x^{3n-1} \sim 
-e^{\frac{2|x|^\frac32}{3\sqrt{3t}}}\sum_{m=0}^{\infty}  \frac{(-1)^{m}3^{\frac{3m}{2}+\frac{3}{4}} t^{\frac{3m}{2}+\frac{1}{4}}}{\sqrt{\pi} |x|^{\frac{3m}{2}+\frac{3}{4}}} f_0^{(m)}(0),~~~\mbox{as}~t\to0^+,
\eeq
\no and, in general, $|I_{f_0}(x,t)| \to \infty$, as $t\to 0^+$. Similarly, for $x<0$, we have
\beq
\label{eqn:If1inf}
\sum_{n=1}^\infty b_{3n-1}(t) x^{3n-1} \sim 
-e^{\frac{2|x|^\frac32}{3\sqrt{3t}}} \sum_{m=0}^{\infty} \frac{3^{\frac{5}{4}+\frac{3m}{2}} t^{\frac{3}{4}+\frac{3m}{2}}}{2 \sqrt{\pi} |x|^{\frac{5}{4}+\frac{3m}{2}}}  f_1^{(m)}(0), ~~~\mbox{as}~t\to 0^+,
\eeq
\no and $|I_{f_1}(x,t)| \to \infty$ as $t\to 0^+$ as well. 

We write 
\beq
I_0(x,t) = I_0^{(0)}(x,t) + I_0^{(1)}(x,t) + I_0^{(2)}(x,t) =  I_0^{(0)}(x,t) + 2\Re\left\{I_0^{(1)}(x,t)\right\},
\eeq
\no where we define
\beq
I_0^{(0)}(x,t) = \frac{1}{2\pi} \int_{-\infty}^\infty e^{ikx-ik^3t} \hat u_0(k) \, dk \to \case{1}{u_0(x), & x>0, \\ 0, & x<0,}
~~~\mbox{as}~t\to0^+,
\eeq
\no However, the other terms in $I_0(x,t)$ diverge to infinity as $t\to 0^+$, as shown next. We have
\beq
I_0^{(1)}(x,t) = -\frac{1}{2\pi} \int_0^\infty dy\, u_0(y) \int_{\gamma_1} e^{ik(x-\alpha^2y)-ik^3t} \, dk = \frac{\alpha}{\sqrt[3]{3t}} \int_0^\infty \Ai\left(\frac{y-\alpha x}{\sqrt[3]{3t}}\right) u_0(y) \, dy.
\eeq
\no Using $z=(y-\alpha x)/\sqrt[3]{3t}$, we find

\figcl{\imsize}{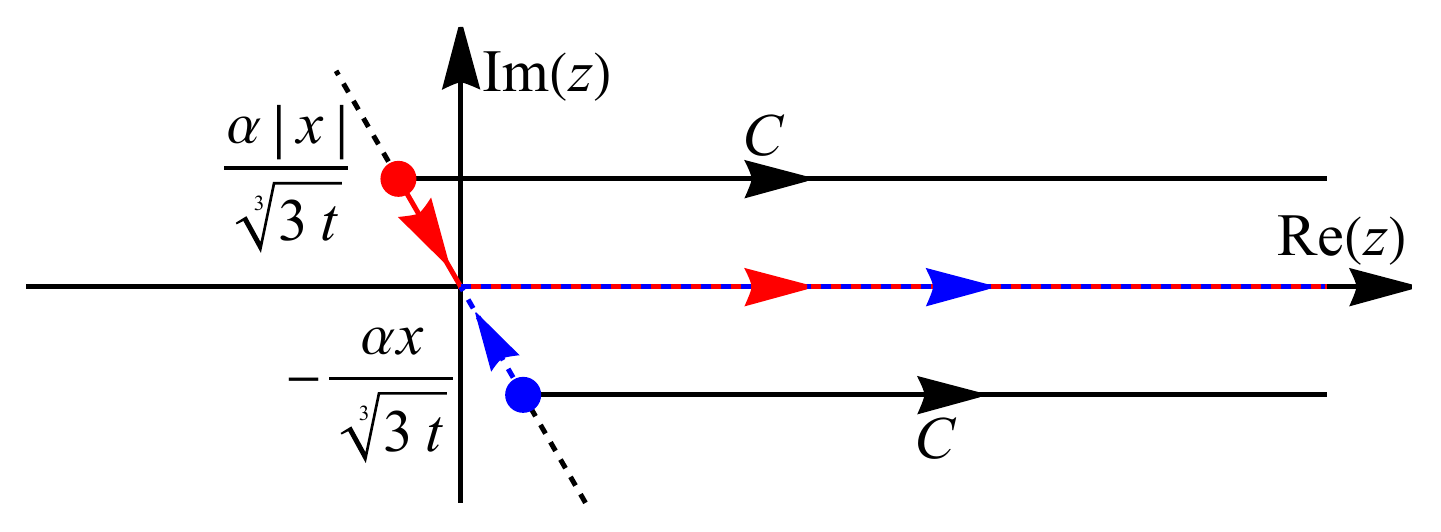}{Contour deformation for the asymptotics of $I_0^{(1)}(x,t)$ for $x>0$ (blue, dashed) and for $x<0$ (red, solid).}{kdvcontourdeform}

\begin{align}\nonumber
I_0^{(1)}(x,t) &= \alpha \int_{C} \Ai(z) u_0\left(\alpha x+z\sqrt[3]{3t}\right) \, dz\\ 
&= \alpha \int_{0}^\infty \Ai(z) u_0\left(\alpha x+z\sqrt[3]{3t}\right) \, dz - \alpha^2 \int_{0}^{-{x}/{\sqrt[3]{3t}}} \Ai(\alpha\rho) u_0\left(\alpha x+\alpha \rho\sqrt[3]{3t}\right) \, d\rho,
\end{align}
\no where the positive real axis is mapped to the contour $C$ as shown in Figure~\ref{fig:kdvcontourdeform}. As $t\to0^+$, the first term limits to $\alpha u_0(\alpha x)/3$. For $x>0$, the second term has the limit $-\alpha u_0(\alpha x)/3$, so that $I_0(x,t)\to 0$ as $t\to 0^+$, as expected. For $x<0$, 
\begin{align}\nonumber
I_0^{(1)}(x,t) 
&\sim \frac{\alpha}{3}u_0(\alpha x)  - \alpha^2\sum_{n=0}^\infty \frac{\alpha^n}{n!} u_0^{(n)}(0)\int_{0}^{{|x|}/{\sqrt[3]{3t}}} \Ai(\alpha\rho)  \left(x+\rho\sqrt[3]{3t}\right)^n  \, d\rho \\\nonumber
&\sim - \alpha^2\sum_{n=0}^\infty \frac{\alpha^n}{n!} u_0^{(n)}(0)\sum_{j=0}^n \binom{n}{j} x^{n-j}(3t)^{\frac{j}{3}}\int_{0}^{{|x|}/{\sqrt[3]{3t}}} \Ai(\alpha\rho) \rho^j  \, d\rho \\\nonumber
&=\sum_{n=0}^\infty \frac{\alpha^{n+2}x^{n+1}}{n!} u_0^{(n)}(0)\sum_{j=0}^n \binom{n}{j}  \frac{(-1)^{j}}{3(j+1)\Gamma\left(\frac23\right)t^{\frac{1}{3}}} {}_1F_2\left(\begin{array}{c}\frac{j+1}{3}\\ \frac23, \, \frac{j+4}{3} \end{array}; \, \frac{|x|^3}{27t}\right) \\
&~~~+\sum_{n=0}^\infty \frac{\alpha^{n}x^{n+2}}{n!} u_0^{(n)}(0)\sum_{j=0}^n \binom{n}{j} \frac{(-1)^{j}}{3(j+2)\Gamma\left(\frac13\right)t^{\frac{2}{3}}}{}_1F_2\left(\begin{array}{c} \frac{j+2}{3}\\ \frac{4}{3}, \, \frac{j+5}{3}\end{array}; \frac{|x|^3}{27t}\right),
\end{align}
\no where the bounded function $\alpha u_0(\alpha x)/3$ was omitted from the asymptotic series because it is not part of the leading-order behavior. The last integral is evaluated by integrating the Maclaurin series for the Airy function and manipulating the resulting series, obtaining the result in terms of hypergeometric functions. Using the asymptotic expansions for the generalized hypergeometric functions ${}_{p}F_q$ from \cite[Section~16.11]{dlmf}, we find
\beq
\label{eqn:I0inf}
I_0(x,t) \sim e^{\frac{2|x|^{\frac32}}{3\sqrt{3t}}} \sum_{n=0}^\infty  \frac{(-1)^{n} 3^{\frac{3n}{2}+ \frac{3}{4}} t^{\frac{3n}{2} +\frac{1}{4}}}{\sqrt{\pi}|x|^{\frac{3n}{2}+\frac34}}u_0^{(3n)}(0)+e^{\frac{2|x|^{\frac32}}{3\sqrt{3t}}}\sum_{n=0}^\infty  \frac{(-1)^{n} 3^{\frac{3n}{2}+ \frac{5}{4}} t^{\frac{3n}{2} +\frac{3}{4}}}{\sqrt{\pi}|x|^{\frac{3n}{2}+\frac54}}u_0^{(3n+1)}(0),
\eeq
\no so that, in general, $|I_0(x,t)| \to \infty$, as $t\to 0^+$.
It is clear that for some initial and boundary data $u_0(x),$ $f_0(t)$, and $f_1(t)$, (for instance, $f_0(t)=f_1(t)=0$, $u_0(x)\neq 0$), a boundary-to-initial map does not exist. However, if the compatibility conditions, 
\beq
f_0^{(n)}(0)=u_0^{(3n)}(0), \qquad \text{ and } \qquad f_1^{(n)}(0)=u_0^{(3n+1)}(0), ~~~~n=0,1,\ldots
\eeq
\no are satisfied, then it can be shown that singular behavior as $t\rightarrow 0^+$ cancels and 
\beq
I_0(x,t)+I_{f_0}^\text{ext}(x,t)+I_{f_1}^\text{ext}(x,t)\to u_0(x),~~~\mbox{as}~~t\to 0^+,
\eeq 
\no where $u_0(x)$ is the analytical continuation of $u_0(x)$ (which is defined for $x>0$) to the whole real line. Indeed, comparing \rf{eqn:I0inf} to \rf{eqn:If0inf} and \rf{eqn:If1inf} is suggestive of this, although a proof requires more work. It is an interesting question to isolate the conditions on the initial and boundary data for which a boundary-to-initial map exists, and the half-line problem can be viewed as the restriction of a whole-line problem, even if the whole-line initial condition is unbounded and discontinuous. 

\subsection{Examples}

We demonstrate our results using two examples. Our first example uses the whole-line solution,
\beq\la{kdv2ex1}
u_{\R}(x,t) = 2e^{-\sqrt{3}x}\cos(x+8t),
\eeq
\no with the inferred boundary and initial conditions. The UTM solution is no longer defined for $x<0$, but the analytic continuation recovers the exact solution on the whole line. The results are shown in Figure~\ref{fig:kdv2}{\textcolor{blue}a}.

\begin{figure}[tb]
\begin{center}
\def \sc {0.475}
\begin{tabular}{cc}
\includegraphics[scale=\sc]{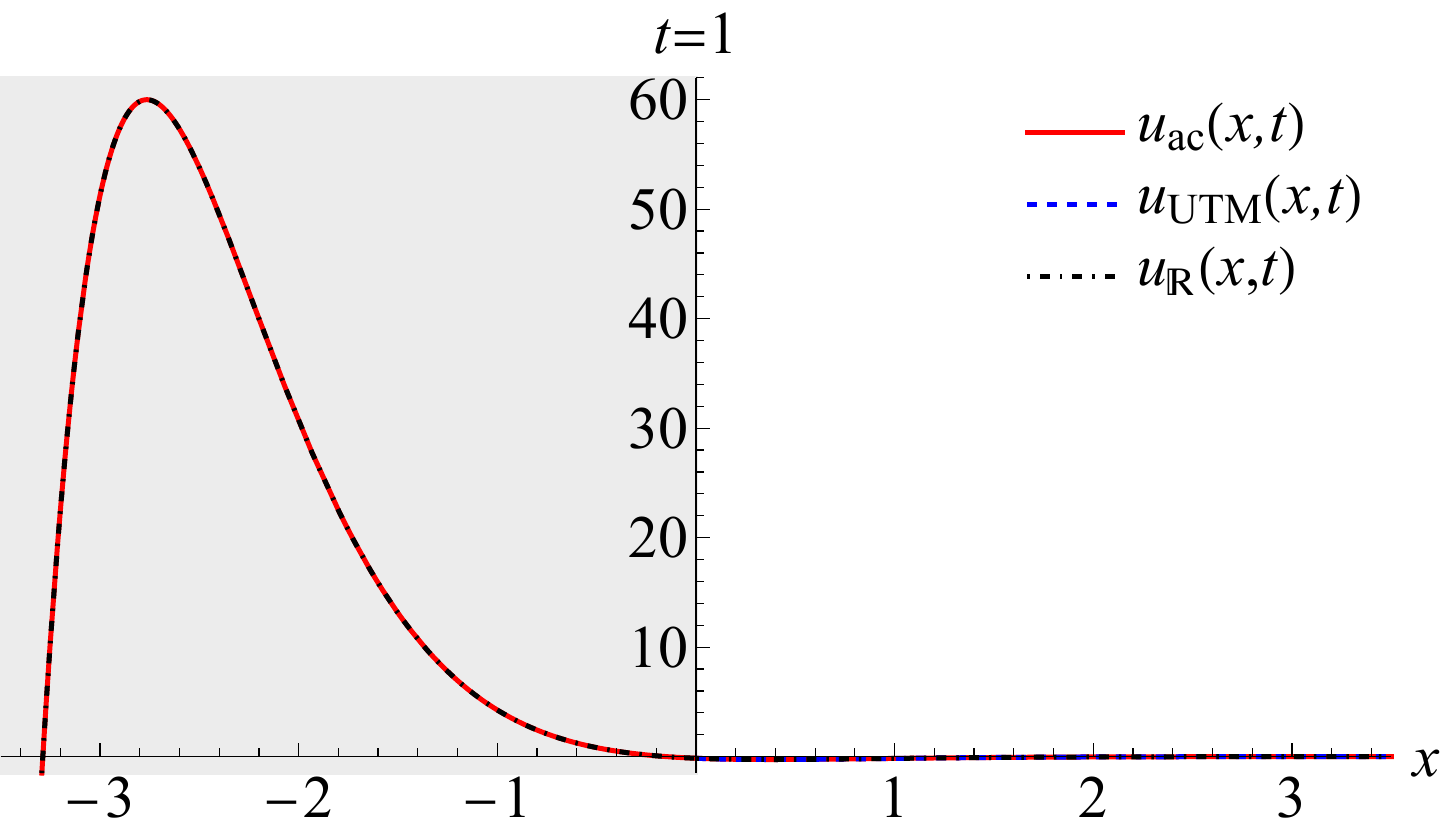} & \includegraphics[scale=\sc]{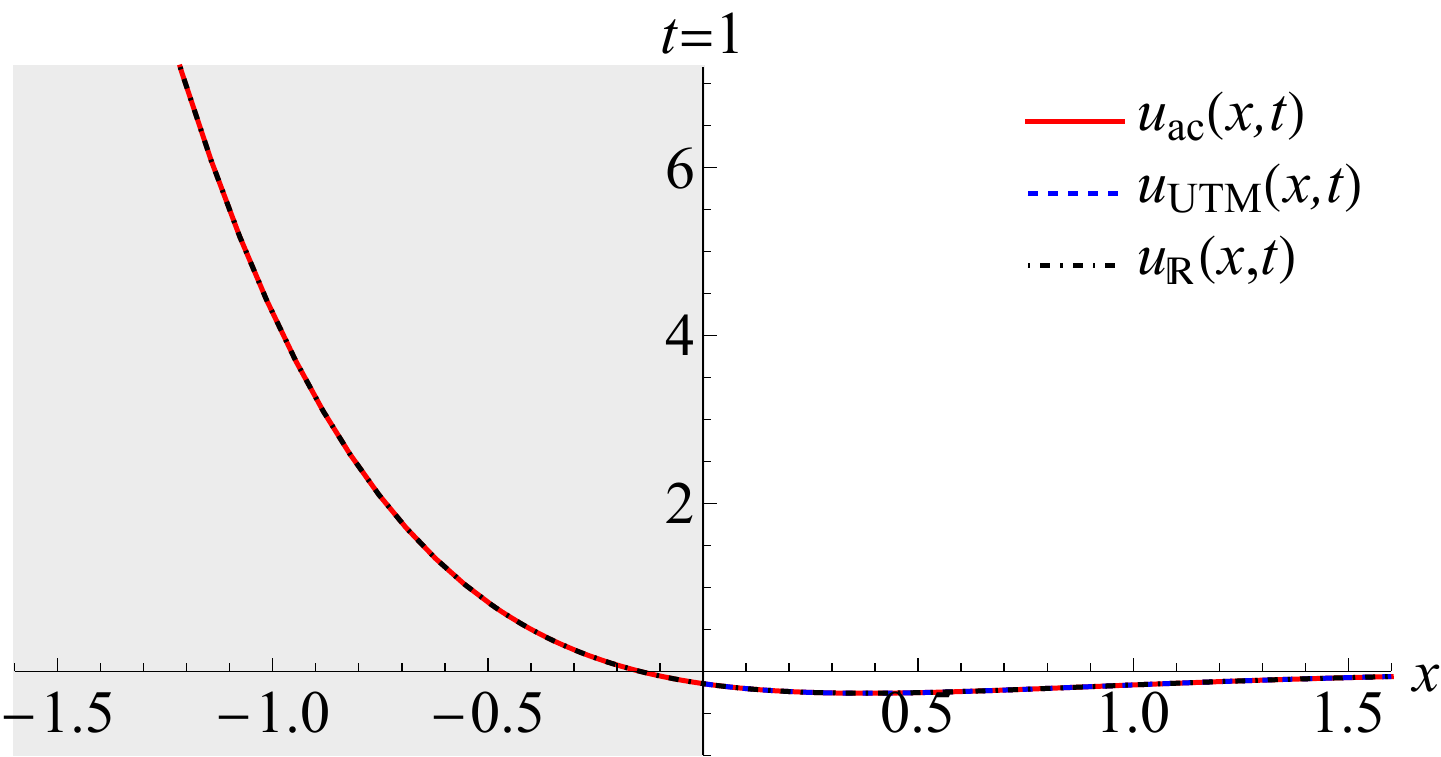}\\
(a) & (b) \\
\includegraphics[scale=\sc]{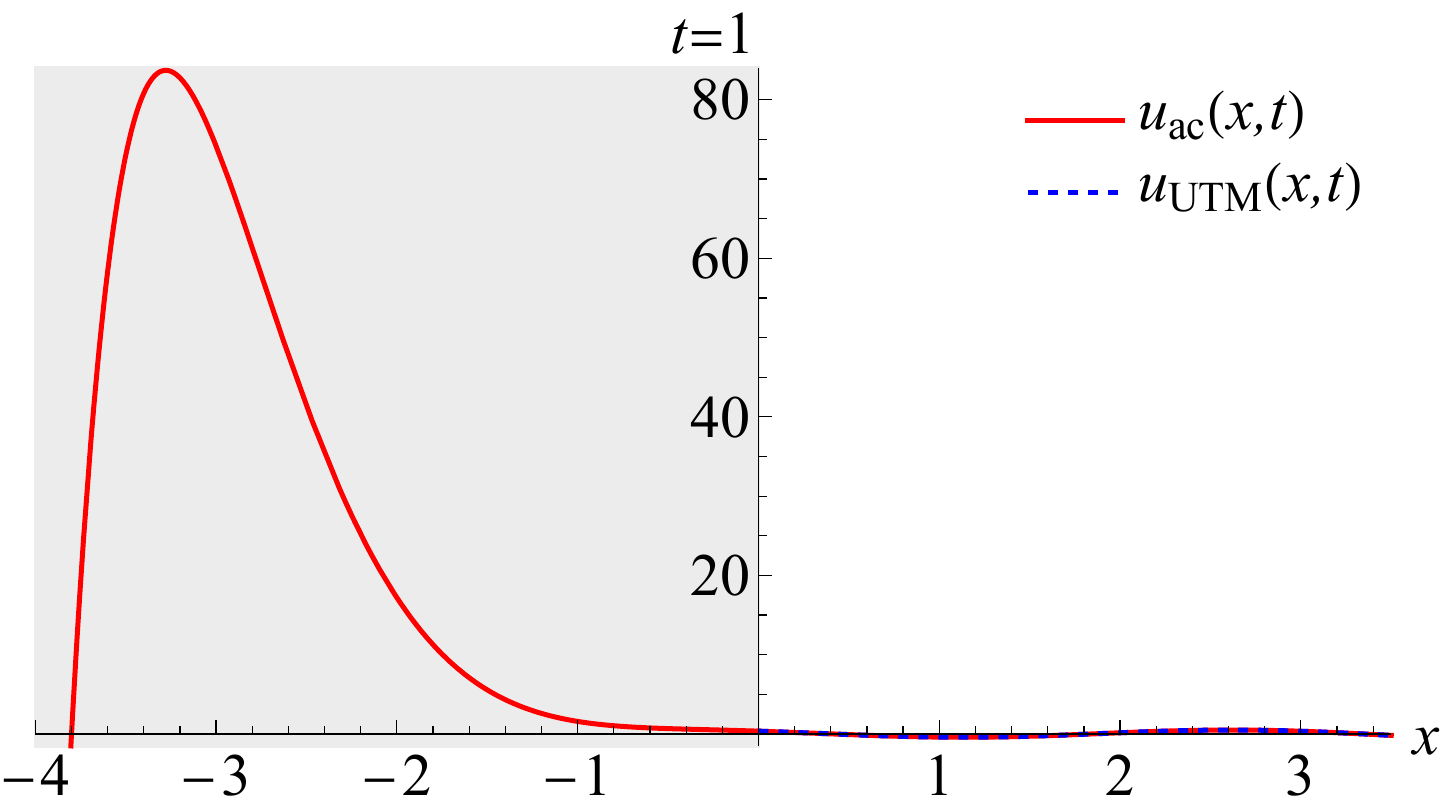} & \includegraphics[scale=\sc]{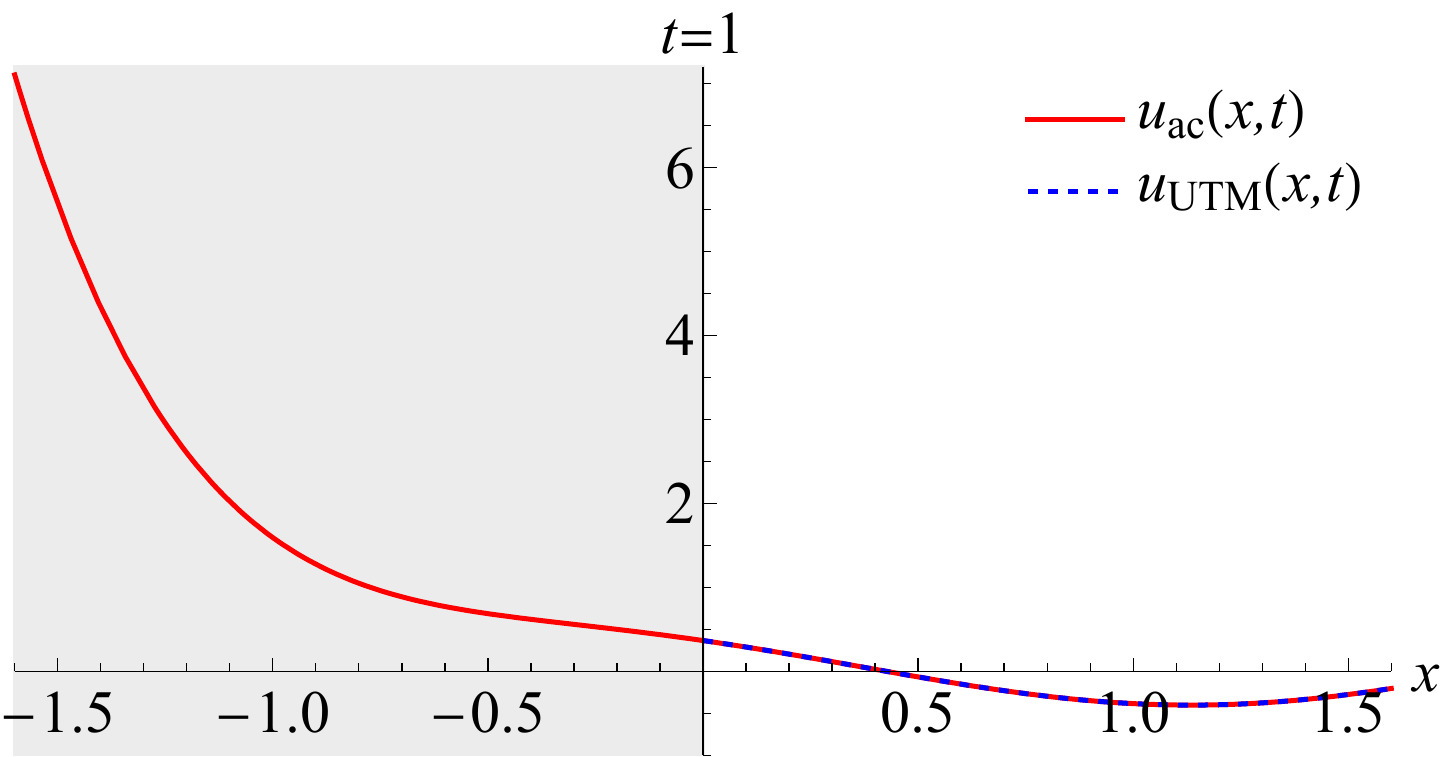}\\
(c) & (d)
\end{tabular}
\end{center}
\caption{(a) The solution $u_{\text{ac}}(x,t)$ at $t=1$ obtained through analytic continuation, shown with the UTM solution $u_\text{UTM}(x,t)$ for $x>0$ (it is not defined for $x<0$) and the whole line solution $u_\R(x,t)$ given by \rf{kdv2ex1}. (b) The same plot zoomed in. (c) The solution $u_{\text{ac}}(x,t)$ at $t=1$, for the incompatible case, obtained through analytic continuation, shown with the UTM solution $u_\text{UTM}(x,t)$ for $x>0$ (it is not defined for $x<0$). (d) The same plot zoomed in. }
\la{fig:kdv2}
\end{figure}

Next, we consider \rf{kdv2ibvp} with $f_0(t) = te^{-t}$, $u_0(x) = u_\R(x,0)$, and $f_1(t)=u_{\R,x}(x,t)|_{x=0}$ for $x>0$. We find the analytic continuation of the solution shown in Figure~\ref{fig:kdv2}{\textcolor{blue}b}.


\section{The discretized advection equation}\label{sec:sdadv}

As discussed above, the analytic continuation of the solution for the continuous-in-space advection equation is immediate through the use of d'Alembert's formula. The situation is more complicated in the discrete space setting, as the discretization used affects how the analytic continuation is done. 

\subsection{Backward discretization}

\label{advec2_backward_halfline_analytcont}
	
Consider, with wave speed $c>0$,
\begin{subequations} \label{advec2_prob}
\begin{align}
		u_t = -c\,u_{x},&~~~~ x > 0,\, t > 0, \\
		u(x,0) = \phi(x),&~~~~ x > 0,\\
		u(0,t) = f_0(t),&~~~~ t > 0.
\end{align}
\end{subequations}
\no Discretizing the spacial derivative $u_x(x,t)$ using the standard backward stencil gives 
\beq \la{advecdiscrete1}
\dot{u}_n(t) = c\,\frac{u_{n-1}(t) - u_{n}(t)}{h},
\eeq 
\no with dispersion relation 
\begin{equation}
	    W(k) = c\,\frac{1 - e^{-ikh}}{h}.
	    \label{advec2_back_W}
\end{equation}

Following \cite{JC_SDUTM_HL}, the solution to this semi-discrete IBVP \eqref{advec2_prob} is
\begin{align}\begin{split}
		u_n(T) &= \frac{1}{2\pi} \int_{-\pi/h}^{\pi/h} e^{iknh} e^{-WT} \hat{u}(k,0)\,dk \,\,+ \,\,\frac{c}{2\pi} \int_{-\pi/h}^{\pi/h} e^{ik(n-1)h} e^{-WT} F_{0}(W,T)\,dk,~~n\geq 1,
		\label{soln_advec2_backward}
\end{split}\end{align}
\no where the Fourier transform
\beq \label{eqn:SDu0hat}
\hat{u}(k,0) =  h \sum_{n=1}^\infty e^{-iknh} u_n(0), \quad\quad \text{Im}(k) \leq 0,
\eeq 
\no begins at $n = 1$ since the Dirichlet boundary condition is given, and
\beq\label{f0}
F_0(W,T) = \int_{0}^T e^{Wt} f_0(t) \, dt, \quad\quad k \in \mathbb{C}.
\eeq

Unlike the solution to the continuous advection equation, the semi-discrete solution \eqref{soln_advec2_backward} couples the initial and boundary conditions, and both contribute at every mesh point. According to Figure~\ref{fig:advec_backward_W}, substituting $n \rightarrow - n$ for $n \in \mathbb{Z}^+$ allows both integral paths to be deformed below the real line where the integrands decay, so that $u_{-n}(T) \equiv 0$ for $n>0$ and for all $T$. In the continuum limit, \textit{i.e.}, as $h \rightarrow 0$, \eqref{soln_advec2_backward} converges to
\begin{align}
	\begin{split}
		u(x,T) &= \frac{1}{2 \pi} \int_{-\infty}^{\infty} e^{ikx} e^{-\tilde{W}T} \hat{u}(k,0) \, dk + \frac{c}{2 \pi} \int_{-\infty}^{\infty} e^{ikx} e^{-\tilde{W}T} \tilde F_0(\tilde{W},T) \, dk,~~x > 0,
		\label{soln_advec2_cont}
	\end{split}
\end{align}
\no with $\tilde{W}(k) = i c k$, $\hat u(k,0)=\int_0^\infty e^{-ikx} u(x,0) dx$, and $\tilde F_0(\tilde W, T)=\int_0^T e^{\tilde W t} u(0,t)dt$. It follows that $u(-x,T) = 0$ for $x > 0$ from this representation. In what follows, we examine how \rf{soln_advec2_backward} may be analytically continued for $n<0$. 

\figcl{\imsize}{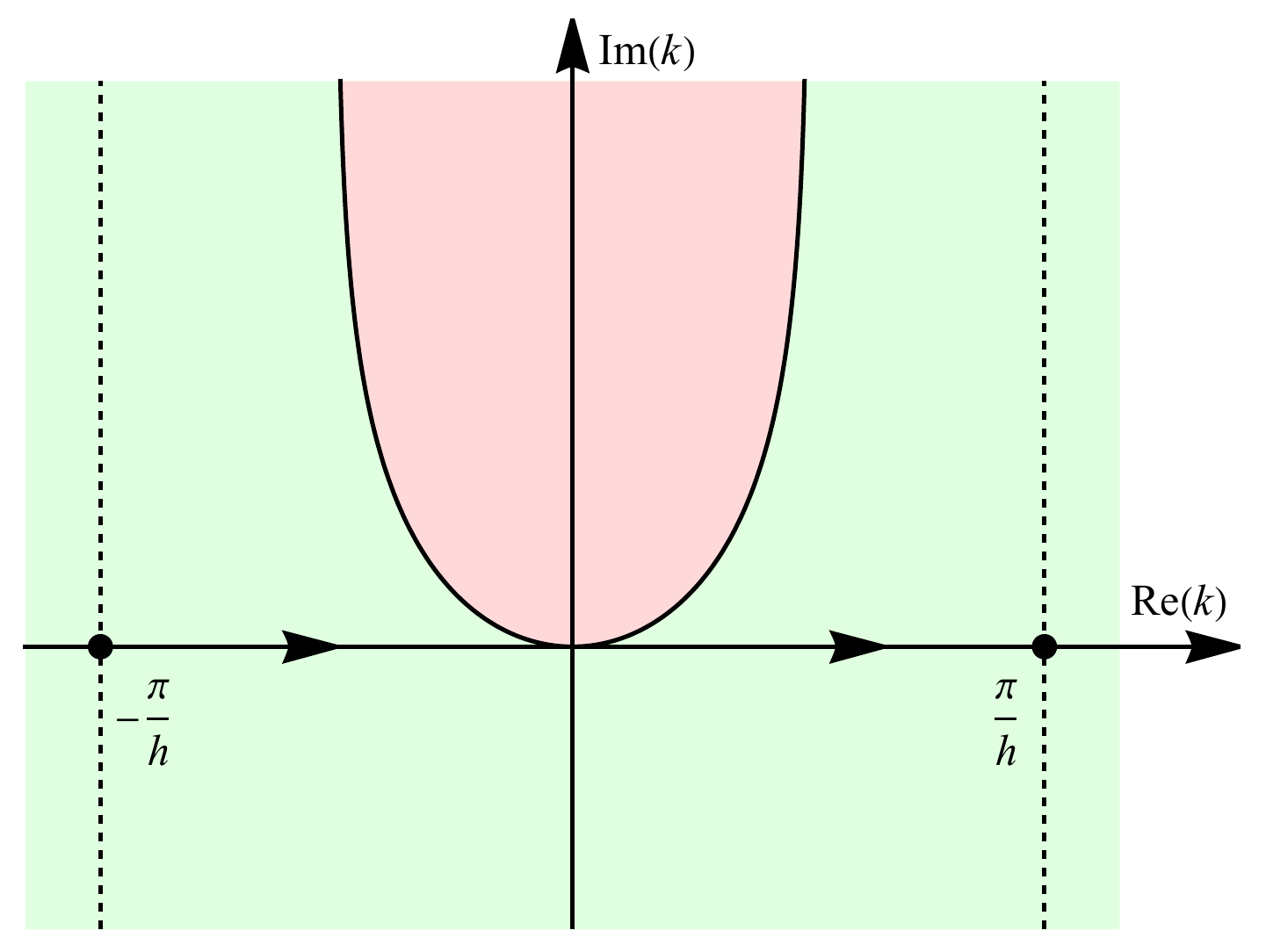}{The green region depicts where $\Re(-W) \leq 0$ and $e^{-W T}$ is decaying with the dispersion relation \eqref{advec2_back_W}.}{advec_backward_W}


For the first integral in \eqref{soln_advec2_backward} determined by the initial condition, we substitute the definition of $\hat{u}(k,0)$:
\begin{align}
		\frac{1}{2\pi} \int_{-\pi/h}^{\pi/h} e^{iknh} e^{-WT} \hat{u}(k,0)\,dk &= \sum_{m = 1}^{\infty} \left[ \frac{h}{2\pi} \int_{-\pi/h}^{\pi/h} e^{ik(n-m)h} e^{-WT} \,dk\right] u_m(0).
\end{align}
\no Using $z = e^{ikh}$, 
\begin{align}\nonumber
		\frac{h}{2\pi} \int_{-\pi/h}^{\pi/h} e^{ik(n-m)h} e^{-WT} \,dk &= \frac{e^{-cT/h}}{2\pi i} \oint_{|z|=1} z^{n-m-1} e^{c T/hz} \,dz = e^{-cT/h} \,\underset{z = 0}{\text{Res}} \left\{ z^{n-m-1} e^{c T/hz} \right\}\\
		&= e^{-cT/h} \left( \frac{c T}{h} \right)^{n - m} \frac{1}{(n - m) !}.
\end{align}
\no Hence,
\begin{align}
		\frac{1}{2\pi} \int_{-\pi/h}^{\pi/h} e^{iknh} e^{-WT} \hat{u}(k,0)\,dk = e^{-cT/h} \sum_{m = 1}^{\infty}  \left( \frac{c T}{h} \right)^{n - m} \frac{u_m(0)}{(n - m) !}.
\end{align}
	
We follow a similar procedure for the second integral of \eqref{soln_advec2_backward} by substituting the definition of $F_0(W,T)$. Denoting this integral by $I(n,t)$, 
\begin{align}\nonumber
		 I(n,t) &= \frac{c}{2\pi} \int_{-\pi/h}^{\pi/h} e^{ik(n-1)h} e^{-WT} F_{0}(W,T)\,dk = c \int_{0}^T \left[ \frac{1}{2\pi} \int_{-\pi/h}^{\pi/h} e^{ik(n-1)h} e^{-W(T-t)} \,dk \right] f_0(t) \, dt\\
		 &= \frac{c}{h(n-1)!} \int_{0}^T e^{-c(T-t)/h} \left( \frac{c (T-t)}{h} \right)^{n - 1} f_0(t) \, dt,
\end{align}
\no using the same steps as above. 
%
Thus \eqref{soln_advec2_backward} is rewritten as
\begin{align}
		\begin{split}
		u_n(T) &= e^{-cT/h} \sum_{m = 1}^{\infty}  \left( \frac{c T}{h} \right)^{n - m} \frac{u_m(0)}{(n - m) !}\,\, + \,\,I(n,T), ~~~n\geq 1.
		\label{soln_advec2_backward_nok1}
	\end{split}\end{align}
	
\no In what follows, we explore how this expression can be analytically continued from $n>0$ to $n<0$. To this end, we manipulate the expression \rf{soln_advec2_backward_nok1} to obtain a different representation which can be evaluated for negative $n$. 

The first term is rewritten using $1/(n-m)! = 1/\Gamma(n-m+1)$, which may be evaluated for $n<0$. Since $1/\Gamma(\alpha)$ has simple zeros at nonpositive integers $\alpha$, the first term does not contribute for $n<0$. We focus on the second term. Substituting $s = \frac{c (T-t)}{h}$ and Taylor expanding about $h = 0$ gives
\begin{align}\nonumber
	I(n,T) &= \frac{1}{(n-1)!} \int_{0}^{cT/h} e^{-s} s^{n - 1} f_0\left(T - \frac{h}{c}s\right) \, ds  \\ \nonumber
	&= \frac{1}{(n-1)!} \int_{0}^{cT/h} e^{-s} s^{n - 1} \sum_{\ell = 0}^{\infty} \frac{f_0^{(\ell)}\left(T\right) (-1)^{\ell}}{\ell !} \left(\frac{h}{c} \right)^{\ell} s^{\ell} \, ds \\
	&= \frac{1}{\Gamma(n)} \sum_{\ell = 0}^{\infty} \frac{f_0^{(\ell)}\left(T\right) (-1)^{\ell}}{\ell !} \left(\frac{h}{c} \right)^{\ell} \gamma\left(n + \ell,\frac{cT}{h}\right),
\end{align}
\no where 
\beq 
\gamma(a,y) = \int_0^y t^{a-1} e^{-t}\,dt
\eeq
\no is the lower incomplete gamma function \cite{dlmf}. The solution \eqref{soln_advec2_backward_nok1} is written as
\begin{align}\begin{split}
		u_n(T) &= e^{-cT/h} \sum_{m = 1}^{\infty}  \left( \frac{c T}{h} \right)^{n - m} \frac{u_m(0)}{(n - m) !} \,\,+\,\, \sum_{\ell = 0}^{\infty} \frac{f_0^{(\ell)}\left(T\right) (-1)^{\ell}}{\ell !} \left(\frac{h}{c} \right)^{\ell} \frac{\gamma\left(n + \ell,\frac{cT}{h}\right)}{\Gamma(n)},
		\label{soln_advec2_backward_nok2}
\end{split}\end{align}
\no for $n \geq 1$. This representation is only valid for $\ell \geq 1 - n$ and nonzero for $n \geq 1$ due to the ratio $\gamma\left(n + \ell,cT/h\right)/\Gamma(n)$. Recursively applying $\Gamma(k+1) = k\, \Gamma(k)$, 
\begin{align}
		\Gamma(a -b) = (a - b -1)(a - b -2)\cdot\ldots\cdot(-b)\,\Gamma(-b) = \frac{(-1)^a\, \Gamma(b+1)\,\Gamma(-b)}{\Gamma(b - a +1)}.
		\label{gammaprop1}
\end{align}

Using the power series of $\gamma(a,y)$ \cite{dlmf}, 
\begin{align}\nonumber
		\frac{\gamma\left(n + \ell,y\right)}{\Gamma(n)} &= \frac{1}{\Gamma(n)} \sum_{k=0}^{\infty} \frac{(-1)^k\,y^{n + \ell+k}\,\Gamma(n + \ell+k)}{k!\,\Gamma(n + \ell+k+1)} \\ \nonumber
		&= (-1)^{\ell} \sum_{k=0}^{\infty} \frac{y^{n + \ell+k}\,\Gamma( 1-n)}{k!\,\Gamma(1+n +k+ \ell)\,\Gamma(1-n-k-\ell)}\\
		&= \frac{(-1)^{\ell} \, y^{n + \ell - n -\ell}\,\Gamma( 1-n)}{(- n -\ell)!}
		= \frac{(-1)^{\ell}\, \Gamma( 1-n)}{\Gamma(1- n -\ell)},
\end{align}

\no where the sum collapses since $\Gamma(1+n+k+l)\Gamma(1-n-k-l)$ is infinite unless $n+k+l=0$. This representation is valid for $n \geq 0$ and nonzero for $\ell \leq -n$. For any $n \in \mathbb{Z}$, 
\begin{align}\nonumber 
u_n(T)&=\lim_{\alpha \rightarrow n} u_{\alpha}(T)\\
&= e^{-cT/h} \sum_{m = 1}^{\infty} \left( \frac{c T}{h} \right)^{n- m} \frac{u_m(0)}{(n- m) !} \,\,+\,\, \sum_{\ell = 0}^{\infty} \frac{f_0^{(\ell)}\left(T\right) (-1)^{\ell}}{\ell !} \left(\frac{h}{c} \right)^{\ell} \lim_{\alpha \rightarrow n} \frac{\gamma\left(\alpha + \ell,\frac{cT}{h}\right)}{\Gamma(\alpha)},
	\label{soln_advec2_backward_nok_lim}
\end{align}
\no where
\beq
\lim_{\alpha \rightarrow n} \frac{\gamma\left(\alpha + \ell,\frac{cT}{h}\right)}{\Gamma(\alpha)} = \begin{dcases} 
		\frac{\gamma\left(n + \ell,\frac{cT}{h}\right)}{\Gamma(n)}, &  \ell \geq 1 - n \quad \text{and} \quad n \geq 1, \\[5pt]
		\frac{(-1)^{\ell} \Gamma(1-n)}{\Gamma(1 - n - \ell)}, &  n \leq 0 \quad \text{and} \quad \ell \leq -n, \\[5pt]
		0, &  \text{otherwise,}
	\end{dcases}
\eeq 
\no so that
\beq  
I(n,T) = \begin{dcases} 
	 \sum_{\ell = 0}^{-n} \frac{f_0^{(\ell)}\left(T\right)}{\ell !} \left(\frac{h}{c} \right)^{\ell} \frac{ \Gamma(1-n)}{\Gamma(1 - n - \ell)}, &  n \leq 0, \\[5pt]
	 \sum_{\ell = 0}^{\infty} \frac{f_0^{(\ell)}\left(T\right) (-1)^{\ell}}{\ell !} \left(\frac{h}{c} \right)^{\ell}\, \frac{\gamma\left(n + \ell,\frac{cT}{h}\right)}{\Gamma(n)}, &  n \geq 1.
	\end{dcases}
\eeq 

Thus,  $u_n(T)$ with $n \geq 1$ is given by \eqref{soln_advec2_backward_nok2}, whereas for $n\leq 0$,   
\begin{align}\begin{split}
		u_{n}(T) &= I(n,T) =  \sum_{\ell = 0}^{n} \frac{f_0^{(\ell)}\left(T\right)}{\ell !} \left(\frac{h}{c} \right)^{\ell} \frac{ \Gamma(1-n)}{\Gamma(1 - n - \ell)}.
		\label{soln_advec2_backward_nok3}
	\end{split}\end{align}

\no Since \eqref{soln_advec2_backward_nok2} originates from \eqref{soln_advec2_backward}, we combine \eqref{soln_advec2_backward} and \eqref{soln_advec2_backward_nok3} to create an expression valid for $n \in \mathbb{Z}$ as
\begin{align}
		\begin{split}
		u_{n}(T) &= \frac{1}{2\pi} \int_{-\pi/h}^{\pi/h} e^{iknh} e^{-WT} \hat{u}(k,0)\,dk \,\,+ \,\,\frac{c}{2\pi} \int_{-\pi/h}^{\pi/h} e^{ik(n-1)h} e^{-WT} F_{0}\,dk \\
		&\quad\, + \sum_{\ell = 0}^{-n} \frac{f_0^{(\ell)}\left(T\right)}{\ell !} \left(\frac{h}{c} \right)^{\ell} \frac{ \Gamma(1-n)}{\Gamma(1 - n - \ell)},
		\label{soln_advec2_backward_nok_alln}
\end{split}\end{align}
\no where the integral terms only contribute for $n \geq 1$ and the sum only contributes for $n \leq 0$. 


\subsubsection{Continuum Limit}

For the continuous case, \eqref{advec2_prob}, 
\begin{equation}
		u(-x,T) = f_0\left(T+\frac{x}{c}\right), \quad x > 0,
		\label{soln_advec2_cont_classical}
\end{equation}
\no and the negative half-line solution \eqref{soln_advec2_backward_nok3} depends only on the boundary condition. The same is true in the semi-discrete case, see \rf{soln_advec2_backward_nok3}. We show that \eqref{soln_advec2_backward_nok3} limits to \eqref{soln_advec2_cont_classical} by taking $h \rightarrow 0$. 

Recursively applying $\Gamma(z+1) = z \Gamma(z)$, 
\beq 
f(n,\ell) = \frac{ \Gamma(1+n)}{\Gamma(1 + n - \ell)} = \prod_{p = 0}^{\ell - 1} (n-p) = \begin{dcases} 1, \quad &\ell = 0, \\ \sum_{p=0}^{\ell-1} a_p n^{p+1}, \quad &\ell \geq 1, \end{dcases}
\eeq 
\no and  $f(n,\ell)$ is a polynomial in $n$ of degree $\ell$ with leading coefficient $a_{\ell-1} = 1$. Hence, with $n\geq 0$, 
\begin{align}
		u_{-n}(T) 
		&= \sum_{\ell = 0}^{n} \frac{f_0^{(\ell)}\left(T\right)}{\ell !} \left[ \left(\frac{nh}{c} \right)^{\ell} + a_{\ell-2}\left(\frac{h}{c} \right) \left(\frac{nh}{c}\right)^{\ell-1} + \ldots \right].
\end{align}
\no In the continuum limit $h\rightarrow 0$ with $\lim_{h \rightarrow 0} n h = x$,
\begin{align}
u(-x,T)=
\sum_{\ell = 0}^{\infty} \frac{f_0^{(\ell)}\left(T\right)}{\ell !} \left(\frac{x}{c} \right)^{\ell}= f_0\left(T+\frac{x}{c}\right).
\end{align}

\subsubsection{Examples}
	
\begin{figure}[tb]
\begin{center}
\def \sc {0.475}
\begin{tabular}{cc}
\includegraphics[scale=\sc]{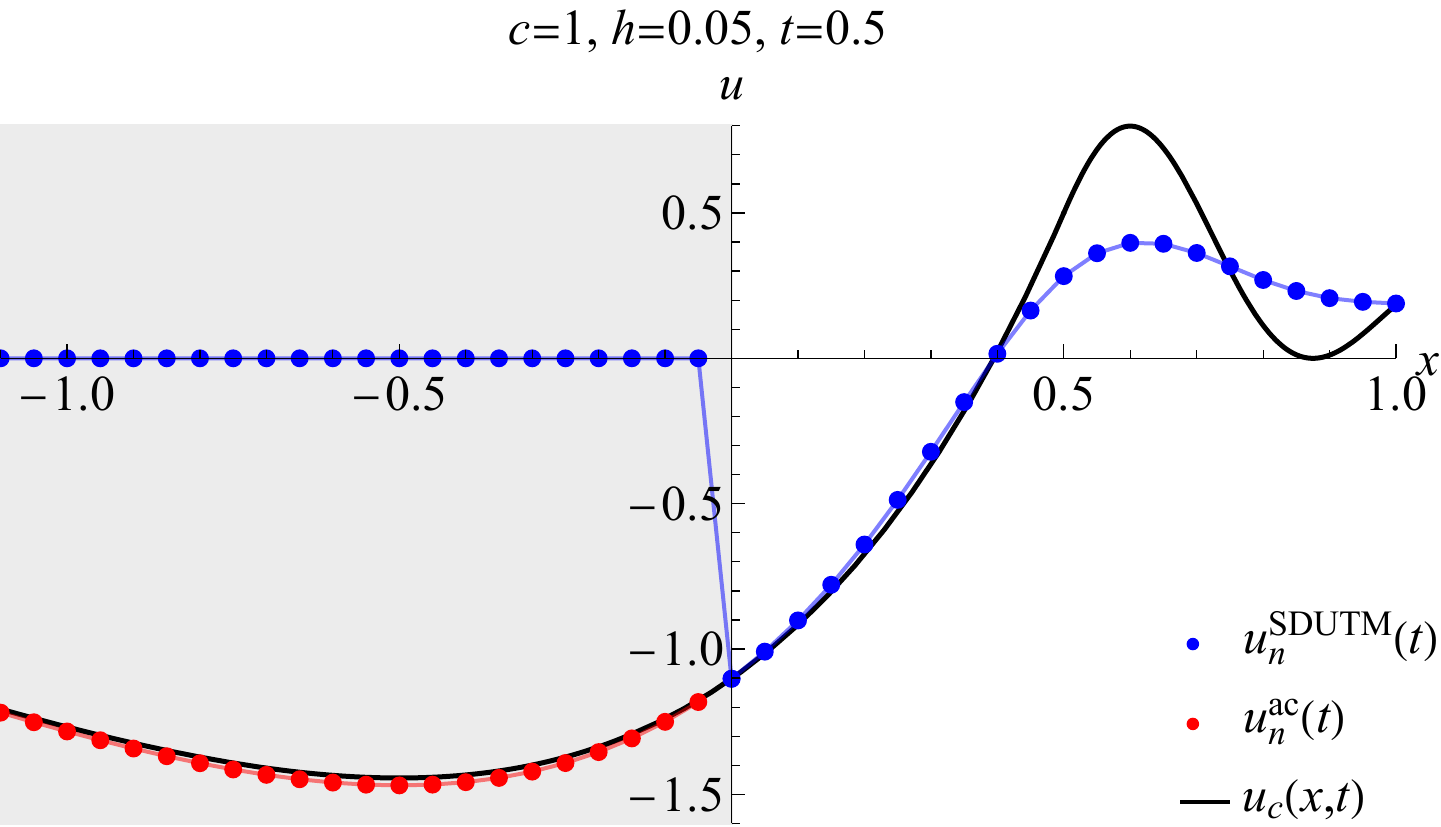} & \includegraphics[scale=\sc]{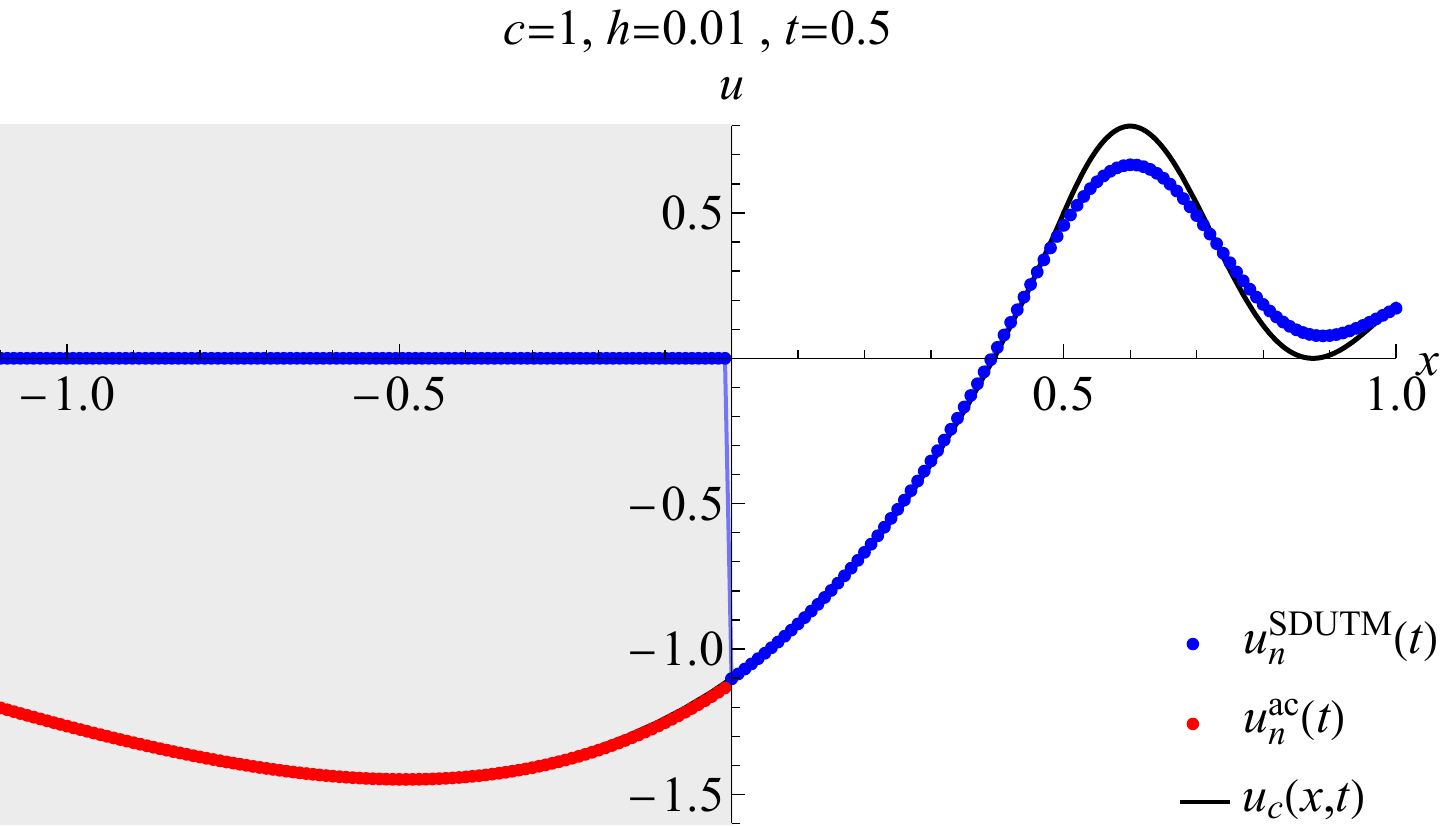}\\
(a) & (b)
\end{tabular}
\end{center}
\caption{(a) For $h=0.05$, the solution $u_n^{\text{ac}}(t)$ given by \rf{soln_advec2_backward_nok_alln} at $t=0.5$ and with $c=1$ obtained through analytic continuation, shown with the SDUTM solution $u_n^\text{SDUTM}(t)$ \rf{soln_advec2_backward_nok2} and the continuum solution $u_c(x,t)$ given by \rf{eqn:soln_adv}. (b) The same plot for $h=0.01$.
}
\label{advec2_backward_analytcont1}
\end{figure}

Figure \ref{advec2_backward_analytcont1} depicts the semi-discrete UTM solutions for $n \in \mathbb{Z}$ using \eqref{soln_advec2_backward} and \eqref{soln_advec2_backward_nok3} on the left and \eqref{soln_advec2_backward_nok_alln} on the right for the IBVP \rf{advec2_prob} with, for $x,t>0$,  
\beq 
u(x,0)=\phi(x) = \left(\frac{\sin(4 \pi x) + 1}{2}\right)e^{-2x},~~\mbox{and}~~
u(0,t) =f_0(t)= \frac{1}{2} + (1 - 2 \pi)t e^{-t}.
		\label{advec2_prob_ex}
\eeq 

\subsection{A second-order backward discretization}

\label{advec2_2backward_halfline_analytcont}


The backward, second-order discretization applied to \eqref{advec2_prob} is
\begin{equation}
		\dot{u}_n(t) = c \,\frac{-u_{n-2}(t) + 4 u_{n-1}(t) - 3u_{n}(t)}{2h},~~~n\geq 1, \, t>0. 
		\label{advec2_backward2}
\end{equation}
\no The details for the standard second-order centered discretization are presented in the Appendix, as they present peculiarities that are inherent to the equation and the discretization, but not the analytic continuation. Following the steps outlined in \cite{JC_SDUTM_HL, JC_SDUTM_FI}, \eqref{advec2_backward2} has the solution 
\begin{align} \nonumber
	u_n(T) &= \frac{1}{2\pi} \int_{-\pi/h}^{\pi/h} e^{iknh} e^{-WT} \hat{u}(k,0)\,dk + \frac{c}{4\pi} \int_{-\pi/h}^{\pi/h} e^{ik(n-1)h} e^{-WT} \left( 3 - e^{-ikh}\right) F_{0}(W,T) \,dk \\
	&\quad\, - \frac{c}{4\pi} \int_{-\pi/h}^{\pi/h} e^{ik(n-1)h} e^{-WT} \left[ \left(\frac{h }{c}\right) U_1 + \left(\frac{h^2}{2c^2}\right) U_2 \right]\,dk, ~~~n\geq 1, \; t>0, 
	\label{soln_advec2_backward2}
\end{align}
\no with dispersion relation
\begin{equation}
		W(k) = c \,\frac{e^{-2ikh} - 4 e^{-ikh} + 3}{2h},
		\label{W_advec2_backward2}
\end{equation}
\no and $F_0$ as before, see \rf{f0}. Further, 
\begin{align}
U_1(W,T)=\int_0^T e^{Wt}\dot f_0(t)dt, \qquad U_2(W,T)=\int_0^T e^{Wt}\ddot f_0(t)dt. 
\end{align} 
\no The dots on $f_0(t)$ on the right-hand sides denote time derivatives. To derive \rf{soln_advec2_backward2}, \eqref{advec2_prob} has been used to derive boundary conditions from the Dirichlet data, 
\begin{align}
u_x(0,t)=-\dot f_0(t)/c, \qquad u_{xx}(0,t)=\ddot f_0(t)/c^2. 
\end{align}
\no These derivative conditions are discretized using centered second-order accurate stencils (to retain the desired second-order accuracy), followed by multiplying by $\exp(Wt)$ and integrating:
\beq 
\frac{F_{1} - F_{-1}}{2h} = \frac{-U_1}{c}, \qquad	\frac{F_{1} - 2 F_0 + F_{-1}}{h^2} = \frac{U_2}{c^2}. 
\label{advec2_neumann_bc}
\eeq 
\no As $h \rightarrow 0$, the semi-discrete solution \eqref{soln_advec2_backward2} converges to \eqref{soln_advec2_cont}, as the semi-discrete solution correctly loses dependence on the derivative boundary conditions in the continuum limit. As before, one shows that \eqref{soln_advec2_backward2} gives $u_{-n}(T) = 0$ when evaluated for $n \in \mathbb{Z}^{+}$, similar to how \eqref{soln_advec2_cont} gives $u(-x,T) = 0$ for $x \in \mathbb{R}^{+}$. The analytical continuation of \eqref{soln_advec2_backward2} for negative values of $n$ is constructed below. 

We introduce
\beq 
B_m(n,T,g) = \frac{1}{4\pi} \int_{-\pi/h}^{\pi/h} e^{ik(n-m)h-WT} G\,dk, \quad\quad G(W,T) = \int_0^T e^{Wt} g(t)\,dt.
\eeq 
\no The semi-discrete solution is rewritten as
\begin{align}\begin{split}
		u_n(T) &= \frac{1}{2\pi} \int_{-\pi/h}^{\pi/h} e^{iknh} e^{-WT} \hat{u}(k,0)\,dk + 3c B_1(n,T,f_0) - c B_2(n,T,f_0) \\
		&\quad\, - h B_1(n,T,\dot{f_0}) - \frac{h^2}{2c} B_1(n,T,\ddot{f_0}).
		\label{soln_B_advec2_backward2}
\end{split}\end{align}

Leaving the initial-condition as is, we consider $B_m(n,T,g)$. We have
\begin{align*}
		B_m(n,T,g) &= \frac{1}{4\pi} \int_{-\pi/h}^{\pi/h} e^{ik(n-m)h-WT} \left[ \int_0^T e^{Wt} g(t)\,dt\right]\,dk = \int_0^T I_m(n,T,t) g(t)\,dt,
\end{align*}
\no with 
\begin{align}\nonumber
		I_m(n,T,t) &= \frac{1}{4\pi} \int_{-\pi/h}^{\pi/h} e^{ik(n-m)h-W(T-t)} \,dk= \frac{1}{4\pi} \oint_{|z|=1} z^{n-m} e^{-W(T-t)} \,\frac{dz}{i h z} \\\nonumber 
		&= \frac{e^{\frac{-3c(T-t)}{2h}}}{2 h}\, \underset{z = 0}{\text{Res}} \left\{ z^{n-m-1} \exp\left[\left(-z^{-2} + 4 z^{-1}\right)\frac{c(T-t)}{2h}\right] \right\}\\
		&= \frac{ e^{\frac{-3c(T-t)}{2h}}}{2 h} \sum_{k = 0}^{\frac{n-m}{2}} \frac{(-1)^k}{2^{2k} k! (n-m-2k)!} \left(\frac{2c(T-t)}{h}\right)^{n-m-k}. 
\end{align}

\no where a tedious calculation shows that there is no contribution if $n-m$ is odd. Using the substitution $s = 3c(T-t)/(2h)$ and expanding about $h = 0$, 
\begin{align}\nonumber
		B_m(n,T,g) &= \int_0^T \left[ \frac{e^{\frac{-3c(T-t)}{2h}}}{2 h} \sum_{k = 0}^{\frac{n-m}{2}} \frac{(-1)^k}{2^{2k} k! (n-m-2k)!} \left(\frac{2c(T-t)}{h}\right)^{n-m-k}\right] g(t)\,dt \\\nonumber
&\hspace*{-0.7in}= \frac{1}{3c} \sum_{k = 0}^{\frac{n-m}{2}} \frac{(-1)^k}{2^{2k} k! (n-m-2k)!} \left(\frac{4}{3}\right)^{n-m-k}\left[ \int_0^{\frac{3cT}{2h}} e^{-s} s^{n-m-k} g\left(T - \frac{2h}{3c}s\right)\,ds \right]  \\\nonumber
		&\hspace*{-0.7in}= \frac{1}{3c} \sum_{k = 0}^{\frac{n-m}{2}} \frac{(-1)^k}{2^{2k} k! (n-m-2k)!} \left(\frac{4}{3}\right)^{n-m-k} \sum_{p=0}^{\infty} \frac{g^{(p)}(T)(-1)^p}{p!} \left( \frac{2h}{3c}\right)^p \int_0^{\frac{3cT}{2h}} e^{-s} s^{n-m-k+p} \,ds  \\
		&\hspace*{-0.7in}= \frac{1}{3c} \sum_{k = 0}^{\frac{n-m}{2}} \frac{(-1)^k}{2^{2k} k! (n-m-2k)!} \left(\frac{4}{3}\right)^{n-m-k} \sum_{p=0}^{\infty} \frac{g^{(p)}(T)(-1)^p}{p!} \left( \frac{2h}{3c}\right)^p \gamma\left(n-m-k+p+1,\frac{3cT}{2h}\right).
\end{align}
\no Applying relations \eqref{gammaprop1} and the power series of the lower incomplete gamma function \cite{dlmf}, 
\beq\lim_{\alpha \rightarrow n} \frac{\gamma\left(\alpha-m-k+p+1,\frac{3cT}{2h}\right)}{\Gamma(\alpha-m-2k+1)} = \begin{dcases} 
		\frac{\gamma\left(n-m-k+p+1,\frac{3cT}{2h}\right)}{\Gamma(n-m-2k+1)},  &p \geq k - n + m, \, k \leq \tfrac{n - m}{2},\\[5pt]
		\frac{(-1)^{k+p}\Gamma(m+2k-n)}{\Gamma(m+k-p-n)},  &k \geq \tfrac{1 + n - m}{2}, \, p \leq m+k-n-1, \\[5pt]
		0,  &\text{otherwise}.
	\end{dcases}
\label{lim_gamma1}
\eeq
\no Collecting the summands in $S(n,T,k,p)$ and relaxing the upper bound on the summation over $k$ (again by adding zero terms), we write
\begin{align}\nonumber 
		B_m(n,T,g) &= \sum_{k = 0}^{\infty} \sum_{p=0}^{\infty} S(n,T,k,p) \\\nonumber 
		&= \left(\sum_{k = 0}^{\tfrac{n - m}{2}} + \sum_{k = \max\left(\tfrac{1 + n - m}{2},0\right)}^{\infty}\right) \left(\sum_{p=0}^{m+k-n-1} + \sum_{p=\max(k - n + m,0)}^{\infty}\right) S(n,T,k,p) \\
		&= \sum_{k = 0}^{\tfrac{n - m}{2}}\sum_{p=\max(k - n + m,0)}^{\infty} S(n,T,k,p) + \sum_{k = \max\left(\tfrac{1 + n - m}{2},0\right)}^{\infty} \sum_{p=0}^{m+k-n-1} S(n,T,k,p) ,
\end{align}
\no where the other terms are shown to vanish. For $n > 0$, the second of the double sums above vanishes, and the first can be rewritten in the integral form from \eqref{soln_advec2_backward2}. For $n<0$, the first pair vanishes. We let $n \rightarrow -n$ for $n>0$ from here on. Relaxing the bounds (re-introducing zero contributions), 
\begin{align}\nonumber 
		B_m(-n,T,g) &= \sum_{k = \max\left(\tfrac{1 - n - m}{2},0\right)}^{\infty}\sum_{p=0}^{m+k+n-1} S(-n,T,k,p) \\\nonumber 
		&\hspace*{-0.7in}= \frac{1}{3c} \sum_{k = 0}^{\infty} \sum_{p=0}^{\infty} \frac{(-1)^k}{2^{2k} k!} \left(\frac{4}{3}\right)^{-n-m-k} \frac{g^{(p)}(T)(-1)^p}{p!} \left( \frac{2h}{3c}\right)^p \frac{(-1)^{k+p}\Gamma(m+2k+n)}{\Gamma(m+k-p+n)}\\\nonumber 
		&\hspace*{-0.7in}= \frac{1}{3c} \left(\frac{3}{4}\right)^{n+m} \sum_{p=0}^{\infty} \frac{g^{(p)}(T)}{p!} \left( \frac{2h}{3c}\right)^p \sum_{k = 0}^{\infty}\frac{1}{2^{2k} k!} \left(\frac{3}{4}\right)^{k} \frac{\Gamma(m+2k+n)}{\Gamma(m+k-p+n)}\\
		&\hspace*{-0.7in}= \frac{1}{3c} \left(\frac{3}{4}\right)^{n+m} \sum_{p=0}^{\infty} \frac{g^{(p)}(T)}{p!} \left( \frac{2h}{3c}\right)^p \Gamma(m+n) \, \ptFq{2}{1}\left(\tfrac{m+n}{2},\tfrac{m+n+1}{2};m+n-p;\tfrac{3}{4}\right),
\end{align}
\no where $\ptFq{2}{1}(a,b;c;z) = \pFq{2}{1}(a,b;c;z)/\Gamma(c)$ is the regularized hypergeometric function \cite{dlmf}. In summary, the solution for $n>0$ is given by \eqref{soln_advec2_backward2}, while the analytic continuation for negative values of $n$,
\begin{align}\begin{split}
		u_{n}(T) &= 3c B_1(n,T,f_0) - c B_2(n,T,f_0) - h B_1(n,T,\dot{f_0}) - \frac{h^2}{2c} B_1(n,T,\ddot{f_0}),~~n\leq 0, 
		\label{soln_AC_advec2_backward2}
\end{split}\end{align}
\no since the initial-condition integral vanishes. A careful limit calculation shows that the analytic continuation obtained converges to $u(-x,T) = f_0(T+x/c)$ as $h\rightarrow 0$. 

\subsubsection{Examples}
	
\begin{figure}[tb]
\begin{center}
\def \sc {0.475}
\begin{tabular}{cc}
\includegraphics[scale=\sc]{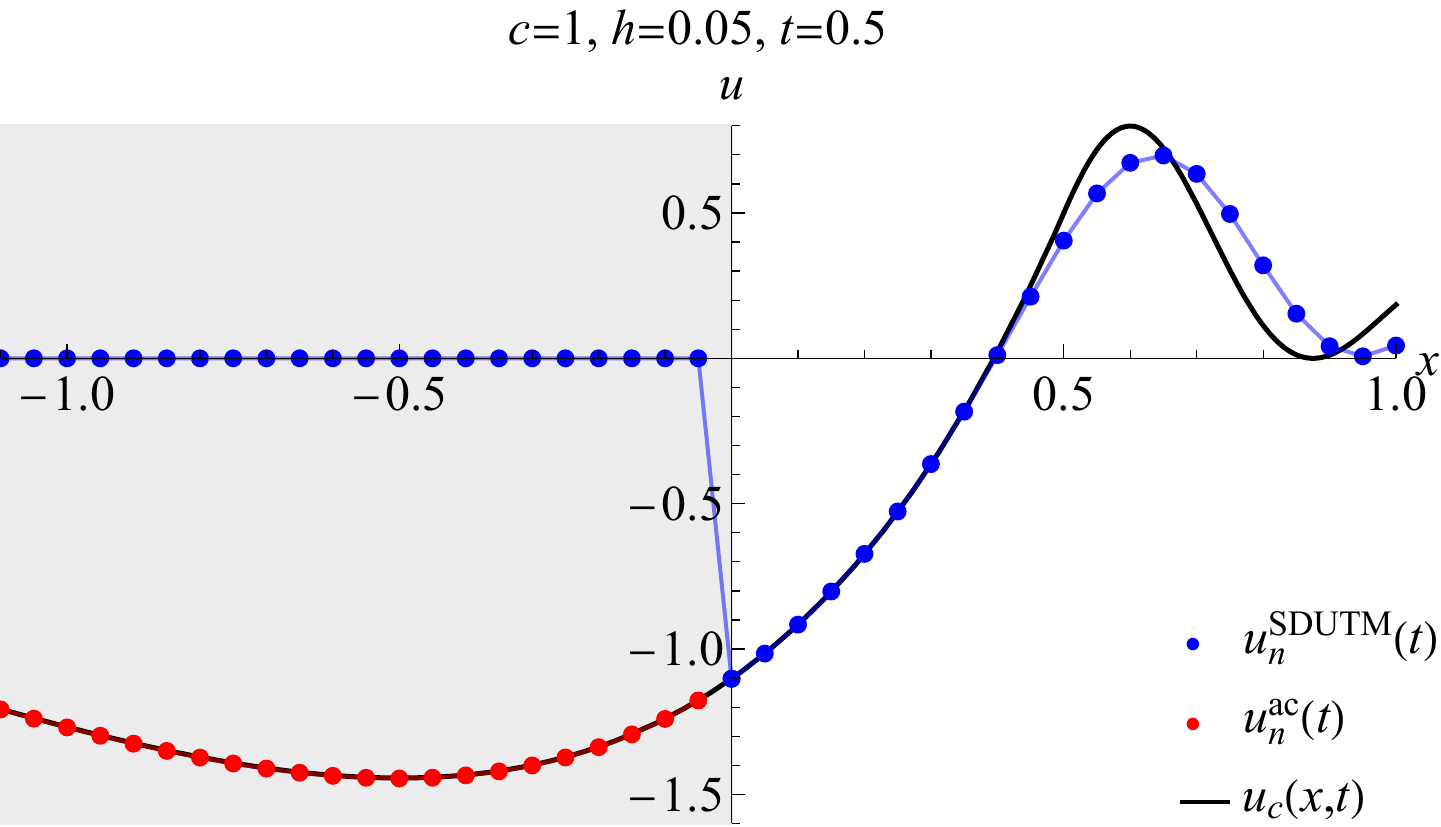} & \includegraphics[scale=\sc]{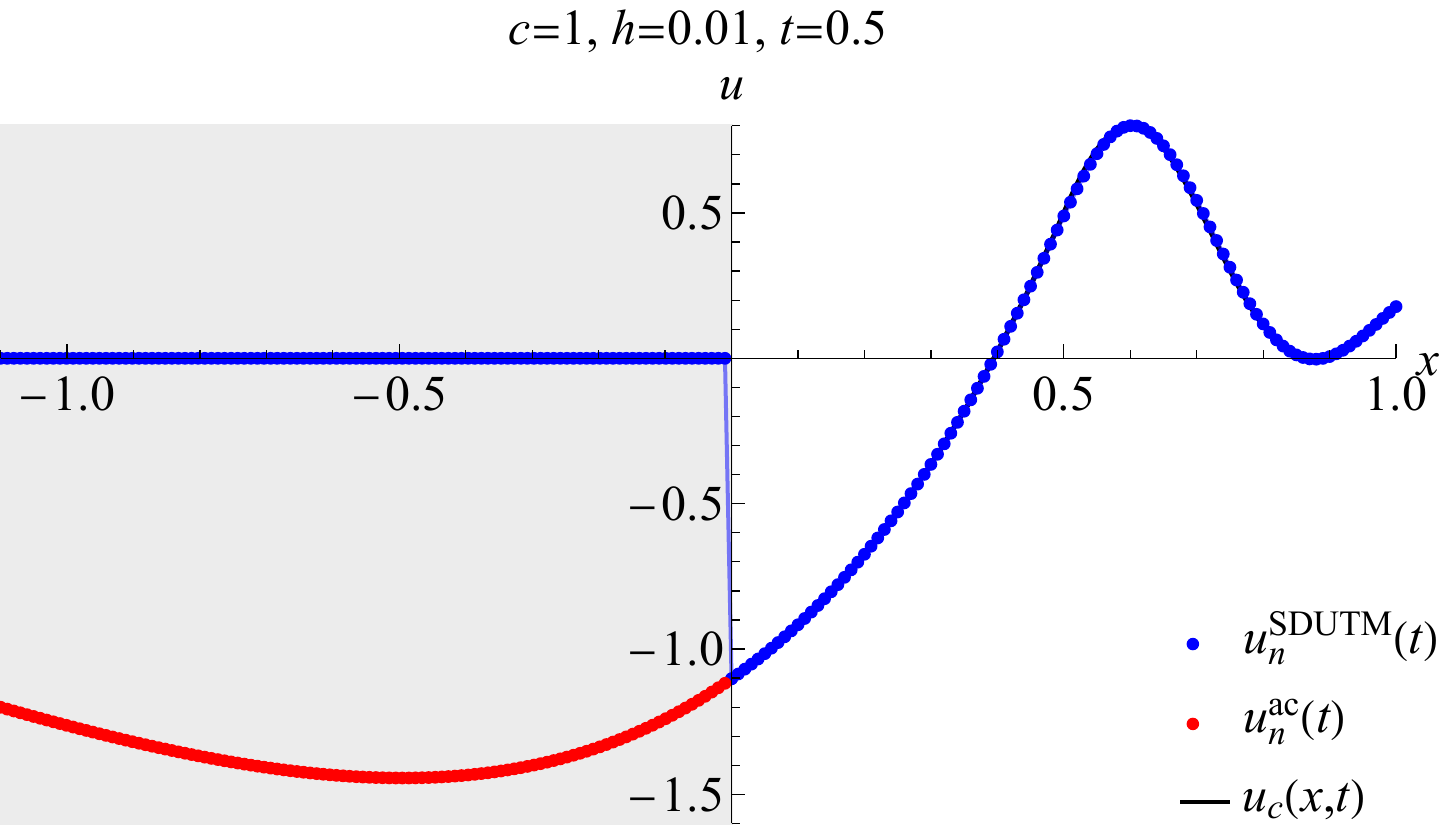}\\
(a) & (b)
\end{tabular}
\end{center}
\caption{(a) For $h=1/20$, the solution $u_n^{\text{ac}}(t)$ given by \rf{soln_AC_advec2_backward2} at $t=0.5$ and with $c=1$ obtained through analytic continuation, shown with the SDUTM solution $u_n^\text{SDUTM}(t)$ \rf{soln_advec2_backward2} and the continuum solution $u_c(x,t)$ given by \rf{eqn:soln_adv}. (b) The same plot for $h=1/150$.
}
\label{advec2_backward_analytcont1_2nd}
\end{figure}

Figure \ref{advec2_backward_analytcont1_2nd} illustrates the semi-discrete solutions for $n \in \mathbb{Z}$ for the IBVP \rf{advec2_backward2} with boundary data given by \rf{advec2_prob_ex}.

\section{The discretized heat equation, Dirichlet boundary conditions, centered discretization} \label{sec:sdheat}



\label{heat_cent_halfline_analytcont}

We examine a discretization of \rf{fig:heatdirichletbc}\footnote{The calculations for the Neumann problem are similar, and the details are left to the reader. The main difference with the Dirichlet problem is the role played by the discretization of the Neumann conditions.}. 
Discretizing the spacial derivative $u_{xx}$ with the standard centered stencil gives 
\beq\la{heatdiscretecentered}
\dot{u}_n(t) = \frac{u_{n-1}(t) - 2u_{n}(t) + u_{n+1}(t)}{h^2},
 \eeq 
\no with dispersion relation 
\begin{equation}
	    W(k) = \frac{2 - e^{ikh} - e^{-ikh}}{h^2}.
	    \label{heat_W_cent2}
\end{equation}
\no	The solution \cite{JC_SDUTM_HL} to this semi-discrete problem is
\begin{align}\begin{split}
		u_n(T) &= \frac{1}{2 \pi} \int_{-\pi/h}^{\pi/h} e^{iknh} e^{-WT}\left[ \hat{u}(k,0) - \hat{u}(-k,0)\right] \,dk  - \frac{i}{ \pi h} \int_{-\pi/h}^{\pi/h} e^{iknh} e^{-WT} \sin(kh) F_{0} \,dk, ~~~n>0,
	\label{soln_heat_centered}
\end{split}\end{align}
\no where $\hat u(k,0)$ and $F_0(W,T)$ are defined in \rf{eqn:SDu0hat} and \rf{f0} respectively. When evaluated for $n\leq 0$, this representation of the solution gives $u_{-n}(T) = -u_{n}(T)$ for $n \geq 1$ and $u_0(T) = 0$. In the continuum limit, \eqref{soln_heat_centered} converges to
\begin{align}
		u(x,T) &= \frac{1}{2 \pi} \int_{-\infty}^{\infty} e^{ikx} e^{-\tilde{W}T}\left[ \hat{u}(k,0) - \hat{u}(-k,0)\right]\, dk - \frac{i}{\pi} \int_{-\infty}^{\infty} e^{ikx} e^{-\tilde{W}T} k \tilde F_0 \,dk,
		\label{soln_heat_cont}
	\end{align}
\no with $\tilde{W} = k^2$ and $u(-x,T) = - u(x,T)$ for $x > 0$.

We leave the first integral term as is, and introduce 
\beq 
B(n,a) := \frac{1}{2 \pi i} \oint_{|z|=1} \exp\left[\frac{a}{2}\left(z + \frac{1}{z} \right)\right] \frac{1}{z^{n + 1}}\, dz \,\equiv \,\sum_{\ell = 0}^{\infty} \frac{1}{\ell! \, \Gamma(\ell+n+1)}\left(\frac{a}{2}\right)^{2\ell + n},
\eeq 
\no the modified Bessel function of the first kind \cite{dlmf} with vital properties 
\beq\la{besselprop}
B(-n,a) = B(n,a) \,\text{ for }\, n \in \mathbb{Z}, \quad \text{ and }\quad B(n,a) = \frac{2(n+1)}{a} B(n+1,a) + B(n+2,a).
\eeq 
\no For the second term of \eqref{soln_heat_centered}, we have
\begin{align}
		\frac{i}{ \pi h} \int_{-\pi/h}^{\pi/h} e^{iknh} e^{-WT} \sin(kh) F_{0} \,dk &= \frac{i}{ \pi h} \int_0^T J(n,T-t) f_0(t) \, dt,
\end{align}
\no	with 
\begin{align}
		J(n,T) &=  \int_{-\pi/h}^{\pi/h} e^{iknh} e^{-WT} \sin(kh) \,dk
= \frac{\pi}{i h} I_1(n,T,-1) - \frac{\pi}{i h} I_1(n,T,1),
\end{align}
\no where
\begin{align}
	    I_1(n,T,m) &= \frac{h}{2 \pi} \int_{-\pi/h}^{\pi/h} e^{ik(n-m)h} e^{-WT} \,dk = \frac{1}{2 \pi i} \oint_{|z|=1} z^{n-m-1} \exp\left[-\left(\frac{2 - z - z^{-1}}{h^2}\right)T\right]\notag\\ 
	    &= e^{-2T/h^2} B\left(m-n,\frac{2T}{h^2}\right),
	    \label{I1_heat}
\end{align}
\no	so that, using \rf{besselprop}, 
\begin{align}
		J(n,T) &= \frac{\pi e^{-2T/h^2} }{i h} \left[  B\left(-1-n,\frac{2T}{h^2}\right) - B\left(1-n,\frac{2T}{h^2}\right)\right]
		= \frac{- n h \pi e^{-2T/h^2} }{iT} B\left(n,\frac{2T}{h^2}\right).
\end{align}

This allows for a rewrite of the second term and 
\eqref{soln_heat_centered} becomes
\begin{align}\begin{split}
		\hspace{-20pt}u_n(T) &= \frac{1}{2 \pi} \int_{-\pi/h}^{\pi/h} e^{iknh} e^{-WT}\left[ \hat{u}(k,0) - \hat{u}(-k,0)\right] \,dk  + n \int_0^T \frac{e^{-2(T-t)/h^2}}{T-t}  B\left(n,\frac{2(T-t)}{h^2}\right) f_0(t) \, dt.
	\label{soln_heat_centered_nok1}
\end{split}\end{align}
\no	As before, \eqref{soln_heat_centered_nok1} gives $u_{-n}(T) = - u_{n}(T)$ for $n \geq 1$ and $u_0(T) = 0$.
	
We denote the second term of \eqref{soln_heat_centered_nok1} by $K(n,T)$. Using the transformation $s = {2(T-t)}/{h^2}$ and a Taylor series expansion about $h = 0$,
\begin{align}\nonumber
		K(n,T) 
		&= n \int_0^{2T/h^2} \frac{e^{-s}}{s} B\left(n,s\right) f_0\left(T - \frac{h^2}{2} s\right) \, ds \\\nonumber 
		&= n \sum_{p = 0}^{\infty} \frac{f_0^{(p)}(T) (-1)^p}{p!} \left( \frac{h^2}{2}\right)^p \int_0^{2T/h^2} e^{-s} B\left(n,s\right) s^{p-1} \, ds \\\nonumber 
		&= n \sum_{p = 0}^{\infty} \frac{f_0^{(p)}(T) (-1)^p}{p!} \left( \frac{h^2}{2}\right)^p \sum_{\ell = 0}^{\infty} \frac{1}{\ell! \, \Gamma(\ell+n+1)\,2^{2\ell + n}} \int_0^{2T/h^2} e^{-s} s^{2\ell + n + p -1}  \, ds \\
		&= \frac{n}{2^n} \sum_{p = 0}^{\infty} \frac{f_0^{(p)}(T) (-1)^p}{p!} \left( \frac{h^2}{2}\right)^p \sum_{\ell = 0}^{\infty} \frac{1}{\ell! \, \Gamma(\ell+n+1)\,2^{2\ell}} \,\gamma\left( 2\ell + n + p, \frac{2T}{h^2}\right).
\end{align} 

To evaluate at any $n \in \mathbb{Z}$, we consider the limit
\begin{align}
		 K(n,T) &= \frac{n}{2^{n}} \sum_{p = 0}^{\infty} \frac{u^{(p)}(T) (-1)^p}{p!} \left( \frac{h^2}{2}\right)^p \sum_{\ell = 0}^{\infty} \frac{1}{\ell! \,2^{2\ell }} \lim_{\alpha\rightarrow n} \frac{\gamma\left( 2\ell + \alpha + p, \frac{2T}{h^2}\right)}{ \Gamma(\ell+\alpha+1)},
\end{align}
\no where \eqref{gammaprop1} and the Taylor series of the incomplete gamma function \cite{dlmf} with $y = 2T/h^2$ give 
\begin{align}\nonumber
		\frac{\gamma\left( 2\ell + n + p, y\right)}{\Gamma(\ell+n+1)} &= \frac{1}{\Gamma(\ell+n+1)} \sum_{k=0}^{\infty} \frac{(-1)^k\,y^{2\ell + n + p+k}}{k!\,\Gamma(1+2\ell + n + p+k)}\cdot\Gamma(2\ell + n + p+k) \\\nonumber 
		&= \frac{\Gamma(1-n-\ell-1)}{(-1)^{\ell+1}\,\Gamma(1-n)\,\Gamma(n)} \sum_{k=0}^{\infty} \frac{(-1)^k\,y^{2\ell + n + p+k}}{k!\,\Gamma(1+2\ell + n + p+k)}\cdot \frac{(-1)^{2\ell+p+k}\,\Gamma(1-n)\,\Gamma(n)}{\Gamma(1 - 2\ell - n - p-k)} \\
		&= (-1)^{\ell+p+1}\, \Gamma(-n-\ell) \cdot \frac{y^{2\ell + n + p-2\ell - n - p}}{(-2\ell - n - p)!}
		=\frac{  (-1)^{\ell+p+1}\, \Gamma(-n-\ell) }{\Gamma(1-2\ell - n - p)},
\end{align}
\no so that
\beq 
\lim_{\alpha\rightarrow n} \frac{\gamma\left( 2\ell + \alpha + p, \frac{2T}{h^2}\right)}{ \Gamma(\ell+\alpha+1)} = \begin{dcases} 
		\frac{\gamma\left( 2\ell + n + p, \frac{2T}{h^2}\right)}{\Gamma(\ell+n+1)}, &  2\ell + p \geq 1-n \quad \text{and} \quad \ell \geq -n, \\[5pt]
		\frac{(-1)^{\ell+p+1} \Gamma(-n-\ell)}{\Gamma(1-n - 2\ell - p)}, &  \ell \leq -n-1 \quad \text{and} \quad 2\ell + p \leq -n, \\[5pt]
		0, &  \text{otherwise.}
	\end{dcases}
\eeq

For brevity, define
\beq 
S(n,T,\ell,p) = \frac{f_0^{(p)}(T) (-1)^p}{p!} \left( \frac{h^2}{2}\right)^p \frac{1}{\ell! \,2^{2\ell }} \lim_{\alpha\rightarrow n} \frac{\gamma\left( 2\ell + \alpha + p, \frac{2T}{h^2}\right)}{ \Gamma(\ell+\alpha+1)},
\eeq
\no so that splitting the $\ell$ and $p$ sums gives
\begin{align}\nonumber
K(n,T) &= \frac{n}{2^{n}} \sum_{p = 0}^{\infty} \sum_{\ell = 0}^{\infty} S(n,T,\ell,p)\\\nonumber  
		&= \frac{n}{2^{n}} \sum_{\ell = 0}^{-n-1} \sum_{p = 0}^{-n-2\ell} S(n,T,\ell,p) + \frac{n}{2^{n}} \sum_{\ell = \max(-n,0)}^{\infty} \sum_{p = 0}^{-n-2\ell} S(n,T,\ell,p) \\
		&\quad\, + \frac{n}{2^{n}} \sum_{\ell = 0}^{-n-1}  \sum_{p = \max(-n-2\ell+1,0)}^{\infty} S(n,T,\ell,p) + \frac{n}{2^{n}} \sum_{\ell = \max(-n,0)}^{\infty} \sum_{p = \max(-n-2\ell+1,0)}^{\infty} S(n,T,\ell,p).
\end{align}
\no Consider $n < 0$, so that the $\ell$-indexed sum from the second pair of sums begins at $\ell = -n$, from which the upper bound of the $p$-indexed sum is $-n - 2\ell = n < 0$. Since the starting index is $p = 0$, this second pair of sums does not contribute for $n < 0$. Now consider $n \geq 0$, so that $\ell = 0$ is the starting index for the first sum, from which $-n - 2\ell = -n \leq 0$ is the upper bound of the $p$-indexed sum. Thus, the second pair of sums does not contribute for any $n \in \mathbb{Z}$. The third pair of sums likewise vanishes for all $n \in \mathbb{Z}$, since $S(n,T,\ell,p) = 0$ for these ranges of $\ell$ and $p$, regardless of $\lim_{\alpha \rightarrow n} \gamma\left( 2\ell + \alpha + p,\, 2T/h^2\right) /\Gamma(\ell+\alpha+1)$. We have
\begin{align}\nonumber 
		K(n,T)&= \frac{-n}{2^{n}} \sum_{\ell = 0}^{-n-1} \frac{(-1)^{\ell}}{\ell! \,2^{2\ell }} \sum_{p = 0}^{-n-2\ell} \frac{f_0^{(p)}(T)}{p!} \left( \frac{h^2}{2}\right)^p \frac{\Gamma(-n-\ell)}{\Gamma(1-n - 2\ell - p)} \\
		&\quad\, + \frac{n}{2^{n}} \sum_{\ell = \max(-n,0)}^{\infty} \frac{1}{\ell! \,2^{2\ell }} \sum_{p = \max(-n-2\ell+1,0)}^{\infty} \frac{f_0^{(p)}(T) (-1)^p}{p!} \left( \frac{h^2}{2}\right)^p \frac{\gamma\left( 2\ell + n + p, \frac{2T}{h^2}\right)}{ \Gamma(\ell+n+1)}.
\end{align}
\no In the first coupled sum, we relax the upper bound of the $p$-indexed sum (adding zeros), allowing us to interchange the sums: 
\begin{align}\nonumber 
		&\frac{-n}{2^{n}} \sum_{\ell = 0}^{-n-1} \frac{(-1)^{\ell}}{\ell! \,2^{2\ell }} \sum_{p = 0}^{\infty} \frac{f_0^{(p)}(T)}{p!} \left( \frac{h^2}{2}\right)^p \frac{\Gamma(-n-\ell)}{\Gamma(1-n - 2\ell - p)}\\ 
		&= \frac{-n}{2^{n}} \sum_{p = 0}^{\infty} \frac{f_0^{(p)}(T)}{p!} \left( \frac{h^2}{2}\right)^p \sum_{\ell = 0}^{-n-1} \frac{(-1)^{\ell}}{\ell! \,2^{2\ell }} \frac{\Gamma(-n-\ell)}{\Gamma(1-n - 2\ell - p)}
		= \frac{-n}{2^{n}} \sum_{p = 0}^{\infty} f_0^{(p)}(T) \left( \frac{h^2}{2}\right)^p L(n,p),
\end{align}
\no with
\begin{align}\nonumber 
L(n,p) &= \frac{1}{p!}\sum_{\ell = 0}^{{(-n-p)}/{2}} \frac{(-1)^{\ell}}{\ell! \,2^{2\ell }} \frac{\Gamma(-n-\ell)}{\Gamma(1-n - 2\ell - p)} = \frac{2^{n+p} (-1)^{n+p}\Gamma(1-n)\,\Gamma(n)\, \Gamma(2p-n-p)}{p! \,\Gamma(1-p)\,\Gamma(2p)\,\Gamma(1-n-p)}\\
&= \frac{2^{n+p+1}\,\Gamma(p-n) }{\Gamma(1-n-p)\,\Gamma(2p+1)},
\end{align} 
\no where the last equality follows from \eqref{gammaprop1}. Substituting $L(n,p)$ and truncating the $p$-indexed sum from $p = 0$ to $p = -n$, 
\begin{align}\nonumber 
K(n,T)&= - 2 n \sum_{p = 0}^{-n} \frac{f_0^{(p)}(T) \, h^{2p}\,\Gamma(p-n)}{\Gamma(1-n-p)\,\Gamma(2p+1)} \\
		&\!\!\!\!\!\!\!\!\!\!\quad\, + \frac{n}{2^{n}} \sum_{\ell = \max(-n,0)}^{\infty} \frac{1}{\ell! \,2^{2\ell }} \sum_{p = \max(-n-2\ell+1,0)}^{\infty} \frac{f_0^{(p)}(T) (-1)^p}{p!} \left( \frac{h^2}{2}\right)^p \frac{\gamma\left( 2\ell + n + p, \frac{2T}{h^2}\right)}{ \Gamma(\ell+n+1)},
\end{align}
\no so that \eqref{soln_heat_centered_nok1} is written as
\begin{align}\nonumber 
u_n(T) &= \frac{1}{2 \pi} \int_{-\pi/h}^{\pi/h} e^{iknh} e^{-WT}\left[ \hat{u}(k,0) - \hat{u}(-k,0)\right] \,dk  \,\,-\,\, 2 n \sum_{p = 0}^{-n} \frac{f_0^{(p)}(T) \, h^{2p}\,\Gamma(p-n)}{\Gamma(1-n-p)\,\Gamma(2p+1)} \\
		&\quad\, \,\,+\,\, \frac{n}{2^{n}} \sum_{\ell = \max(-n,0)}^{\infty} \frac{1}{\ell! \,2^{2\ell }} \sum_{p = \max(-n-2\ell+1,0)}^{\infty} \frac{f_0^{(p)}(T) (-1)^p}{p!} \left( \frac{h^2}{2}\right)^p \frac{\gamma\left( 2\ell + n + p, \frac{2T}{h^2}\right)}{ \Gamma(\ell+n+1)}.
	\label{soln_heat_centered_nok2}
\end{align} 
	
For $n>0$, the second term above does not contribute and the remaining terms are rewritten as in the original representation \eqref{soln_heat_centered}. For $n \leq 0$, 
\begin{align}
u_{n}(T) &= \frac{1}{2 \pi} \int_{-\pi/h}^{\pi/h} e^{iknh} e^{-WT}\left[ \hat{u}(k,0) - \hat{u}(-k,0)\right] \,dk  \,\,-\,\, 2 n \sum_{p = 0}^{-n} \frac{f_0^{(p)}(T) \, h^{2p}\,\Gamma(p-n)}{\Gamma(1-n-p)\,\Gamma(2p+1)} \notag\\
		&\quad\, \,\,+\,\, \frac{n}{2^{n}} \sum_{\ell =-n}^{\infty} \frac{1}{\ell! \,2^{2\ell }} \sum_{p = \max(-n-2\ell+1,0)}^{\infty} \frac{f_0^{(p)}(T) (-1)^p}{p!} \left( \frac{h^2}{2}\right)^p \frac{\gamma\left( 2\ell + n + p, \frac{2T}{h^2}\right)}{ \Gamma(\ell+n+1)} \notag\\[10pt]
		&= \frac{-1}{2 \pi} \int_{-\pi/h}^{\pi/h} e^{-iknh} e^{-WT}\left[ \hat{u}(k,0) - \hat{u}(-k,0)\right] \,dk  \,\,-\,\, 2 n \sum_{p = 0}^{-n} \frac{f_0^{(p)}(T) \, h^{2p}\,\Gamma(p-n)}{\Gamma(1-n-p)\,\Gamma(2p+1)} \notag\\
		&\quad\, \,\,+\,\, \frac{n}{2^{-n}} \sum_{\ell =0}^{\infty} \frac{1}{\ell! \,2^{2\ell}} \sum_{p = 0}^{\infty} \frac{f_0^{(p)}(T) (-1)^p}{p!} \left( \frac{h^2}{2}\right)^p \frac{\gamma\left( 2\ell - n + p, \frac{2T}{h^2}\right)}{ \Gamma(\ell-n+1)}, 
\end{align}
\no or
\beq 
	u_{n}(T) = -2 n \sum_{p = 0}^{-n} \frac{f_0^{(p)}(T) \, h^{2p}\,\Gamma(p-n)}{\Gamma(1-n-p)\,\Gamma(2p+1)} \,\,-\,\, u_{-n}(T), ~~~~n\leq 0. 
\label{soln_heat_centered_negn}
\eeq 

To recover the boundary condition at $n = 0$, we extract the first term of the sum above:
\begin{align}\nonumber 
u_{n}(T) &= 2 n \left[ \frac{f_0(T)}{n} + \sum_{p = 1}^{-n} \frac{f_0^{(p)}(T) \, h^{2p}\,\Gamma(p-n)}{\Gamma(1-n-p)\,\Gamma(2p+1)}\right] - u_{-n}(T)\\
&= 2 f_0(T) \,+\, 2 n \sum_{p = 1}^{-n} \frac{f_0^{(p)}(T) \, h^{2p}\,\Gamma(p-n)}{\Gamma(1-n-p)\,\Gamma(2p+1)} - u_{-n}(T),
\end{align}
\no so that for $n = 0$, $u_{0}(T) =  2 f_0(T)  - u_0(T) = f_0(T)$.

\subsection{Continuum Limit}

The analytic extension \eqref{heat_cont_full} of the continuous solution to the IBVP \rf{heatdirichletibvp} is
	\begin{equation}
		u(x,T) = 2 \sum_{p = 0}^{\infty} \frac{f_0^{(p)}(T)}{(2p)!} x^{2p} - u(-x,T), ~~~x>0. 
		\label{soln_cont_negx}
	\end{equation}

\no It is possible to recover \eqref{soln_cont_negx} from the continuum limit of \eqref{soln_heat_centered_negn} by noting that
\beq 
f(n,p) = \frac{n\Gamma(p+n)}{\Gamma(1+n-p)} = \prod_{\ell = 0}^{p-1}(n-\ell)(n+\ell) = \begin{dcases} 1, \quad &p = 0 \\ \sum_{\ell=1}^{p} a_{\ell} \,n^{2 \ell}, \quad &p \geq 1, \end{dcases}
\eeq 
\no and $f(n,p)$ is a polynomial in $n$ of degree $2p$ with leading coefficient $a_{p} = 1$. Hence,
\begin{align}\nonumber 
		u_{n}(T) &= 2 \sum_{p = 0}^{-n} \frac{f_0^{(p)}(T) \, h^{2p}}{(2p)!}\,f(-n,p) - u_{-n}(T)\\
		&= 2 \sum_{p = 0}^{-n} \frac{f_0^{(p)}(T) }{(2p)!}\left[  (nh)^{2p} +  a_{p-1} h^{2} (nh)^{2p-2} + \ldots \right] - u_{-n}(T),
\end{align}
\no which gives the desired continuum limit. 

\subsubsection{Examples}
	
We solve the IBVP
\begin{equation} \case{1}{
    u_t = u_{xx}, &  x > 0,\, t > 0, \\
    u(x,0) = 3 x e^{-x},& x > 0,\\
    u(0,t) = \sin(4 \pi t), & t > 0,
	\label{heat_prob1}
}\end{equation}

\begin{figure}[tb]
\begin{center}
\def \sc {0.475}
\begin{tabular}{cc}
\includegraphics[scale=\sc]{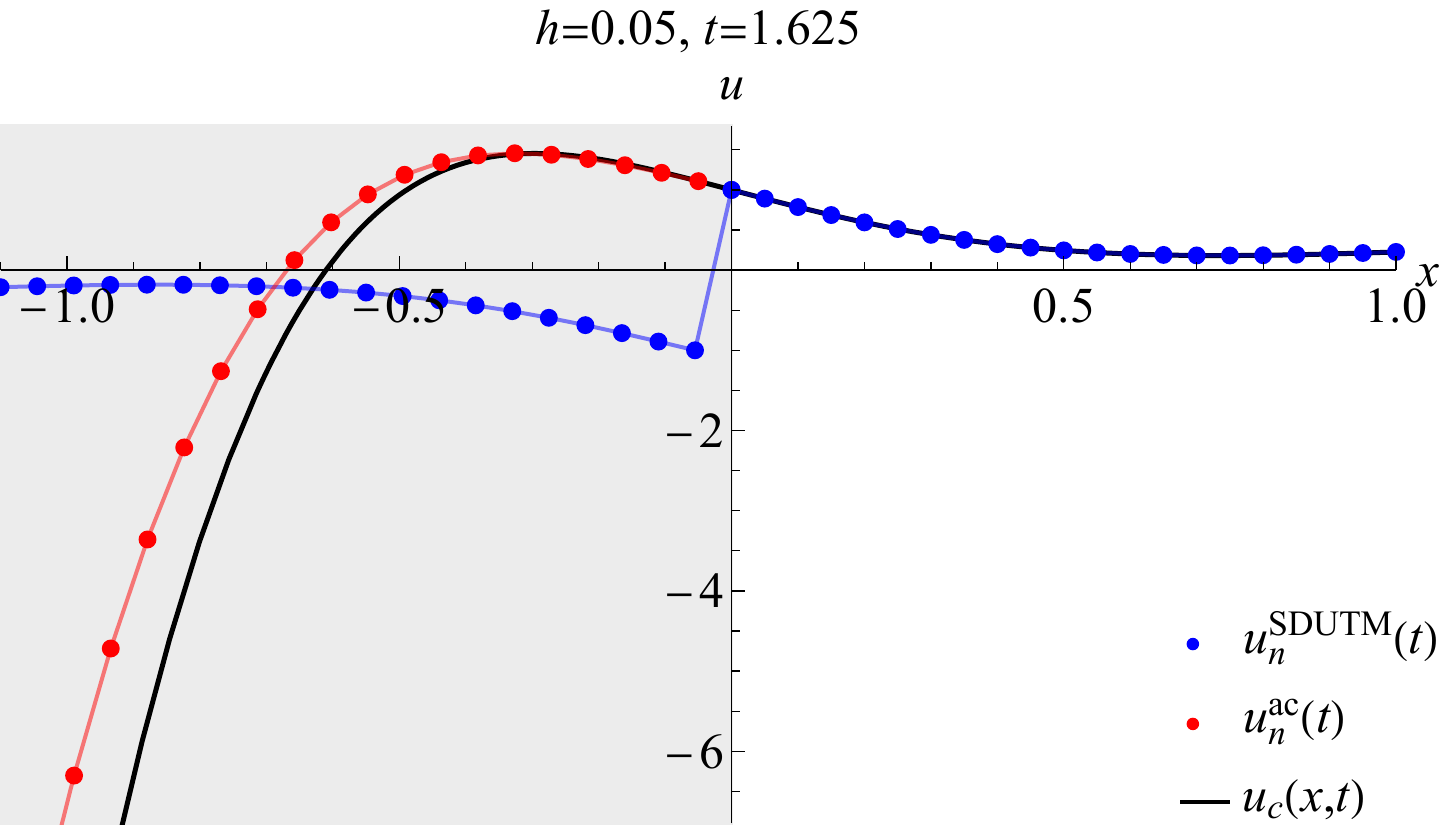} & \includegraphics[scale=\sc]{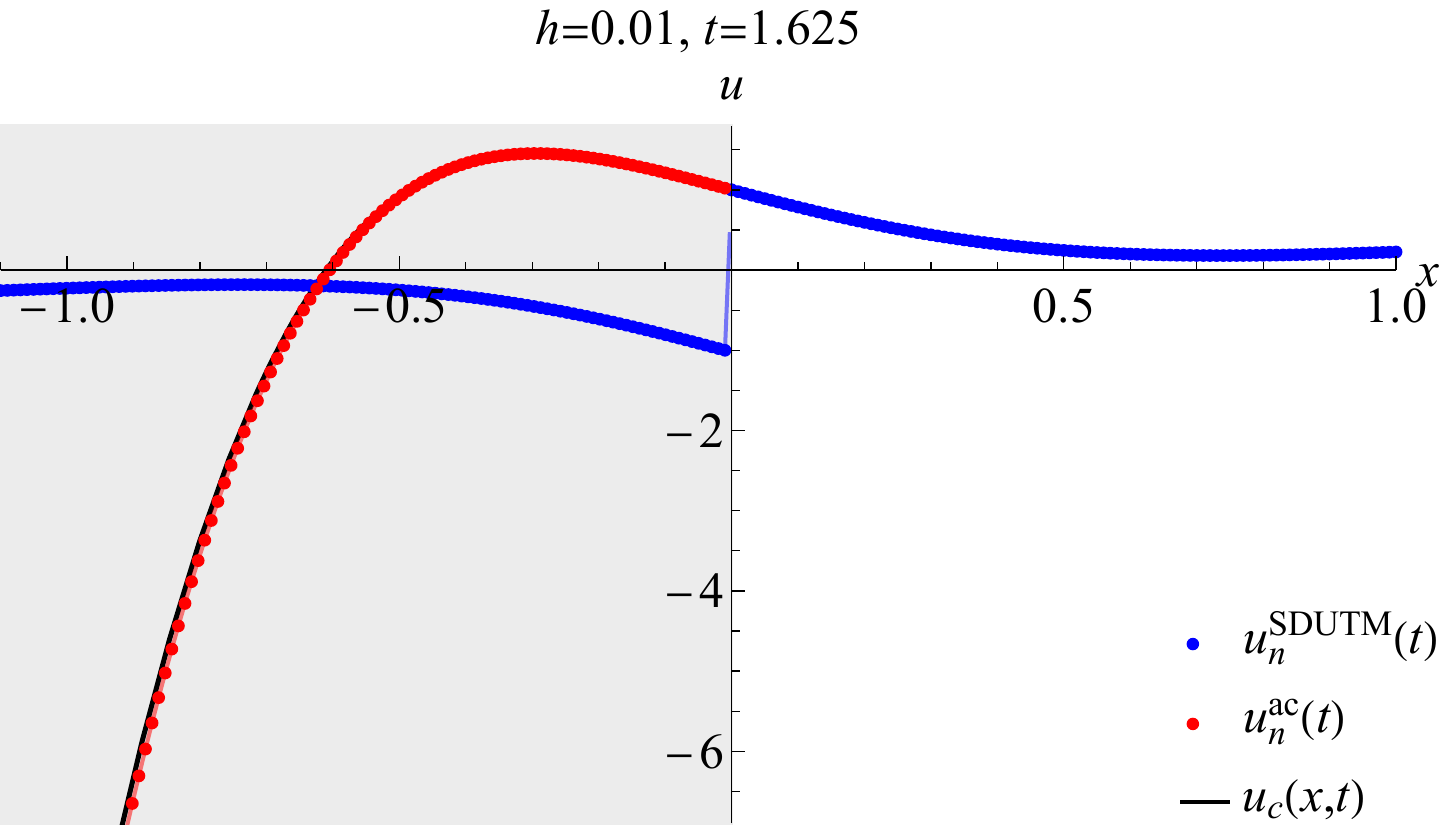}\\
(a) & (b)
\end{tabular}
\end{center}
\caption{(a) For $h=1/20$, the solution $u_n^{\text{ac}}(t)$ given by \rf{soln_heat_centered_negn} at $t=0.5$ obtained through analytic continuation, shown with the SDUTM solution $u_n^\text{SDUTM}(t)$ \rf{soln_heat_centered} and the continuum solution $u_c(x,t)$ given by \rf{eqn:soln_heat}. (b) The same plot for $h=1/100$.
}
\label{heat_centered_analytcont1}
\end{figure}

\no for all $n \in \mathbb{Z}$. The results are presented in Figure~\ref{heat_centered_analytcont1}.


\section{Conclusion}

We have demonstrated how the solutions of linear, constant-coefficient IBVPs can be analytically extended outside their original spatial domain of definition, using the UTM as a method for doing so. A number of representative examples were used to illustrate our approach. The results most useful for computational purposes for these examples are collected here. The reader should refer to the sections where these examples are treated in detail for additional results, and for the introduction of the notation used below. 

\begin{itemize}
    
\item {\bf The heat equation on $x>0$ with Dirichlet boundary conditions}, see \rf{heatdirichletibvp}. The solution for $x\in \mathbb{R}$ may be represented as 

\begin{equation} \label{eqn:soln_heat}
u(x,t)=I_0(x,t) + I_{f_0}^{\text{ext}}(x,t), 
\end{equation}

where $I_0(x,t)$ and $I_{f_0}^{\text{ext}}(x,t)$ are defined in \rf{I0heat1} and \rf{heat_cont_full}, respectively. 

\item {\bf The heat equation on $x>0$ with Neumann boundary conditions}, see \rf{heatneumannibvp}. The solution for $x\in \mathbb{R}$ may be represented as

\begin{equation}
u(x,t)=I_0(x,t) + I_{f_1}^{\text{ext}}(x,t), 
\end{equation}

where $I_0(x,t)$ and $I_{f_1}^{\text{ext}}(x,t)$ are defined in \rf{heatneumanni0} and \rf{eqn:If1_ext}, respectively.

\item {\bf The heat equation on $x\in (0,L)$ with Dirichlet boundary conditions}, see \rf{heatfiibvp}. The solution for $x\in \mathbb{R}$ may be represented as

\begin{equation}
u(x,t) = I_0(x,t) + I_{f_0}^{\text{ext}}(x,t) + I_{g_0}^{\text{ext}}(x,t),
\end{equation}
where $I_0(x,t)$, $I_{f_0}^{\text{ext}}(x,t)$ and $I_{g_0}^{\text{ext}}(x,t)$ are defined in \rf{heatfii0}, \rf{heatif0ext} and \rf{heatig0ext}, respectively. 

\item {\bf The advected heat equation on $x>0$ with Dirichlet boundary conditions}, see \rf{advheatibvp}. The solution for $x\in \mathbb{R}$ may be represented as

\begin{equation}
u(x,t) = I_0(x,t) + I_{f_0}^{\text{ext}}(x,t),
\end{equation}

where $I_0(x,t)$ and $I_{f_0}^{\text{ext}}(x,t)$ are defined in \rf{I0heatadv1} and \rf{advheatbext}, respectively.

\item {\bf The linear KdV equation $u_t+u_{xxx}=0$ on $x>0$ with Dirichlet boundary conditions}, see \rf{kdvprob1}. The solution for $x\in \mathbb{R}$ may be represented as

\begin{equation}
u(x,t) = I_0(x,t) + I_{f_0}^{\text{ext}}(x,t),
\end{equation}

where $I_0(x,t)$ and $I_{f_0}^{\text{ext}}(x,t)$ are defined in \rf{kdv1I0} and \rf{advheatIf0ext}, respectively.

\item {\bf The linear KdV equation $u_t-u_{xxx}=0$ on $x>0$ with Dirichlet boundary conditions}, see \rf{kdvprob2}. The solution for $x\in \mathbb{R}$ may be represented as

\begin{equation}
u(x,t) = I_0(x,t) + I_{f_0}^{\text{ext}}(x,t)+ I_{f_1}^{\text{ext}}(x,t),
\end{equation}

where $I_0(x,t)$ is defined in \rf{kdv2I0}, $I_{f_0}^{\text{ext}}(x,t)$ and $I_{f_1}^{\text{ext}}(x,t)$ are defined in \rf{eqn:kdv2ext0} and \rf{eqn:kdv2ext1}, respectively.

\item {\bf The backward-discretized advection equation on $n>0$ with Dirichlet boundary conditions}, see \rf{advecdiscrete1}. 
The solution for $n\in \mathbb{Z}$ may be represented by \rf{soln_advec2_backward_nok_alln}. 

\item {\bf The second-order backward-discretized advection equation on $n>0$ with Dirichlet boundary conditions}, see \rf{advec2_backward2}. The solution for $n\in \mathbb{Z}$ may be represented by \rf{soln_advec2_backward2} for $n>0$, and by \rf{soln_AC_advec2_backward2} for $n<0$. 

\item {\bf The second-order centered-discretized heat equation on $n>0$ with Dirichlet boundary conditions}, see \rf{heatdiscretecentered}. The solution for $n\in \mathbb{Z}$ may be represented by \rf{soln_heat_centered} for $n>0$, and by \rf{soln_heat_centered_negn} for $n<0$. 

\end{itemize}

\begin{appendix}

\section*{Appendix. The second-order centered discretized advection equation}
\label{sec:appendix}

\label{advec2_centered_halfline_analytcont}

This section demonstrates 
the discretization of \eqref{advec2_prob} using the standard centered stencil, which is known to be a poor choice \cite{JC_randy}. The SDUTM may still be used to solve the discretized system and the analytic continuation method still gives a continuation formula. However, since the method is entirely dispersive, it should not be surprising that undesirable dispersive behavior occurs.

Discretizing \eqref{advec2_prob} using the standard centered stencil gives

\beq \label{eqn:cent_disc} \dot{u}_n(t) = c\,\frac{u_{n-1}(t) - u_{n+1}(t)}{2h},\eeq
\no with dispersion relation 
\beq W(k) = c\,\frac{e^{ikh} - e^{-ikh}}{2h},\eeq
\no and nontrivial symmetry $\nu_1(k) = - k - {\pi}/{h}$. The semi-discrete solution is
\beq u_n(T) = \frac{1}{2\pi} \int_{-\pi/h}^{\pi/h} e^{iknh} e^{-WT} \left[ \hat{u}(k,0) -\hat{u}\left(\nu_1,0\right)\right] \,dk + \frac{c}{2\pi} \int_{-\pi/h}^{\pi/h} e^{iknh} e^{-WT} \cos\left( kh \right)F_0(W,T) \,dk. 	\label{soln_advec2_centered} \eeq
\no Using this representation, it follows that $u_{-n}(T) = (-1)^{n+1}u_n(T)$ for $n \in \mathbb{Z}^+$. Repeating similar steps as before, the correct analytic extension satisfies
\beq \label{eqn:cent_ac} u_{-n}(T) = (-1)^{n+1} u_n(T) \,\,\,+\,\,\,n \, \sum_{\substack{k = 0\\n-k\text{ even}}}^{n} f_0^{(k)}(T) \,\left(\frac{2h}{c}\right)^{k} \,\frac{\Gamma\left(\frac{n+k}{2}\right)}{k!\, \left(\frac{n-k}{2}\right)!}, \qquad n>0.\eeq

Due to the purely dispersive nature of the discretization \rf{eqn:cent_disc}, if the boundary function $f_0(t)$ and the initial condition $\phi(x)$ are not compatible (\textit{i.e.}, $f_0(t)\neq\phi(-ct)$), this solution exhibits a dispersive fan, see Figure~\ref{fig:advcentered}{\textcolor{blue}a} and \ref{fig:advcentered2}{\textcolor{blue}a}. The dispersive fan seen in this figure travels upwind at speed $c$. Its envelope is determined by a combination of the initial and the boundary conditions, see \eqref{eqn:envelope}. 
For $c<0$, \rf{advec2_prob} does not have a boundary condition, but the discrete solution \rf{soln_advec2_centered} requires one. If we consider the discretized system in its own right and apply an incompatible boundary condition, we observe the dispersive fan for $n>0$, demonstrating the appearance of the dispersive fan is not an analytic continuation issue, but one of the discretization itself. If we choose the compatible boundary condition $f_0(t) = \phi(-ct)$ with analytic $\phi(x)$, we do not see the dispersive fan, see Figure~\ref{fig:advcentered2}{\textcolor{blue}b}. In general, discontinuous 
initial data $\phi(x)$ 
gives rise to dispersive effects (as in \cite{trogdon}, see Figure~\ref{fig:advcenteredsquare}), further demonstrating that \rf{eqn:cent_disc}
is an inappropriate discretization  \cite{JC_randy} and that the dispersive fan in Figure~\ref{fig:advcentered}{\textcolor{blue}a} and \ref{fig:advcentered2}{\textcolor{blue}a} is not caused by discontinuity or the analytic continuation formula \rf{eqn:cent_ac}, but by the incompatibility of the boundary and initial data in conjunction with the dispersive nature of \rf{eqn:cent_disc}. 

It is interesting to note that one can ``average out'' the dispersive fan wave to get the extension formula
\beq \label{eqn:cent_ac_avg} u_{-n}(T) = \frac n2 \, \sum_{k = 0}^{n} f_0^{(k)}(T) \,\left(\frac{2h}{c}\right)^{k} \,\frac{\Gamma\left(\frac{n+k}{2}\right)}{k!\, \left(\frac{n-k}{2}\right)!}, \eeq
\no which is second-order accurate, but does not solve the discretized equation \rf{eqn:cent_disc}. Considering \eqref{eqn:cent_ac} in its continuum limit,
\begin{align}\nonumber \label{eqn:envelope}
u_{-n}(T) &\sim (-1)^{n+1} u_n(T) \,\,\,+\,\,\, 2 \, \sum_{\substack{k = 0\\n-k\text{ even}}}^{n} \frac{f_0^{(k)}(T)}{k!} \,\left(-\frac{x}{c}\right)^{k} \\
&\sim (-1)^{n+1} u_n(T) \,\,\,+\,\,\, f_0\left(T-\frac{x}{c}\right) +(-1)^n f_0\left(T+\frac{x}{c}\right),
\end{align}

\no and since for $x<0$,

\beq 
u_n(T) \to \case{1}{f_0(T+x/c), & x>-cT,\\ \phi(-x-ct), & x<-cT, }
\eeq

\no which has a limit as $n\to \infty$ ($h\to0^+$ and $nh\to -x$) if and only if $u_n(T)=f_0\left(T+{x}/{c}\right)$, \textit{i.e.}, if the initial and boundary conditions are compatible. Further, these previous two equations give the envelope of the dispersive fan.

\begin{figure}[tb]
\begin{center}
\def \sc {0.475}
\begin{tabular}{cc} 
\includegraphics[scale=\sc]{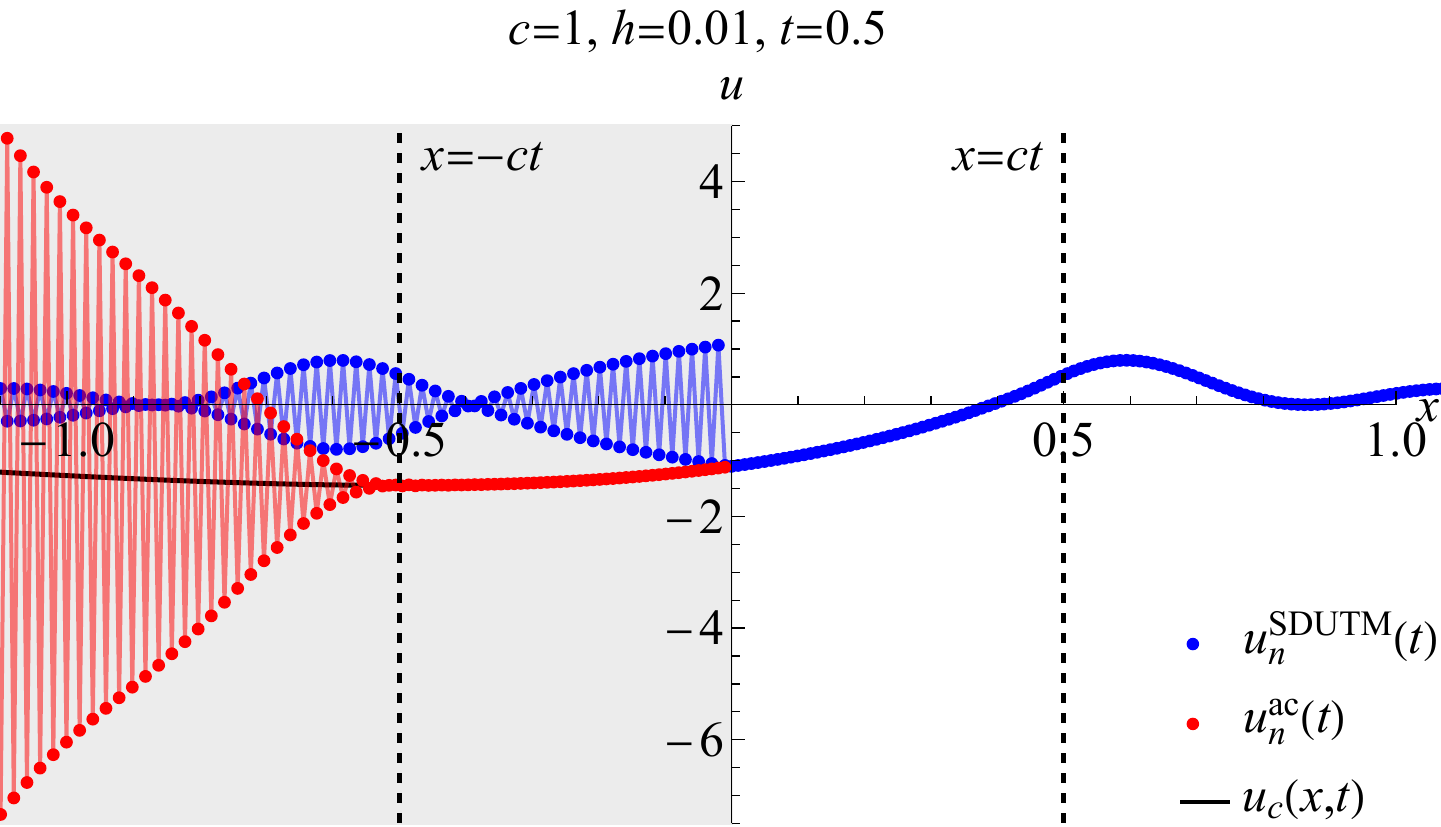} &\includegraphics[scale=\sc]{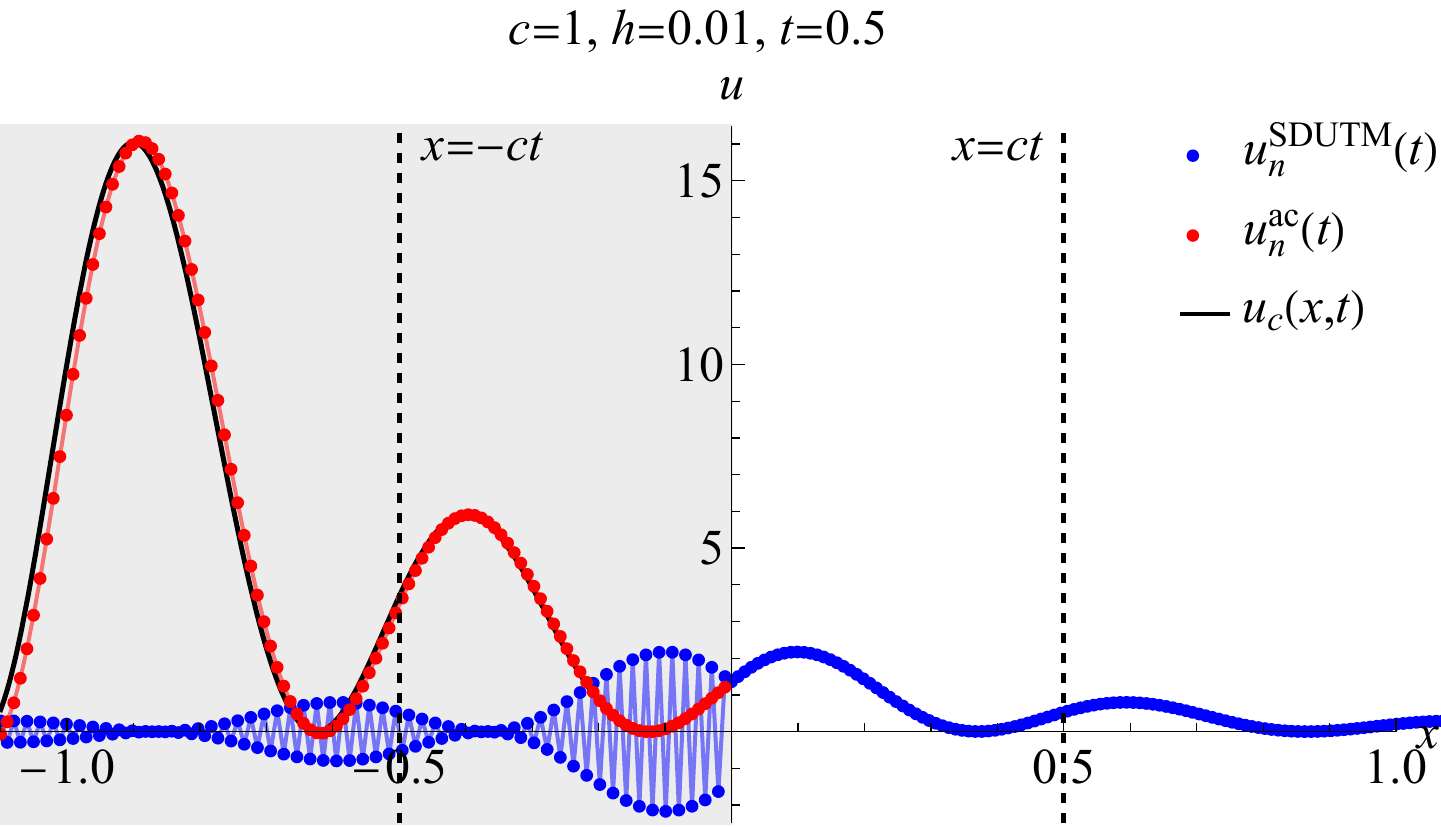}\\
(a) & (b)
\end{tabular}
\end{center}
\caption{For $c=1$: (a) The SDUTM solution \rf{soln_advec2_centered} for incompatible boundary and initial conditions given by \rf{advec2_prob_ex}, shown with the analytic continuation given by \rf{eqn:cent_ac} and the continuum solution $u_c(x,t)$ given by \rf{eqn:soln_adv}. (b) The SDUTM solution, \rf{soln_advec2_centered}, for initial condition \rf{advec2_prob_ex} and compatible boundary conditions, shown with the analytic continuation given by \rf{eqn:cent_ac} and the continuum solution $u_c(x,t)$ given by \rf{eqn:soln_adv}.}
\la{fig:advcentered}
\end{figure}

\begin{figure}[tb]
\begin{center}
\def \sc {0.475}
\begin{tabular}{cc} 
\includegraphics[scale=\sc]{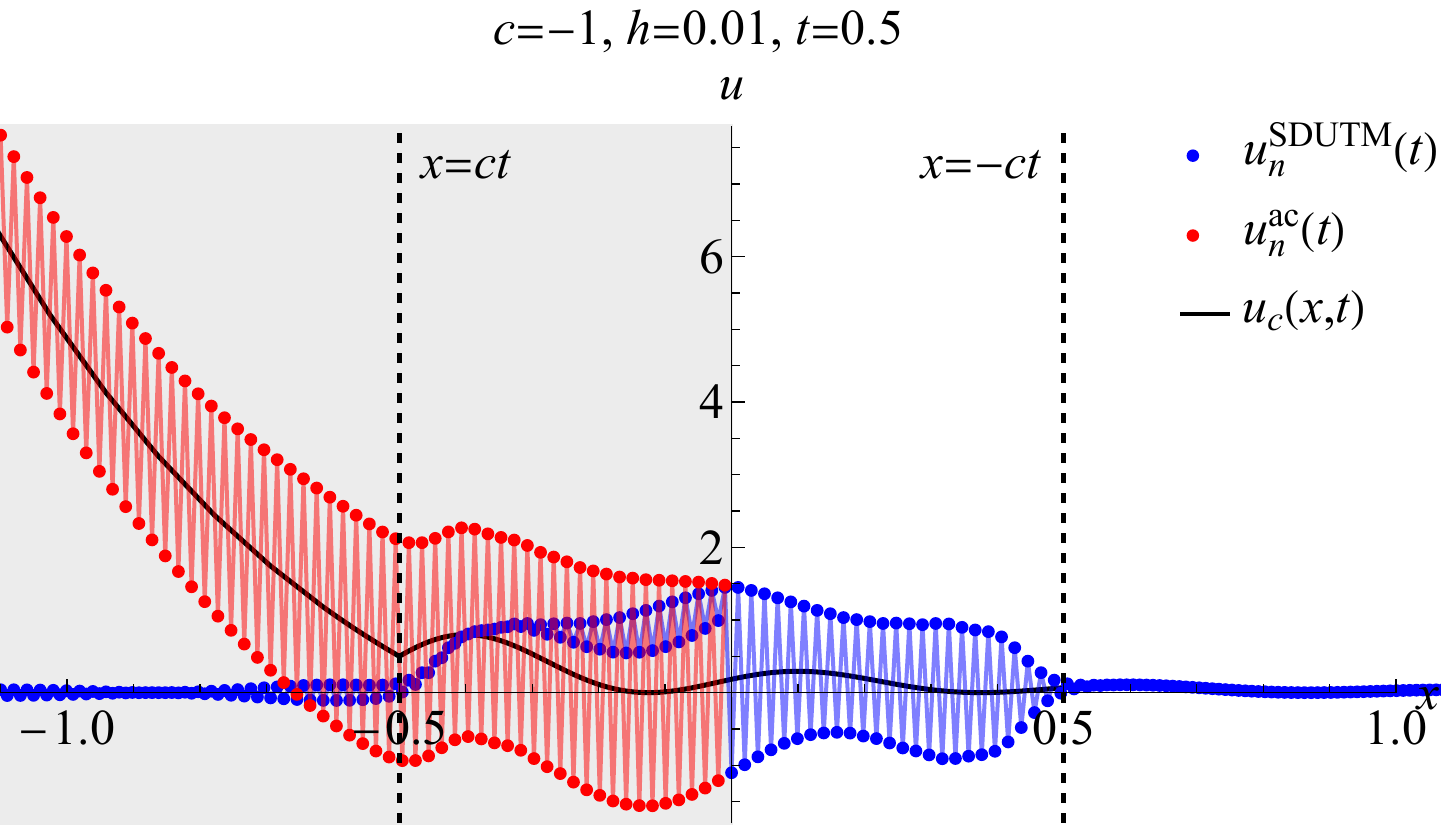} 
& \includegraphics[scale=\sc]{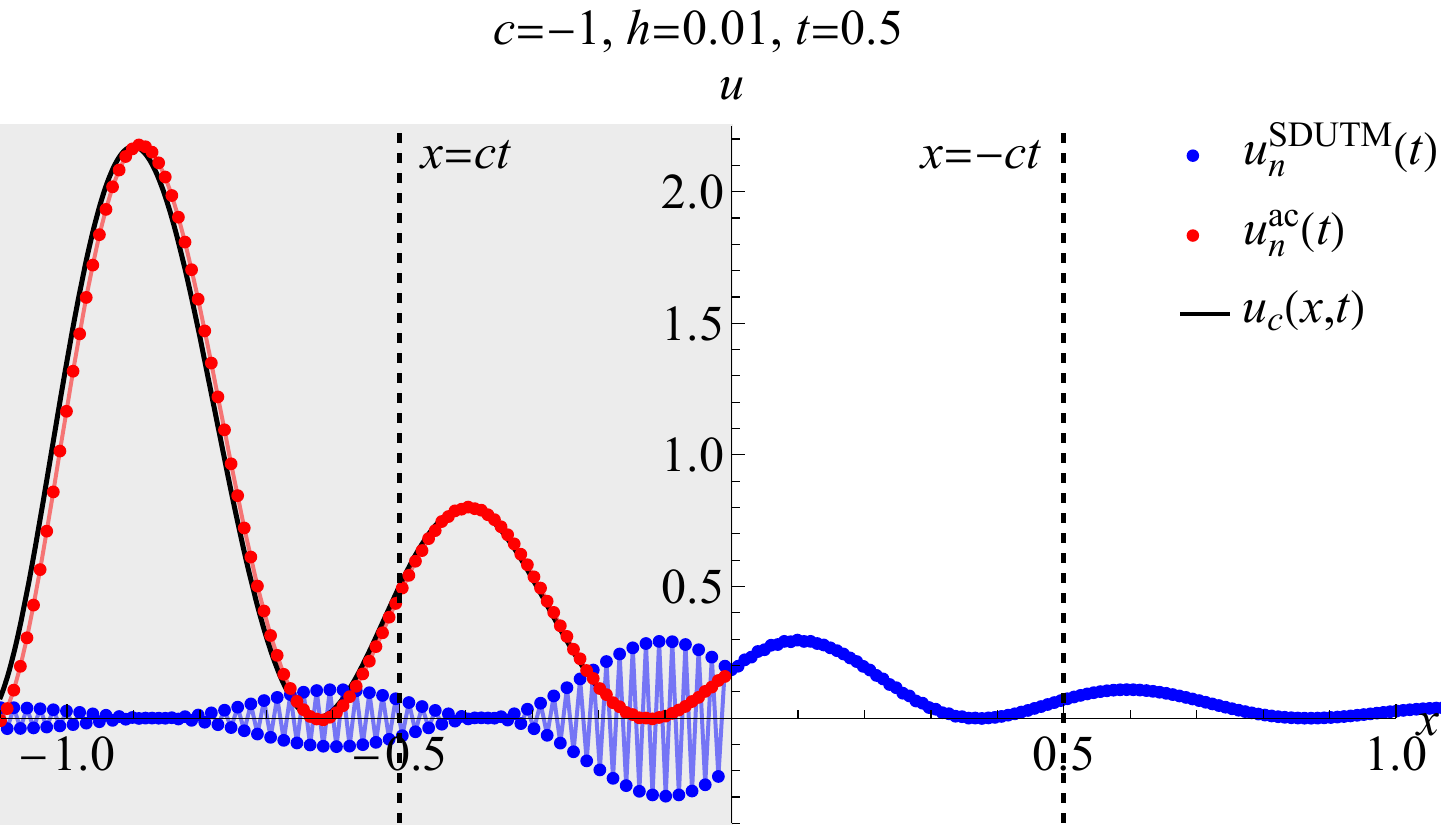} \\
(a) & (b)
\end{tabular}
\end{center}
\caption{For $c=-1$: (a) The SDUTM solution \rf{soln_advec2_centered} for incompatible boundary and initial conditions given by \rf{advec2_prob_ex}, shown with the analytic continuation given by \rf{eqn:cent_ac} and the continuum solution $u_c(x,t)$ given by \rf{eqn:soln_adv}. (b) The SDUTM solution, \rf{soln_advec2_centered}, for initial condition \rf{advec2_prob_ex} and compatible boundary conditions, shown with the analytic continuation given by \rf{eqn:cent_ac} and the continuum solution $u_c(x,t)$ given by \rf{eqn:soln_adv}.}
\la{fig:advcentered2}
\end{figure}

\begin{figure}[tb]
\begin{center}
\def \sc {0.46}
\begin{tabular}{cc} 
\includegraphics[scale=\sc]{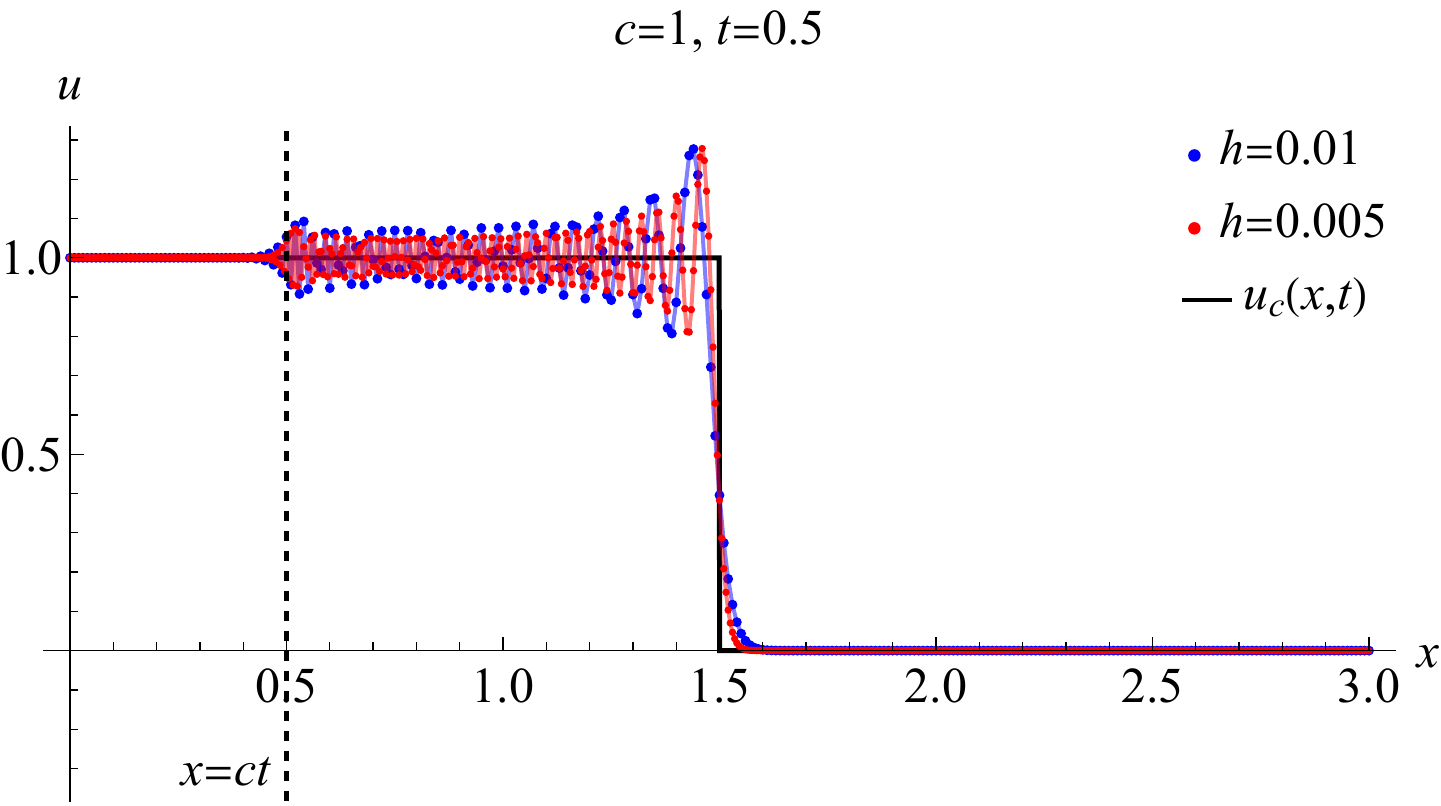} & \hspace{0.1in} \includegraphics[scale=\sc]{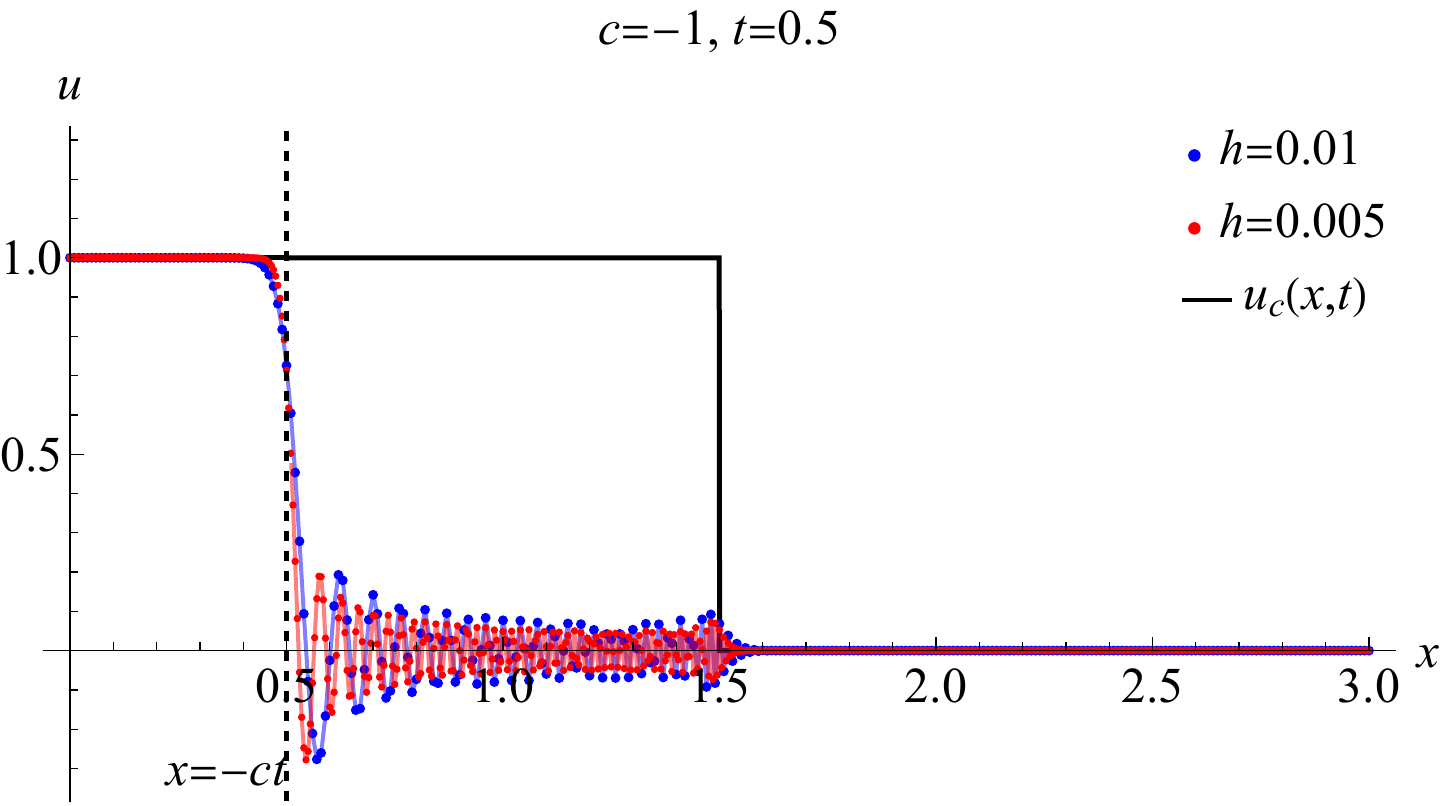}\\
(a) & (b)
\end{tabular}
\end{center}
\caption{The centered discretization with step function initial data, resulting in Gibbs phenomenon-like behavior for (a) $c=1$, and (b) $c=-1$.}
\la{fig:advcenteredsquare}
\end{figure}

\end{appendix}

\section*{Acknowledgements} The authors thank Tom Trogdon for useful conversations. This work was supported by the Graduate Opportunities \& Minority Achievement Program Fellowship
from the University of Washington and the Ford Foundation Predoctoral Fellowship (JC). Any opinions, findings, and conclusions or recommendations expressed in this material are those of the authors and do not necessarily reflect the views of the funding sources.


{\small
\bibliographystyle{abbrv}
\bibliography{references}
}

\end{document}